\documentclass[hidelinks, 10pt]{article}
\usepackage[UKenglish]{babel}
\usepackage[utf8]{inputenc}
\usepackage{amsmath}
\usepackage{amsthm}
\usepackage{amsfonts}
\usepackage{amssymb}
\usepackage{mathtools}
\usepackage{dsfont} % for identity
\usepackage[dvipsnames]{xcolor}
\usepackage[breaklinks]{hyperref}
\usepackage{lmodern}
\usepackage{fullpage}
\usepackage{microtype}
\usepackage{newpxmath,newpxtext}
\usepackage{booktabs}
\usepackage{enumitem} % to customize description lists
\usepackage{comment}
\usepackage{tikz}
\usetikzlibrary{positioning}
\usetikzlibrary{arrows,backgrounds}
\usetikzlibrary{decorations.pathmorphing}
\usepackage{tikz-cd}
\usepackage{adjustbox}
\usepackage[autostyle,italian=guillemets]{csquotes}
\usepackage[style=alphabetic,backend=biber]{biblatex}
\usepackage{url}

\usepackage{breakurl}
\tikzset{node distance=2cm, auto}
\tikzcdset{row sep/normal=2.7em,column sep/normal=3.5em}
\allowdisplaybreaks
\setlist[description]{font=\normalfont\bfseries\space}
\renewcommand\labelenumi{(\roman{enumi})}

\renewcommand\theenumi\labelenumi
\hypersetup{
    colorlinks=true,
    linkcolor=gray,
    filecolor=gray,      
    urlcolor=gray,
    citecolor=gray
    }
\numberwithin{equation}{section}
\newtheorem{theorem}{Theorem}[section]
\newtheorem{proposition}[theorem]{Proposition}
\newtheorem{lemma}[theorem]{Lemma}
\newtheorem{corollary}[theorem]{Corollary}

\theoremstyle{definition}
\newtheorem{definition}[theorem]{Definition}
\newtheorem{example}[theorem]{Example}
\newtheorem{notation}[theorem]{Notation}
\newtheorem{remark}[theorem]{Remark}
\newenvironment{acknowledgements}{\vskip .5cm\noindent\small\textbf{Acknowledgements}.\quad\noindent}{}

\newcommand{\Font}[1]{\mathsf{#1}}
\renewcommand{\d}{\mathrm{d}}
\renewcommand{\,}{{,\dots,}}
\renewcommand{\-}{{-1}}
\newcommand{\x}{\otimes}
\newcommand{\1}{\mathds{1}}
\newcommand{\CC}{{\mathbb{K}}}
\renewcommand{\=}{\colon\kern-1ex=}
\renewcommand{\epsilon}{\varepsilon}
\renewcommand{\o}{\circ}
\newcommand{\s}{{\scriptscriptstyle\square}}
\renewcommand{\b}{{\scriptscriptstyle\blacksquare}}
\newcommand{\<}{\langle}
\renewcommand{\>}{\rangle}
\newcommand{\op}{\Font{op}}
\newcommand{\co}{\Font{co}}
\newcommand{\lnr}{\Font{lnr}}
\newcommand{\strb}{\Font{strb}}
\newcommand{\pp}{\Font{pp}}
\renewcommand{\^}[1]{^{\left(#1\right)}}
\renewcommand{\*}{\ast}
\newcommand{\nto}{\nrightarrow}

\newcommand{\Comp}{\I^2}
\newcommand{\Unit}{\I^0}
\newcommand{\El}{\Font{EL}}

\newcommand{\mathdash}{\relbar\mkern-3mu\relbar\mkern-3mu\relbar---}
\newcommand{\N}{\mathbb{N}}
\newcommand{\Z}{\mathbb{Z}}

\newcommand{\id}{\Font{id}}
\newcommand{\End}{\Font{END}}
\newcommand{\cRig}{\Font{CRIG}}
\newcommand{\cRing}{\Font{CRING}}
\newcommand{\X}{\mathbb{X}}
\newcommand{\Y}{\mathbb{Y}}
\newcommand{\Mod}{\Font{MOD}}
\newcommand{\Cat}{\Font{CAT}}

\newcommand{\Mnd}{\Font{MND}}

\newcommand{\MonCat}{\Font{MONCAT}}
\newcommand{\Enrch}{\Font{ENRCH}}
\newcommand{\twoCat}{2\Cat}
\newcommand{\Fib}{\Font{FIB}}
\newcommand{\Indx}{\Font{INDX}}
\newcommand{\trans}{\Font{trans}}
\newcommand{\FF}{\mathscr{F}}
\newcommand{\I}{\mathfrak{I}}
\newcommand{\II}{\mathscr{I}}
\newcommand{\Sym}{\mathrm{Sym}}
\newcommand{\T}{\mathrm{T}}
\newcommand{\TT}{\mathbb{T}}
\newcommand{\TngCat}{\Font{TNGCAT}}
\newcommand{\TTngCat}{\mathbb{T}\Font{NGCAT}}
\newcommand{\Tng}{\Font{TNG}}
\newcommand{\TTng}{\mathbb{T}\Font{NG}}
\newcommand{\DBnd}{\Font{DBND}}

\newcommand{\TngFib}{\Font{TNGFIB}}
\newcommand{\IndxTng}{\Font{INDXTNG}}
\newcommand{\TngIndx}{{\Font{TNGINDX}}}
\newcommand{\Dsply}{\mathscr{D}}
\newcommand{\Weil}{\Font{WEIL}_1}
\newcommand{\Leung}{\mathrm L}
\title{The Grothendieck construction\\ in the context of tangent categories}
\author{Marcello Lanfranchi}
\date{}
\makeatother
\makeatletter
\begin{document}

\maketitle

\begin{abstract}\noindent
The Grothendieck construction establishes an equivalence between fibrations, a.k.a. fibred categories, and indexed categories, and is one of the fundamental results of category theory. Cockett and Cruttwell introduced the notion of fibrations into the context of tangent categories and proved that the fibres of a tangent fibration inherit a tangent structure from the total tangent category. The main goal of this paper is to provide a Grothendieck construction for tangent fibrations. Our first attempt will focus on providing a correspondence between tangent fibrations and indexed tangent categories, which are collections of tangent categories and tangent morphisms indexed by the objects and morphisms of a base tangent category. We will show that this construction inverts Cockett and Cruttwell's result but it does not provide a full equivalence between these two concepts. In order to understand how to define a genuine Grothendieck equivalence in the context of tangent categories, inspired by Street's formal approach to monad theory we introduce a new concept: tangent objects. We show that tangent fibrations arise as tangent objects of a suitable $2$-category and we employ this characterization to lift the Grothendieck construction between fibrations and indexed categories to a genuine Grothendieck equivalence between tangent fibrations and tangent indexed categories.
\end{abstract}

\begin{acknowledgements}
We would like to thank Geoffrey Cruttwell for introducing us to the world of fibrations. We also would like to thank Dorette Pronk and Geoffrey Vooys for the useful discussions on this topic and for the great suggestions for future work. We would like to thank again Geoffrey Cruttwell and Dorette Pronk for their support, suggestions, and corrections, and Rory Lucyshyn-Wright for an important correction in the definition of tangent objects (Remark~\ref{remark:pointwise-limits}). We are also grateful to the anonymous referees for their suggestions.
\end{acknowledgements}

%__________________________________________________________________________
%__________________________________________________________________________

\tableofcontents

\section{Introduction}
\label{section:introduction}
In topology, the Fibre Bundle Construction Theorem establishes an equivalence between fibre bundles on a given base space and a collection of fibres, equipped with a suitable family of transition functions (cf.~\cite[Section~3]{steenrod:topology-fibre-bundles}). The category-theoretic analog of the Fibre Bundle Construction is captured by the Grothendieck construction, which establishes an equivalence between (cloven) fibrations, also known as fibred categories, and indexed categories. Informally, the latter can be interpreted as a collection of categories indexed by the objects of a fixed category, together with a collection of functors between them, indexed by the morphisms of the index category, in a compatible way with the indexing of the objects, while the former is the category-theoretic analog of fibre bundles.
\par More concretely, a fibration consists of a functor from a category, known as the total category of the fibration, to another category, known as the base category of the fibration, satisfying a universal property for which morphisms of the base category can be lifted to the total category suitably. This universal property of fibrations is in every aspect the category-theoretic analog of the path-lifting property of the universal covering space of a topological space.
\par Fibrations are suitable to describe objects that depend on another object of another category. An example of a fibration is given by the category of pairs $(R,M)$ formed by a ring $R$ and an $R$-module $M$, and morphisms $(f,\xi)\colon(R,M)\to(S,N)$ defined by a morphism $f\colon R\to S$ of rings and a morphism of $R$-modules $\xi\colon M\to f^\*N$, where $f^\*N$ denotes the $R$-module induced by $f\colon R\to S$ via scalar inclusion. The fibration corresponds to the projection which sends each pair $(R,M)$ to the underlying ring $R$ and each morphism $(f,\xi)$ to the underlying morphism of rings $f$.
\par Equivalently, one can display the same information by defining a pseudofunctor which sends each ring $R$ to the category $\Mod_R$ of $R$-modules and each morphism of rings $f\colon R\to S$ to the functor $f^\*\colon\Mod_S\to\Mod_R$. The Grothendieck construction establishes that every (cloven) fibration is equivalent to a fibration whose total category is the category of pairs $(A,E)$ formed by an object $A$ of the base category together with an object $E$ on the fibre over $A$.
\par Cockett and Cruttwell in~\cite[Section~5]{cockett:differential-bundles}, introduced the notion of fibrations into the context of tangent categories. Tangent categories, first defined by Rosick\`y in~\cite{rosicky:tangent-cats} and recently revisited and generalized by Cockett and Cruttwell in~\cite{cockett:tangent-cats}, provide a minimal categorical framework to axiomatize the tangent bundle functor of differential geometry.
\par Concretely, a tangent category consists of a category $\X$ equipped with an endofunctor $\T$, called the tangent bundle functor, together with a collection of natural transformations: a projection $p\colon\T\Rightarrow\id_\X$ which makes $\T A$ into a bundle for every object $A$ of $\X$, a zero section $z\colon\id_\X\Rightarrow\T$ of the projection, known as the zero morphism, a sum morphism $s\colon\T_2\Rightarrow\T$, whose domain is the pullback of the projection along itself, that allows one to sum tangent vectors of the same fibre, a vertical lift $l\colon\T\Rightarrow\T^2$ which introduces a notion of local linearity, and a canonical flip $c\colon\T^2\Rightarrow\T^2$, which axiomatizes the commutativity of partial derivatives for smooth functions.
\par In their paper~\cite{cockett:differential-bundles}, Cockett and Cruttwell proved an important result: the fibres of a tangent fibration, i.e., a fibration between tangent categories which preserves the tangent structures, inherit a tangent structure from the total tangent category, strongly preserved by the substitution functors (cf.~\cite[Theorem~5.3]{cockett:differential-bundles}). One would like to see this result as part of a correspondence between tangent fibrations and indexed categories whose fibres are equipped with a tangent structure and whose substitution functors are compatible with these tangent structures, which will be called indexed tangent categories.
\par The first goal of this paper is to explore this idea, by showing that the Grothendieck construction of such indexed categories gives rise to a tangent fibration. We refer to this as the \textit{reduced Grothendieck construction for tangent categories}. This, together with Cockett and Cruttwell's construction, defines an adjunction which, however, does not provide an equivalence of categories.
\par Even if this ``fibrewise'' Grothendieck construction does not provide a full equivalence, it is still very relevant in tangent category theory. In particular, Cockett and Cruttwell have shown how the construction of the slice tangent category can be obtained by looking at the tangent structures on the fibres of a given tangent fibration. Furthermore, this construction plays a crucial role in the theory of differential bundles (cf.~\cite{cockett:differential-bundles}).
\par The main goal of this paper is to achieve a genuine Grothendieck equivalence in the context of tangent categories and clarify why the reduced Grothendieck construction fails in this manner. Our approach relies on a new concept: the notion of a \textit{tangent object} of a $2$-category. The idea is to introduce a formal approach for tangent category theory, inspired by Street's formal approach to monad theory~\cite{street:formal-theory-monads}. After introducing this definition, we show that tangent categories are precisely tangent objects in the $2$-category of categories. Similarly, we show that tangent fibrations are precisely tangent objects in the $2$-category of fibrations. We take advantage of this fact to lift the equivalence between (cloven) fibrations and indexed categories to a genuine Grothendieck equivalence between (cloven) tangent fibrations and tangent indexed categories.
\par The desire to explore a Grothendieck construction between certain fibrations and indexed categories equipped with extra structure is not new. For instance, Moeller and Vasilakopoulou in~\cite{moeller:monoidal-grothendieck} investigated a Grothendieck construction in the context of monoidal categories. They encountered a similar phenomenon: a monoidal structure on a fibration can be either fibrewise or global, i.e., such that the total category becomes monoidal. In the general case, these two structures are not equivalent, however, under certain assumptions they become compatible.
\par As part of our investigation, we also compare the notion of internal fibrations introduced by Street in~\cite{street:fibrations} with Cockett and Cruttwell's notion of tangent fibrations. In particular, we show that tangent fibrations can be described as a suitable subclass of pseudoalgebras of a $2$-monad.

%__________________________________________________________________________
\subsection{Outline}
\label{subsection:outline}
In Section~\ref{section:background}, we recall the notions of a fibration and an indexed category, and we extend the classical Grothendieck construction to an equivalence between fibrations and indexed categories over a non-fixed base. In Section~\ref{subsection:tangent-cats-background}, we recall the main definitions of tangent category theory.
\par Section~\ref{section:tangent-fibrations} is dedicated to recalling Cockett and Cruttwell's notion of a tangent fibration and to reformulating their characterization of the fibres of a tangent fibration in terms of a functor which sends a tangent fibration to an indexed tangent category.
\par In Section~\ref{subsection:reduced-Grothendieck}, we show that this functor admits a left adjoint, which however does not provide an equivalence of categories. We conclude Section~\ref{section:tangent-fibrations} with a comparison between internal fibrations in the $2$-category of tangent categories and tangent fibrations.
\par In order to construct a genuine Grothendieck construction for tangent fibrations, in Section~\ref{section:tangent-objects} we introduce a new concept: the notion of a tangent object. In Section~\ref{section:full-Grothendieck}, we employ this formal approach to prove that tangent fibrations are equivalent to tangent indexed categories, i.e., tangent objects in the $2$-category of indexed categories.
\par Finally, Section~\ref{section:conclusions} is dedicated to summarize the story of this paper and to discuss some directions for future work.

%__________________________________________________________________________
\subsection{Notation}
\label{subsection:notation}
Categories are denoted by capitalized names, e.g., $\Cat$ for the category of (small) categories. To denote identity morphisms we adopt the notation $\id_A\colon A\to A$ and for the composition of morphisms, functors, and natural transformations we adopt the functional convention, e.g., $f\o g\colon A\to C$, for $f\colon B\to C$ and $g\colon A\to B$. However, when a composition of functors $F\o G$ is evaluated on an object $A$ or on a morphism $f$, we adopt the applicative notation, e.g., $FGA$ or $FGf$ to denote $(F\o G)(A)$ or $(F\o G)(f)$, respectively.

%__________________________________________________________________________
%__________________________________________________________________________

\section{Background}
\label{section:background}
In this section, we recall some of the main results of fibration theory and tangent category theory. In particular, in Section~\ref{subsection:fibrations-background} we recollect the notion of fibrations, indexed categories, and the Grothendieck construction, while in Section~\ref{subsection:tangent-cats-background} we recall some of the main concepts and results of tangent category theory. For fibration theory we refer to~\cite{borceux:fibrations} and~\cite{streicher:fibrations}; for tangent category theory we refer to~\cite{cockett:tangent-cats}, and~\cite{cockett:differential-bundles}.

%__________________________________________________________________________
\subsection{Fibrations and Grothendieck construction}
\label{subsection:fibrations-background}
A fibration, introduced by Grothendieck (see~\cite[Expos\'e~VI.12 Section~12]{grothendieck:fibrations}), is the category-theoretic analog of a fibre bundle. Concretely, a fibration consists of a functor from a category, called the total category, to another category, called the base category, which satisfies a universal property analog to the path-lifting property of the universal covering space of a topological space. An equivalent way of understanding a fibration is to regard it as a collection of categories, called \textit{fibres}, indexed by the objects of the base category, ``glued'' together.
\par This interpretation is made precise by the Grothendieck construction, which establishes an equivalence between (cloven) fibrations and indexed categories. Concretely, the indexed category associated with a (cloven) fibration is the pseudofunctor which sends an object of the base category of the fibration to the fibre over the given object.
\par In this section, we recall the definitions of a (cloven) fibration, of an indexed category, and we recall the Grothendieck construction as presented in~\cite[Theorem~8.3.1]{borceux:fibrations}.
\par Let us briefly establish some of the jargon we adopt in the rest of the paper.

\begin{notation}
\label{notation:fibres-lifts}
Given a functor $\Pi\colon\X'\to\X$, we refer to $\X$ as the \textbf{base category} of $\Pi$ and to $\X'$ as the \textbf{total category} of $\Pi$. For an object $A$ of $\X$, an object $E$ of $\X'$ is \textbf{over} $A$ if $\Pi(E)=A$ and similarly, a morphism $\varphi\colon E\to E'$ of $\X'$ is \textbf{over} a morphism $f\colon A\to A'$ if $\Pi(\varphi)=f$. Equivalently, we also say that $\varphi$ is a \textbf{lift} of $f$. A morphism of $\X'$ is \textbf{vertical} if it is over the identity morphism of an object of $\X$. The \textbf{fibre} over an object $A$ of $\X$ refers to the category $\Pi^{-1}(A)$ whose objects are objects of $\X'$ over $A$ and morphisms are morphisms of $\X'$ over $\id_A$.
\end{notation}

To introduce the notion of a fibration, we first need to recall the definition of a cartesian morphism.

\begin{definition}
\label{definition:cartesian-morphism}
Given a functor $\Pi\colon\X'\to\X$, a morphism $\varphi\colon\tilde E\to E'$ of $\X'$ over a morphism $f\colon\tilde A\to A'$, is \textbf{cartesian} if for any given morphism $\psi\colon E\to E'$ of $\X$ over a morphism $g\colon A\to A'$, and any morphism $h\colon A\to\tilde A$ of $\X$ making the following diagram commutes
\begin{equation*}
% https://q.uiver.app/#q=WzAsMyxbMCwxLCJcXHRpbGRlIEEiXSxbMSwxLCJBJyJdLFswLDAsIkEiXSxbMiwwLCJoIiwyXSxbMCwxLCJmIiwyXSxbMiwxLCJnIl1d
\begin{tikzcd}
A \\
{\tilde A} & {A'}
\arrow["h"', from=1-1, to=2-1]
\arrow["g", from=1-1, to=2-2]
\arrow["f"', from=2-1, to=2-2]
\end{tikzcd}
\end{equation*}
there is a unique lift $\xi\colon E\to\tilde E$ of $h$ such that:
\begin{equation*}
% https://q.uiver.app/#q=WzAsMyxbMCwxLCJcXHRpbGRlIEUiXSxbMSwxLCJFJyJdLFswLDAsIkUiXSxbMiwwLCJcXHhpIiwyLHsic3R5bGUiOnsiYm9keSI6eyJuYW1lIjoiZGFzaGVkIn19fV0sWzAsMSwiXFx2YXJwaGkiLDJdLFsyLDEsIlxccHNpIl1d
\begin{tikzcd}
E \\
{\tilde E} & {E'}
\arrow["\xi"', dashed, from=1-1, to=2-1]
\arrow["\psi", from=1-1, to=2-2]
\arrow["\varphi"', from=2-1, to=2-2]
\end{tikzcd}
\end{equation*}
The map $\varphi$ is also called a \textbf{cartesian lift} of the morphism $f$ on $E'$.
\end{definition}

\begin{definition}
\label{definition:fibration}
A \textbf{fibration} consists of a functor $\Pi\colon\X'\to\X$ for which for every pair $(f,E')$ formed by a morphism $f\colon A\to A'$ of $\X$ and by an object $E'$ of $\X'$ over $A'$, there is a cartesian lift of $f$ whose codomain is $E'$. A \textbf{cleavage} of a fibration is a choice of a cartesian lift $\varphi\^f_{E'}$ for each pair $(f,E')$ formed by a morphism $f\colon A\to A'$ of $\X$ and of an object $E'$ over $A'$. Finally, a \textbf{cloven fibration} is a fibration equipped with a cleavage.
\end{definition}

\begin{notation}
\label{notation:fibration}
In the following, a fibration $\Pi\colon\X'\to\X$ is denoted by $(\X,\X',\Pi)$. Moreover, for a cloven fibration $(\X,\X',\Pi)$, the cartesian lift of a morphism $f\colon A\to A'$ specified by the cleveage, corresponding to an object $E'$ over $A'$, is denoted by:
\begin{align*}
&\varphi\^f_{E'}\colon f^\*E'\to E'
\end{align*}
For the sake of simplicity, when the object $E'$ is clear by the context, we omit the subscript $_{E'}$. When the symbol adopted for a cloven fibration is decorated with a subscript or with a superscript, the same decoration is applied to the cartesian lifts. For example, for a cloven fibration $(\X_\s,\X_\s',\Pi_\s)$, a morphism $f\colon A\to A'$ of $\X_\s$, and an object $E'$ over $A'$, the corresponding cartesian lift is denoted by $\varphi_\s\^f$.
\end{notation}

Every morphism $\psi\colon E\to E'$ of the total category $\X'$ of a fibration $(\X,\X',\Pi)$ admits a decomposition into a vertical morphism $\xi\colon E\to\tilde E$ and a cartesian morphism $\varphi\colon\tilde E\to E'$. To see why, consider the morphism $(f\colon A\to A')\=\Pi(\psi)$. Since $\Pi$ is a fibration, there is a cartesian lift $\varphi\colon\tilde E\to E'$ of $f$ whose codomain is $E'$. In particular, since $\varphi$ is a lift of $f$, $\Pi(\varphi)\o\id_A=\Pi(\varphi)=f$. For the cartesian property of $\varphi$, there is a unique morphism $\xi\colon E\to\tilde E$ of $\X$ such that $\Pi(\xi)=\id_A$ and $\varphi\o\xi=\psi$. Since $\xi$ is a lift of $\id_A$, $\xi$ is vertical.
\par This decomposition is not unique, since there could be more than one cartesian lift of the morphism $f$. However, when the fibration is equipped with a cleavage, there is a canonical choice of such a decomposition. In particular, for a cloven fibration, each morphism $\psi\colon E\to E'$ of the total category can be decomposed as follows:
\begin{equation*}
% https://q.uiver.app/#q=WzAsMyxbMCwwLCJFIl0sWzAsMSwiZl5cXCpFJyJdLFsxLDEsIkUnIl0sWzAsMSwiXFx4aV9mIiwyXSxbMSwyLCJcXHZhcnBoaV9mIiwyXSxbMCwyLCJcXHBzaSJdXQ==
\begin{tikzcd}
E \\
{f^\*E'} & {E'}
\arrow["{\xi\^f}"', from=1-1, to=2-1]
\arrow["\psi", from=1-1, to=2-2]
\arrow["{\varphi\^f}"', from=2-1, to=2-2]
\end{tikzcd}
\end{equation*}
However, since the cartesian lift $\varphi\^f\colon f^\*E'\to E'$ is determined by the cleavage, by $f$, and by $E'$, each morphism $\psi\colon E\to E'$ of the total category is fully specified by a pair $(f,\xi\^f)$ formed by a morphism $f\=\Pi(\psi)$ of $\X$ and by a vertical morphism $\xi\^f\colon E\to f^\*E'$. This is precisely the intuition underpinning the Grothendieck construction.
\par The first step is to ``split'' a fibration into the collection of its fibres. Such a collection has the structure of an indexed category. Informally, an indexed category on a category $\X$ consists of a collection of categories indexed by the objects of $\X$ together with a collection of functors between these categories, indexed by the morphisms of $\X$, and for which the indexing of the categories and of the functors are compatible with the composition and the identity morphisms of $\X$.

\begin{definition}
\label{definition:indexed-category}
An \textbf{indexed category} on a category $\X$ is a pseudofunctor $\I\colon\X^\op\to\Cat$, where $\Cat$ denotes the $2$-category of categories. Concretely, this consists of the following data:
\begin{description}
\item[categories] For each object $A$ of $\X$, a category $\X\^A$;

\item[substitution functors] for each morphism $f\colon A\to A'$ of $\X$, a functor $f^\ast\colon\X\^{A'}\to\X\^A$, called a \textbf{substitution functor};

\item[unitor] For each object $A$ of $\X$, a natural isomorphism $\Unit\colon\id_{\X\^A}\cong\id_A^\*$, called the \textbf{unitor};

\item[compositor] For each pair of composable morphisms $f\colon A\to A'$ and $g\colon B\to C$ of $\X$, a natural isomorphism $\Comp\colon f^\*\o g^\*\cong(g\o f)^\*$, called the \textbf{compositor};
\end{description}
satisfying some compatibility conditions. We invite the reader to consult~\cite[Section~I,3.2]{gray:formal-cats-fibrations} for details.
\end{definition}

\begin{notation}
\label{notation:indexed-category}
In the following, an indexed category $\I\colon\X^\op\to\Cat$ is denoted by $(\X,\I)$, the category $\I(A)$ associated with an object $A$ of $\X$ is denoted by $\X\^A$, the substitution functor $\I(f)$ induced by a morphism $f\colon A\to A'$ of $\X$ is denoted by $f^\*\colon\X\^{A'}\to\X\^{A}$,and  the unitor and the compositor of $\I$ are denoted by $\Unit$ and $\Comp$, respectively. Moreover, when the symbol used to denote an indexed category is decorated with a subscript or with a superscript, the same decoration is applied to the categories and the substitution functors. For example, for an indexed category $(\X_\s,\I_\s)$, the corresponding categories are denoted by $\X_\s\^A$ for each $A$ of $\X_\s$, and the substitution functors are denoted by $f_\s^\*$. The same convention extends to the unitor and the compositor, e.g., $\Unit_\s$ and $\Comp_\s$.
\end{notation}

Each cloven fibration $(\X,\X',\Pi)$ is associated with an indexed category, denoted by $\II(\Pi)\colon\X^\op\to\Cat$, defined as follows:
\begin{description}
\item[categories] Each object $A$ of $\X$ is associated with the corresponding fibre $\Pi^\-(A)$;

\item[substitution functors] Each morphism $f\colon A\to A'$ of $\X$ is associated with the functor $f^\*\colon\Pi^\-(A')\to\Pi^\-(A)$ which sends each object $E'$ over $A'$ to the object $f^\*E'$, domain of the cartesian lift $\varphi\^f\colon f^\*E'\to E'$ of $f$, and which sends any vertical morphism $\xi\colon E'\to E''$ over $A'$ to the unique vertical morphism $f^\*\xi\colon f^\*E'\to f^\*E''$ over $A$ such that $\varphi\^f_{E''}\o f^\*\xi=\varphi\^f_{E'}$, defined by the universality of the cartesian morphism $\varphi\^f_{E''}$;

\item[unitor] Given an object $E$ over $A$, the unitor is the unique vertical natural isomorphism $\Unit\colon E\to\id_A^\*E$ such that $\Unit\varphi\^{\id_A}=\id_E$, induced by the universality of the cartesian morphism $\varphi\^{\id_A}$;

\item[compositor] Given an object $E''$ over $A''$ and two morphisms $f\colon A\to A'$ and $g\colon A'\to A''$ of $\X$, the compositor is the unique vertical natural isomorphism $\Comp\colon f^\*g^\*E''\to(g\o f)^\*E''$ such that $\varphi\^{g\o f}\o\Comp=\varphi\^g\o\varphi\^f$, induced by the universality of the cartesian morphism $\varphi\^{g\o f}$.
\end{description}

The process of ``separating'' a cloven fibration into the indexed category $\II(\Pi)$ of its fibres is analogous to associating a fibre bundle $\pi\colon E\to M$ over a topological space $M$ with the assignment which sends each point $x$ of $M$ to the fibre $\pi^{-1}(x)$ over $x$. This process can be reversed via the Grothendieck construction.
\par Let us start by considering an indexed category $(\X,\I)$, which sends each object $A$ of $\X$ to the category $\X\^A$ and each morphism $f\colon A\to A'$ to the functor $f^\*\colon\X\^{A'}\to\X\^A$. Then, the \textit{category of elements} of $\I$ is the category $\El(\X,\I)$ defined as follows:
\begin{description}
\item[objects] The objects of $\El(\X,\I)$ are pairs $(A,E)$ formed by an object $A$ of $\X$ together with an object $E$ of $\X\^A$;

\item[morphisms] The morphisms of $\El(\X,\I)$ are pairs $(f,\xi\^f)\colon(A,E)\to(A',E')$ formed by a morphism $f\colon A\to A'$ of $\X$ together with a morphism $\xi\^f\colon E\to f^\*E'$ of $\X\^A$;

\item[identities] The identity morphism of a pair $(A,E)\in\El(\X,\I)$ is given by the pair $(\id_A,\xi\^{\id_A})$, where $\xi\^{\id_A}\colon E\to\id_A^\*E$ is the unitor of $\I$;

\item[composition] Given two composable morphisms $(f,\xi\^f)\colon(A,E)\to(A',E')$ and $(g,\xi\^g)\colon(A',E')\to(A'',E'')$ of $\El(\X,\I)$, their composition is the pair $(g\o f,\xi\^{g\o f})$ where
\begin{align*}
&\xi\^{g\o f}\colon E\xrightarrow{\xi\^f}f^\*E'\xrightarrow{f^\*\xi\^g}f^\*g^\*E''\xrightarrow{\Comp}(g\o f)^\*E''
\end{align*}
$\Comp$ denoting the compositor of $\I$.
\end{description}
The functor $\FF(\I)\colon\El(\X,\I)\to\X$ which sends each object $(A,E)$ of $\El(\X,\I)$ to the object $A$ of $\X$ and each morphism $(f,\xi\^f)$ of $\El(\X,\I)$ to the morphism $f$ of $\X$ is a cloven fibration.
\par The assignments $\II$ and $\FF$ which send a cloven fibration to an indexed category and vice versa, respectively, extend to an equivalence of $2$-categories. To see this, let us recall the definition of the $1$-morphisms and $2$-morphisms of indexed categories and fibrations.

\begin{definition}
\label{definition:1-morphisms-fibrations}
A \textbf{$1$-morphism of cloven fibrations} $(F,F')\colon(\X_\s,\X_\s',\Pi_\s)\to(\X_\b,\X_\b',\Pi_\b)$ consists of a pair of functors $F\colon\X_\s\to\X_\b$ and $F'\colon\X_\s'\to\X_\b'$ such that the following diagram commutes:
\begin{equation*}
% https://q.uiver.app/#q=WzAsNCxbMCwwLCJcXFhfXFxvJyJdLFswLDEsIlxcWF9cXG8iXSxbMSwwLCJcXFhfXFxiJyJdLFsxLDEsIlxcWF9cXGIiXSxbMCwxLCJcXFBpX1xcbyIsMl0sWzIsMywiXFxQaV9cXGIiXSxbMCwyLCJGJyJdLFsxLDMsIkYiLDJdXQ==
\begin{tikzcd}
{\X_\s'} & {\X_\b'} \\
{\X_\s} & {\X_\b}
\arrow["{F'}", from=1-1, to=1-2]
\arrow["{\Pi_\s}"', from=1-1, to=2-1]
\arrow["{\Pi_\b}", from=1-2, to=2-2]
\arrow["F"', from=2-1, to=2-2]
\end{tikzcd}
\end{equation*}
Moreover, $(F,F')$ is required to preserve the universality of the cartesian lifts defined by the cleavages of the fibrations. Concretely, this last condition means that, given a morphism $f\colon A\to A'$ of $\X_\s$ and an object $E'$ of the fibre over $A'$ w.r.t. $\Pi_\s$, the unique morphism $\kappa\^f\colon F'f_\s^\*(E')\to(Ff)_\b^\*(F'E')$ making the following diagram commutes
\begin{equation*}
% https://q.uiver.app/#q=WzAsMyxbMCwwLCJGJyhmX1xcb15cXCpFJykiXSxbMCwxLCIoRmYpXlxcKl9cXGIoRidFJykiXSxbMSwxLCJGJ0UnIl0sWzAsMiwiRlxcdmFycGhpX3tcXG8gZn0iXSxbMSwyLCJcXHZhcnBoaV97XFxiIEZmfSIsMl0sWzAsMSwiXFxrYXBwYV9mIiwyLHsic3R5bGUiOnsiYm9keSI6eyJuYW1lIjoiZGFzaGVkIn19fV1d
\begin{tikzcd}
{F'(f_\s^\*E')} \\
{(Ff)^\*_\b(F'E')} & {F'(E')}
\arrow["{\kappa\^f}"', dashed, from=1-1, to=2-1]
\arrow["{F\varphi\^f_\s}", from=1-1, to=2-2]
\arrow["{\varphi\^{Ff}_\b}"', from=2-1, to=2-2]
\end{tikzcd}
\end{equation*}
is an isomorphism, where $\varphi\^f_\s\colon f_\s^\*E'\to E'$ and $\varphi\^{Ff}_\b\colon(Ff)^\*_\b(F'E')\to F'E'$ are the cartesian lifts of $f$ on $E'$ w.r.t. $\Pi_\s$ and the cartesian lift of $Ff$ on $F'E'$ w.r.t. $\Pi_\b$, respectively. The isomorphisms $\kappa\^f$, indexed by the morphisms $f$ of $\X_\s$, are called the \textbf{distributors} of the morphism $(F,F')$.
\end{definition}

\begin{remark}
\label{remark:distributor-name}
The name \textit{distributor} was adopted in analogy with distributive laws and it should not be confused with the notion of profunctors.
\end{remark}

\begin{definition}
\label{definition:2-morphisms-fibrations}
A \textbf{$2$-morphism of fibrations} $(\theta,\theta')\colon(F,F')\to(G,G')$ between two $1$-morphisms $(F,F'),\linebreak(G,G')\colon(\X_\s,\X_\s',\Pi_\s)\to(\X_\b,\X_\b',\Pi_\b)$ of fibrations consists of two natural transformations $\theta\colon F\Rightarrow G$ and $\theta'\colon F'\Rightarrow G'$ such that:
\begin{align*}
&\Pi_\b\theta'=\theta_{\Pi_\s}
\end{align*}
\end{definition}

\begin{notation}
\label{notation:category-FIB}
We denote by $\Fib$ the $2$-category of cloven fibrations, and corresponding $1$-morphisms and $2$-morphisms. Given a category $\X$, we denote by $\Fib(\X)$ the $2$-category of cloven fibrations whose base category is $\X$, $1$-morphisms are $1$-morphisms of fibrations whose base functor $F\colon\X\to\X$ is the identity functor, and $2$-morphisms are $2$-morphisms of fibrations whose base natural transformation $\theta\colon\id_\X\Rightarrow\id_\X$ is the identity.
\end{notation}

\begin{definition}
\label{definition:1-morphisms-indexed-cats}
A \textbf{$1$-morphism of indexed categories} $(F,F',\kappa)\colon(\X_\s,\I_\s)\to(\X_\b,\I_\b)$ consists of a functor $F\colon\X_\s\to\X_\b$, a collection $F'$ of functors $F\^A\colon\X_\s\^A\to\X_\b\^{FA}$ indexed by the objects $A$ of $\X_\s$, from the category $\X_\s\^A\=\I_\s(A)$ to the category $\X_\b\^{FA}\=\I_\b(FA)$, and a collection $\kappa$ of natural isomorphisms $\kappa\^f\colon F\^A\o f^\*_\s\Rightarrow(Ff)^\*_\b\o F\^{A'}$ called \textbf{distributors}, indexed by the morphisms $f\colon A\to A'$ of $\X_\s$, where $f^\*_\s\=\I_\s(f)$ and $(Ff)^\*_\b\=\I_\b(Ff)$. Moreover, the distributors are compatible with the unitors $\Unit_\s,\Unit_\b$ and the compositors $\Comp_\s,\Comp_\b$ of $\I_\s$ and $\I_\b$, respectively, as follows
\begin{equation*}
% https://q.uiver.app/#q=WzAsNCxbMCwxLCJGXFxeQVxcaWRfe0FcXG99XlxcKkUiXSxbMSwxLCIoXFxpZF97RkF9KV5cXCpfXFxiIEZcXF5BRSJdLFswLDAsIkZcXF5BRSJdLFsxLDAsIihGXFxpZF9BKV5cXCpfXFxiIEZcXF5BRSJdLFswLDEsIlxca2FwcGFcXF57XFxpZF9BfSIsMl0sWzIsMCwiRlxcXkFcXFVuaXRfXFxvIiwyXSxbMiwzLCJcXFVuaXRfXFxiIEZcXF5BIl0sWzMsMSwiIiwwLHsibGV2ZWwiOjIsInN0eWxlIjp7ImhlYWQiOnsibmFtZSI6Im5vbmUifX19XV0=
\begin{tikzcd}
{F\^AE} & {(F\id_A)^\*_\b F\^AE} \\
{F\^A\id_{A\s}^\*E} & {(\id_{FA})^\*_\b F\^AE}
\arrow["{\Unit_\b F\^A}", from=1-1, to=1-2]
\arrow["{F\^A\Unit_\s}"', from=1-1, to=2-1]
\arrow[Rightarrow, no head, from=1-2, to=2-2]
\arrow["{\kappa\^{\id_A}}"', from=2-1, to=2-2]
\end{tikzcd}\hfill
% https://q.uiver.app/#q=WzAsNixbMCwxLCJGXFxeQShmZyleXFwqRScnIl0sWzIsMSwiKEZmRmcpXlxcKl9cXGIgRlxcXntBJyd9RScnIl0sWzAsMCwiRlxcXkFmXlxcKl9cXG8gZ15cXCpfXFxvIEUnJyJdLFsyLDAsIihGZileXFwqX1xcYiAoRmcpXlxcKl9cXGIgRlxcXntBJyd9RScnIl0sWzEsMCwiKEZmKV5cXCpfXFxiIEZcXF57QSd9Z15cXCpfXFxvIEUnJyJdLFsxLDEsIihGKGZnKSleXFwqX1xcYiBGXFxee0EnJ31FJyciXSxbMiwwLCJGXFxeQVxcQ29tcF9cXG8iLDJdLFszLDEsIlxcQ29tcF9cXGIgRlxcXntBJyd9Il0sWzQsMywiKEZmKV5cXCpfXFxiXFxrYXBwYVxcXmciXSxbMiw0LCJcXGthcHBhXFxeZiBnXlxcKl9cXG8iXSxbMCw1LCJcXGthcHBhXFxee2ZnfSIsMl0sWzUsMSwiIiwyLHsibGV2ZWwiOjIsInN0eWxlIjp7ImhlYWQiOnsibmFtZSI6Im5vbmUifX19XV0=
\begin{tikzcd}
{F\^Af^\*_\s g^\*_\s E''} & {(Ff)^\*_\b F\^{A'}g^\*_\s E''} & {(Ff)^\*_\b (Fg)^\*_\b F\^{A''}E''} \\
{F\^A(g\o f)^\*E''} & {(F(g\o f))^\*_\b F\^{A''}E''} & {(FfFg)^\*_\b F\^{A''}E''}
\arrow["{\kappa\^f g^\*_\s}", from=1-1, to=1-2]
\arrow["{F\^A\Comp_\s}"', from=1-1, to=2-1]
\arrow["{(Ff)^\*_\b\kappa\^g}", from=1-2, to=1-3]
\arrow["{\Comp_\b F\^{A''}}", from=1-3, to=2-3]
\arrow["{\kappa\^{g\o f}}"', from=2-1, to=2-2]
\arrow[Rightarrow, no head, from=2-2, to=2-3]
\end{tikzcd}
\end{equation*}
for every object $A$ of $\X_\s$, morphisms $f\colon A\to A'$ and $g\colon A'\to A''$ of $\X_\s$, and objects $E,E'$, and $E''$ objects of $\X_\s'$ over $A,A'$, and $A''$, respectively.
\end{definition}

\begin{definition}
\label{definition:2-morphisms-indexed-cats}
A \textbf{$2$-morphism of indexed categories} $(\theta,\theta')\colon(F,F',\kappa)\Rightarrow(G,G',\lambda)$ between two $1$-morphisms $(F,F',\kappa),(G,G',\lambda)\colon(\X_\s,\I_\s)\to(\X_\b,\I_\b)$ consists of a natural transformation $\theta\colon F\Rightarrow G$ and a collection of natural transformations $\theta\^A\colon F\^A\Rightarrow\theta^\*\o G\^A$ satisfying the following condition:
\begin{equation*}
% https://q.uiver.app/#q=WzAsNSxbMCwwLCJGXFxeQVxcbyBmX1xcb15cXCoiXSxbMSwwLCJcXHRoZXRhXlxcKlxcbyBHXFxeQVxcbyBmX1xcb15cXCoiXSxbMSwxLCJcXHRoZXRhXlxcKlxcbyhHZilfXFxiXlxcKlxcbyBHXFxee0EnfSJdLFsxLDIsIihGZilfXFxiXlxcKlxcb1xcdGhldGFeXFwqXFxvIEdcXF57QSd9Il0sWzAsMiwiKEZmKV9cXGJeXFwqXFxvIEZcXF57QSd9Il0sWzAsMSwiXFx0aGV0YVxcXkFmX1xcb15cXCoiXSxbMSwyLCJcXHRoZXRhXlxcKlxcbGFtYmRhXFxeZiJdLFsyLDMsIlxcZ2FtbWEiXSxbMCw0LCJcXGthcHBhXFxeZiIsMl0sWzQsMywiKEZmKV9cXGJeXFwqXFx0aGV0YVxcXntBJ30iLDJdXQ==
\begin{tikzcd}
{F\^A\o f_\s^\*} & {\theta^\*\o G\^A\o f_\s^\*} \\
& {\theta^\*\o(Gf)_\b^\*\o G\^{A'}} \\
{(Ff)_\b^\*\o F\^{A'}} & {(Ff)_\b^\*\o\theta^\*\o G\^{A'}}
\arrow["{\theta\^Af_\s^\*}", from=1-1, to=1-2]
\arrow["{\kappa\^f}"', from=1-1, to=3-1]
\arrow["{\theta^\*\lambda\^f}", from=1-2, to=2-2]
\arrow["\Comp", from=2-2, to=3-2]
\arrow["{(Ff)_\b^\*\theta\^{A'}}"', from=3-1, to=3-2]
\end{tikzcd}
\end{equation*}
\end{definition}

\begin{notation}
\label{notation:category-INDX}
We denote by $\Indx$ the $2$-category of indexed categories and corresponding $1$-morphisms and $2$-morphisms. Given a category $\X$, we denote by $\Indx(\X)$ the $2$-category of indexed categories whose base category is $\X$, $1$-morphisms are $1$-morphisms of indexed categories whose base functor $F\colon\X\to\X$ is the identity functor, and $2$-morphisms are $2$-morphisms of indexed categories whose base natural transformation $\theta\colon\id_\X\to\id_\X$ is the identity.
\par The Grothendieck construction is usually stated in terms of an equivalence between cloven fibrations and indexed categories over a fixed base category $\X$. We suggest the reader to consult~\cite[Section~8.3]{borceux:fibrations}. For the purposes of this paper, we need to extend this equivalence between cloven fibrations and indexed categories over any given base category.
\end{notation}

\begin{remark}
\label{remark:indexed-categories-vs-pseudofunctors}
In the literature, indexed categories are also described as pseudofunctors from the $1$-category $\Cat$ of categories to the $2$-category $\Cat$ of categories. Similarly, $1$-morphisms of indexed categories are known as pseudonatural transformations, and $2$-morphisms of indexed categories as modifications.
\end{remark}

\begin{proposition}
\label{proposition:classic-Grothendieck}
The assignment $\FF$ which sends an indexed category $(\X,\I)$ to the cloven fibration $(\X,\El(\X,\I),\FF(\I))$ and the assignment $\II$ which sends a cloven fibration $(\X,\X',\Pi)$ to the indexed category $(\X,\II(\Pi))$ extend to an equivalence of $2$-categories:
\begin{align*}
&\Fib\simeq\Indx
\end{align*}
\begin{proof}
Let us start by considering a $1$-morphism of fibrations $(F,F')\colon(\X_\s,\X_\s',\Pi_\s)\to(\X_\b,\X_\b',\Pi_\b)$ and let us define the corresponding $1$-morphism of indexed categories. The base functor is $F\colon\X_\s\to\X_\b$. The collection of functors $F\^A\colon\Pi_\s^\-(A)\to\Pi_\b^\-(FA)$ is given by the restriction of $F'$ to the fibres $\Pi_\s^\-(A)$ of $\Pi_\s$. Indeed, since $F\o\Pi_\s=\Pi_\b\o F'$, each $E\in\Pi_\s^\-(A)$ is sent to $F'E\in\Pi_\b^\-(FA)$ by $F'$. The distributor $\kappa\^f\colon F\^A\o f^\*_\s\Rightarrow(Ff)^\*_\b\o F\^{A'}$ is defined by the distributor
\begin{align*}
&\kappa\^f\colon F'(f^\*_\s E')\to(Ff)^\*_\b(F'E')
\end{align*}
induced by the universality of the cartesian lift $\varphi\^{Ff}_\b\colon(Ff)^\*_\b(F'E')\to F'(f^\*_\s E')$. To prove that $\kappa\^f$ is natural observe that, for a vertical morphism $\psi\colon E'\to E''$, the following diagram commutes:
\begin{equation*}
% https://q.uiver.app/#q=WzAsOCxbMCwzLCIoRmYpXlxcKl9cXGIoRidFJykiXSxbMSwyLCIoRmYpXlxcKl9cXGIoRidFJycpIl0sWzEsMSwiRicoZl5cXCpfXFxvIEUnJykiXSxbMCwwLCJGJyhmXlxcKl9cXG8gRScpIl0sWzIsMSwiRidFJyciXSxbMiwyLCJGJ0UnJyJdLFszLDAsIkYnRSciXSxbMywzLCJGJ0UnIl0sWzIsMSwiXFxrYXBwYVxcXmZfe0UnJ30iLDJdLFszLDAsIlxca2FwcGFcXF5mX3tFJ30iLDJdLFs0LDUsIiIsMix7ImxldmVsIjoyLCJzdHlsZSI6eyJoZWFkIjp7Im5hbWUiOiJub25lIn19fV0sWzYsNywiIiwyLHsibGV2ZWwiOjIsInN0eWxlIjp7ImhlYWQiOnsibmFtZSI6Im5vbmUifX19XSxbMyw2LCJGJ1xcdmFycGhpXFxeZl97XFxvIEUnfSJdLFsyLDQsIlxcdmFycGhpXFxee0ZmfV97XFxiIEUnJ30iXSxbMSw1LCJcXHZhcnBoaVxcXntGZn1fe1xcYiBFJyd9IiwyXSxbMCw3LCJGJ1xcdmFycGhpXFxeZl97XFxvIEUnfSIsMl0sWzYsNCwiRidcXHBzaSIsMV0sWzcsNSwiRidcXHBzaSIsMV0sWzMsMiwiRicoZl5cXCpfXFxvXFxwc2kpIiwxXSxbMCwxLCIoRmYpXlxcKl9cXGIoRidcXHBzaSkiLDFdXQ==
\begin{tikzcd}
{F'(f^\*_\s E')} &&& {F'E'} \\
& {F'(f^\*_\s E'')} & {F'E''} \\
& {(Ff)^\*_\b(F'E'')} & {F'E''} \\
{(Ff)^\*_\b(F'E')} &&& {F'E'}
\arrow["{F'\varphi\^f_{\s E'}}", from=1-1, to=1-4]
\arrow["{F'(f^\*_\s\psi)}"{description}, from=1-1, to=2-2]
\arrow["{\kappa\^f_{E'}}"', from=1-1, to=4-1]
\arrow["{F'\psi}"{description}, from=1-4, to=2-3]
\arrow[Rightarrow, no head, from=1-4, to=4-4]
\arrow["{\varphi\^{Ff}_{\b E''}}", from=2-2, to=2-3]
\arrow["{\kappa\^f_{E''}}"', from=2-2, to=3-2]
\arrow[Rightarrow, no head, from=2-3, to=3-3]
\arrow["{\varphi\^{Ff}_{\b E''}}"', from=3-2, to=3-3]
\arrow["{(Ff)^\*_\b(F'\psi)}"{description}, from=4-1, to=3-2]
\arrow["{F'\varphi\^f_{\s E'}}"', from=4-1, to=4-4]
\arrow["{F'\psi}"{description}, from=4-4, to=3-3]
\end{tikzcd}
\end{equation*}
So, in particular:
\begin{align*}
&\varphi\^{Ff}_{\b E''}\o\kappa\^f_{E''}\o F'(f_\s^\*\psi)=\varphi\^{Ff}_{\b E''}\o(Ff)^\*_\b(F'\psi)\o\kappa\^f_{E'}
\end{align*}
From the universality of the cartesian lift $\varphi\^{Ff}_{\b E''}$, we obtain:
\begin{align*}
&\kappa\^f_{E''}\o F'(f_\s^\*\psi)=(Ff)^\*_\b(F'\psi)\o\kappa\^f_{E'}
\end{align*}
Similarly, to prove the compatibility between the unitors, compositors and distributors, first notice that the following diagrams commute:
\begin{equation*}
% https://q.uiver.app/#q=WzAsNSxbMCwwLCJGJ0UiXSxbMSwwLCJGJ1xcaWRfe0FcXG99XlxcKiJdLFsyLDAsIihcXGlkX3tGQX0pXlxcKl9cXGIgRidFIl0sWzIsMSwiRidFIl0sWzAsMSwiKFxcaWRfe0ZBfSleXFwqX1xcYiBGJ0UiXSxbMCwxLCJGJ1xcVW5pdF9cXG8iXSxbMSwyLCJcXGthcHBhXFxee1xcaWRfQX0iXSxbMiwzLCJcXHZhcnBoaVxcXntcXGlkX3tGQX19X1xcYiJdLFsxLDMsIkYnXFx2YXJwaGlcXF57XFxpZF9BfV9cXG8iLDFdLFswLDQsIlxcVW5pdF9cXGIgRiciLDJdLFs0LDMsIlxcdmFycGhpXFxee1xcaWRfe0ZBfX1fXFxiIiwyXSxbMCwzLCIiLDIseyJjdXJ2ZSI6MSwibGV2ZWwiOjIsInN0eWxlIjp7ImhlYWQiOnsibmFtZSI6Im5vbmUifX19XV0=
\begin{tikzcd}
{F'E} & {F'\id_{A\s}^\*} & {(\id_{FA})^\*_\b F'E} \\
{(\id_{FA})^\*_\b F'E} && {F'E}
\arrow["{F'\Unit_\s}", from=1-1, to=1-2]
\arrow["{\Unit_\b F'}"', from=1-1, to=2-1]
\arrow[bend right=10, Rightarrow, no head, from=1-1, to=2-3]
\arrow["{\kappa\^{\id_A}}", from=1-2, to=1-3]
\arrow["{F'\varphi\^{\id_A}_\s}"{description}, from=1-2, to=2-3]
\arrow["{\varphi\^{\id_{FA}}_\b}", from=1-3, to=2-3]
\arrow["{\varphi\^{\id_{FA}}_\b}"', from=2-1, to=2-3]
\end{tikzcd}
\end{equation*}
\begin{equation*}
% https://q.uiver.app/#q=WzAsMTAsWzAsMCwiRidmXlxcKl9cXG8gZ15cXCpfXFxvIEUnJyJdLFsxLDAsIihGZileXFwqX1xcYiBGJ2deXFwqX1xcbyBFJyciXSxbMiwwLCIoRmYpXlxcKl9cXGIoRmcpXlxcKl9cXGIgRidFJyciXSxbMywwLCIoRihmZykpXlxcKl9cXGIgRidFJyciXSxbMywyLCJGJ0UnJyJdLFsyLDEsIihGZyleXFwqX1xcYiBGJ0UnJyJdLFsxLDEsIkYnZ15cXCpfXFxvIEUnJyJdLFsyLDIsIkYnRScnIl0sWzAsMiwiRicoZmcpXlxcKl9cXG8gRScnIl0sWzEsMiwiKEYoZmcpKV5cXCpfXFxiIEYnRScnIl0sWzIsMywiXFxDb21wX1xcYiBGJyJdLFszLDQsIlxcdmFycGhpXFxee0YoZmcpfV9cXGIiXSxbMiw1LCJcXHZhcnBoaVxcXntGZn1fXFxiKEZnKV5cXCpGJyJdLFsxLDIsIihGZileXFwqX1xcYlxca2FwcGFcXF5nIl0sWzAsMSwiXFxrYXBwYVxcXmZnXlxcKl9cXG8iXSxbMSw2LCJcXHZhcnBoaVxcXntGZn1fXFxiIEYnZ15cXCpfXFxvIl0sWzYsNSwiXFxrYXBwYVxcXmciLDFdLFs1LDcsIlxcdmFycGhpXFxee0ZnfV9cXGIgRiciXSxbNyw0LCIiLDIseyJsZXZlbCI6Miwic3R5bGUiOnsiaGVhZCI6eyJuYW1lIjoibm9uZSJ9fX1dLFs2LDcsIkYnXFx2YXJwaGlcXF5nX1xcbyIsMV0sWzAsNiwiRidcXHZhcnBoaVxcXmZfXFxvIGdeXFwqX1xcbyIsMV0sWzAsOCwiRidcXENvbXBfXFxvIiwyXSxbOCw5LCJcXGthcHBhXFxee2ZnfSIsMl0sWzksNywiXFx2YXJwaGlcXF57RihmZyl9X1xcYiIsMl1d
\begin{tikzcd}
{F'f^\*_\s g^\*_\s E''} & {(Ff)^\*_\b F'g^\*_\s E''} & {(Ff)^\*_\b(Fg)^\*_\b F'E''} & {(F(g\o f))^\*_\b F'E''} \\
& {F'g^\*_\s E''} & {(Fg)^\*_\b F'E''} \\
{F'(g\o f)^\*_\s E''} & {(F(g\o f))^\*_\b F'E''} & {F'E''} & {F'E''}
\arrow["{g^\*_\s\o\kappa\^f}", from=1-1, to=1-2]
\arrow["{F'\varphi\^f_\s g^\*_\s}"{description}, from=1-1, to=2-2]
\arrow["{F'\Comp_\s}"', from=1-1, to=3-1]
\arrow["{(Ff)^\*_\b\kappa\^g}", from=1-2, to=1-3]
\arrow["{\varphi\^{Ff}_\b F'g^\*_\s}", from=1-2, to=2-2]
\arrow["{\Comp_\b F'}", from=1-3, to=1-4]
\arrow["{\varphi\^{Ff}_\b(Fg)^\*F'}", from=1-3, to=2-3]
\arrow["{\varphi\^{F(g\o f)}_\b}", from=1-4, to=3-4]
\arrow["{\kappa\^g}"{description}, from=2-2, to=2-3]
\arrow["{F'\varphi\^g_\s}"{description}, from=2-2, to=3-3]
\arrow["{\varphi\^{Fg}_\b F'}", from=2-3, to=3-3]
\arrow["{\kappa\^{g\o f}}"', from=3-1, to=3-2]
\arrow["{\varphi\^{F(g\o f)}_\b}"', from=3-2, to=3-3]
\arrow[Rightarrow, no head, from=3-3, to=3-4]
\end{tikzcd}
\end{equation*}
From the universality of the cartesian lifts we conclude:
\begin{align*}
&\kappa\^{\id_A}\o F\^A\Unit_\s={\Unit_\b}_{F\^A}\\
&\Comp_{\b F\^{A''}}\o(Ff)^\*_\b\kappa\^g\o\kappa\^f_{g^\*_\s}=\kappa\^{g\o f}\o F\^A\Comp_\s
\end{align*}
Let us now consider a $2$-morphism $(\theta,\theta')\colon(F,F')\Rightarrow(G,G')$ of fibrations between two $1$-morphisms $(F,F'),\linebreak(G,G')\colon(\X_\s,\X_\s',\Pi_\s)\to(\X_\b,\X_\b',\Pi_\b)$ of fibrations and let us define the corresponding $2$-morphism
\begin{align*}
&\II(\theta,\theta')\colon\II(F,F')\to\II(G,G')
\end{align*}
of indexed categories between the two $1$-morphisms $(F,F',\kappa),(G,G',\lambda)\colon\II(\X_\s,\X_\s',\Pi_\s)\to\II(\X_\b,\X_\b',\Pi_\b)$. The base natural transformation is $\theta\colon F\Rightarrow G$. To define the natural transformations $\theta\^A\colon F\^A\Rightarrow\theta^\*\o G\^A$, first notice that, since $\Pi_\b\theta'=\theta_{\Pi_\s}$, $\theta'\colon F'E\to G'E$ is over $\theta\colon FA\to GA$, for $A=\Pi_\s E$. Therefore, the universality of the cartesian lift $\varphi\^\theta\colon\theta^\*G'E\to G'E$ defines a unique morphism $\theta\^A_E\colon F'E\to\theta^\*G'E$ making the following diagram commutes:
\begin{equation*}
% https://q.uiver.app/#q=WzAsMyxbMCwxLCJcXHRoZXRhXlxcKkcnRSJdLFsxLDEsIkcnRSJdLFswLDAsIkYnRSJdLFswLDEsIlxcdmFycGhpXFxeXFx0aGV0YSIsMl0sWzIsMSwiXFx0aGV0YSciXSxbMiwwLCJcXHRoZXRhXFxeQSIsMix7InN0eWxlIjp7ImJvZHkiOnsibmFtZSI6ImRhc2hlZCJ9fX1dXQ==
\begin{tikzcd}
{F'E} \\
{\theta^\*G'E} & {G'E}
\arrow["{\theta\^A}"', dashed, from=1-1, to=2-1]
\arrow["{\theta'}", from=1-1, to=2-2]
\arrow["{\varphi\^\theta}"', from=2-1, to=2-2]
\end{tikzcd}
\end{equation*}
To show that $\theta\^A$ is natural, consider a vertical morphism $\psi\colon E\to E'$:
\begin{equation*}
% https://q.uiver.app/#q=WzAsOCxbMCwzLCJcXHRoZXRhXlxcKkcnRSJdLFs0LDMsIkcnRSJdLFswLDAsIkYnRSJdLFs0LDAsIkcnRSJdLFsxLDEsIkYnRSciXSxbMywxLCJHJ0UnIl0sWzMsMiwiRydFJyJdLFsxLDIsIlxcdGhldGFeXFwqRydFJyJdLFswLDEsIlxcdmFycGhpXFxeXFx0aGV0YV9FIiwyXSxbMiwwLCJcXHRoZXRhXFxeQV9FIiwyLHsic3R5bGUiOnsiYm9keSI6eyJuYW1lIjoiZGFzaGVkIn19fV0sWzIsMywiXFx0aGV0YSdfRSJdLFszLDEsIiIsMCx7ImxldmVsIjoyLCJzdHlsZSI6eyJoZWFkIjp7Im5hbWUiOiJub25lIn19fV0sWzcsNiwiXFx2YXJwaGlcXF5cXHRoZXRhX3tFJ30iXSxbNCw1LCJcXHRoZXRhJ197RSd9IiwyXSxbNCw3LCJcXHRoZXRhXFxeQV97RSd9IiwwLHsic3R5bGUiOnsiYm9keSI6eyJuYW1lIjoiZGFzaGVkIn19fV0sWzUsNiwiIiwyLHsibGV2ZWwiOjIsInN0eWxlIjp7ImhlYWQiOnsibmFtZSI6Im5vbmUifX19XSxbMyw1LCJHJ1xccHNpIiwxXSxbMiw0LCJGJ1xccHNpIiwxXSxbMCw3LCJcXHRoZXRhXlxcKkcnXFxwc2kiLDFdLFsxLDYsIkcnXFxwc2kiLDFdXQ==
\begin{tikzcd}
{F'E} &&&& {G'E} \\
& {F'E'} && {G'E'} \\
& {\theta^\*G'E'} && {G'E'} \\
{\theta^\*G'E} &&&& {G'E}
\arrow["{\theta'_E}", from=1-1, to=1-5]
\arrow["{F'\psi}"{description}, from=1-1, to=2-2]
\arrow["{\theta\^A_E}"', dashed, from=1-1, to=4-1]
\arrow["{G'\psi}"{description}, from=1-5, to=2-4]
\arrow[Rightarrow, no head, from=1-5, to=4-5]
\arrow["{\theta'_{E'}}"', from=2-2, to=2-4]
\arrow["{\theta\^A_{E'}}", dashed, from=2-2, to=3-2]
\arrow[Rightarrow, no head, from=2-4, to=3-4]
\arrow["{\varphi\^\theta_{E'}}", from=3-2, to=3-4]
\arrow["{\theta^\*G'\psi}"{description}, from=4-1, to=3-2]
\arrow["{\varphi\^\theta_E}"', from=4-1, to=4-5]
\arrow["{G'\psi}"{description}, from=4-5, to=3-4]
\end{tikzcd}
\end{equation*}
In particular
\begin{align*}
&\varphi\^\theta_{E'}\o\theta^\*G'\psi\o\theta\^A_E=\varphi\^\theta_{E'}\o\theta\^A_{E'}\o F'\psi
\end{align*}
and from the universality of $\varphi\^\theta_{E'}\colon\theta^\*G'E'\to G'E'$ we deduce:
\begin{align*}
&\theta^\*G'\psi\o\theta\^A_E=\theta\^A_{E'}\o F'\psi
\end{align*}
To prove the compatibility between the distributors and $\theta'$, observe that the following diagrams commute:
\begin{equation*}
% https://q.uiver.app/#q=WzAsOSxbMCwwLCJGJ2ZeXFwqX1xcbyBFJyJdLFsxLDAsIlxcdGhldGFeXFwqRydmXlxcKl9cXG8gRSciXSxbMiwwLCJcXHRoZXRhXlxcKihHZileXFwqX1xcYiBHJ0UnIl0sWzMsMCwiKEZmKV5cXCpfXFxiXFx0aGV0YV5cXCogRydFJyJdLFszLDEsIlxcdGhldGFeXFwqIEcnRSciXSxbMywyLCJHJ0UnIl0sWzIsMSwiKEdmKV5cXCpfXFxiIEcnRSciXSxbMiwyLCJHJ0UnIl0sWzEsMSwiRydmXlxcKl9cXG8gRSciXSxbMCwxLCJcXHRoZXRhXFxeQWZeXFwqX1xcbyJdLFsxLDIsIlxcdGhldGFeXFwqXFxsYW1iZGFcXF5mIl0sWzIsMywiIiwwLHsibGV2ZWwiOjIsInN0eWxlIjp7ImhlYWQiOnsibmFtZSI6Im5vbmUifX19XSxbMyw0LCJcXHZhcnBoaVxcXntGZn1fXFxiXFx0aGV0YV5cXCpHJyJdLFs0LDUsIlxcdmFycGhpXFxeXFx0aGV0YSBHJyJdLFsyLDYsIlxcdmFycGhpXFxeXFx0aGV0YShHZileXFwqX1xcYiBHJyIsMV0sWzYsNywiXFx2YXJwaGlcXF57RmZ9X1xcYiBHJyIsMV0sWzcsNSwiIiwxLHsibGV2ZWwiOjIsInN0eWxlIjp7ImhlYWQiOnsibmFtZSI6Im5vbmUifX19XSxbOCw2LCJcXGxhbWJkYVxcXmYiLDFdLFsxLDgsIlxcdmFycGhpXFxeXFx0aGV0YSBHJ2ZeXFwqX1xcbyIsMV0sWzgsNywiRydcXHZhcnBoaVxcXmZfXFxvIiwyXSxbMCw4LCJcXHRoZXRhJ2ZeXFwqX1xcbyIsMl1d
\begin{tikzcd}
{F'f^\*_\s E'} & {\theta^\*G'f^\*_\s E'} & {\theta^\*(Gf)^\*_\b G'E'} & {(Ff)^\*_\b\theta^\* G'E'} \\
& {G'f^\*_\s E'} & {(Gf)^\*_\b G'E'} & {\theta^\* G'E'} \\
&& {G'E'} & {G'E'}
\arrow["{\theta\^Af^\*_\s}", from=1-1, to=1-2]
\arrow["{\theta'f^\*_\s}"', from=1-1, to=2-2]
\arrow["{\theta^\*\lambda\^f}", from=1-2, to=1-3]
\arrow["{\varphi\^\theta G'f^\*_\s}"{description}, from=1-2, to=2-2]
\arrow[Rightarrow, no head, from=1-3, to=1-4]
\arrow["{\varphi\^\theta(Gf)^\*_\b G'}"{description}, from=1-3, to=2-3]
\arrow["{\varphi\^{Ff}_\b\theta^\*G'}", from=1-4, to=2-4]
\arrow["{\lambda\^f}"{description}, from=2-2, to=2-3]
\arrow["{G'\varphi\^f_\s}"', from=2-2, to=3-3]
\arrow["{\varphi\^{Ff}_\b G'}"{description}, from=2-3, to=3-3]
\arrow["{\varphi\^\theta G'}", from=2-4, to=3-4]
\arrow[Rightarrow, no head, from=3-3, to=3-4]
\end{tikzcd}
\end{equation*}
\begin{equation*}
% https://q.uiver.app/#q=WzAsNyxbMCwwLCJGJ2ZeXFwqX1xcbyBFJyJdLFsxLDAsIihGZileXFwqX1xcYiBGJ0UnIl0sWzIsMCwiKEZmKV5cXCpfXFxiXFx0aGV0YV5cXCpHJ0UnIl0sWzIsMSwiXFx0aGV0YV5cXCpHJ0UnIl0sWzIsMiwiRydFJyJdLFsxLDEsIkYnRSciXSxbMCwyLCJHJ2ZeXFwqX1xcbyBFJyJdLFsyLDMsIlxcdmFycGhpXFxee0ZmfV9cXGJcXHRoZXRhXlxcKkcnIl0sWzMsNCwiXFx2YXJwaGlcXF5cXHRoZXRhIEcnIl0sWzEsMiwiKEZmKV5cXCpfXFxiXFx0aGV0YVxcXntBJ30iXSxbMCwxLCJcXGthcHBhXFxeZiJdLFsxLDUsIlxcdmFycGhpXFxee0ZmfV9cXGIgRiciLDFdLFs1LDMsIlxcdGhldGFcXF57QSd9IiwxXSxbNSw0LCJcXHRoZXRhJyIsMV0sWzAsNSwiRidcXHZhcnBoaVxcXmZfXFxvIiwxXSxbMCw2LCJcXHRoZXRhJ2ZeXFwqX1xcbyIsMl0sWzYsNCwiRydcXHZhcnBoaVxcXmZfXFxvIiwyXV0=
\begin{tikzcd}
{F'f^\*_\s E'} & {(Ff)^\*_\b F'E'} & {(Ff)^\*_\b\theta^\*G'E'} \\
& {F'E'} & {\theta^\*G'E'} \\
{G'f^\*_\s E'} && {G'E'}
\arrow["{\kappa\^f}", from=1-1, to=1-2]
\arrow["{F'\varphi\^f_\s}"{description}, from=1-1, to=2-2]
\arrow["{\theta'f^\*_\s}"', from=1-1, to=3-1]
\arrow["{(Ff)^\*_\b\theta\^{A'}}", from=1-2, to=1-3]
\arrow["{\varphi\^{Ff}_\b F'}"{description}, from=1-2, to=2-2]
\arrow["{\varphi\^{Ff}_\b\theta^\*G'}", from=1-3, to=2-3]
\arrow["{\theta\^{A'}}"{description}, from=2-2, to=2-3]
\arrow["{\theta'}"{description}, from=2-2, to=3-3]
\arrow["{\varphi\^\theta G'}", from=2-3, to=3-3]
\arrow["{G'\varphi\^f_\s}"', from=3-1, to=3-3]
\end{tikzcd}
\end{equation*}
So, we have
\begin{align*}
&\varphi\^\theta_{G'}\o\varphi\^{Ff}_{\b\theta^\*G'}\o\theta^\*\lambda\^f\o\theta\^A_{f^\*_\s}=\varphi\^\theta_{G'}\o\varphi\^{Ff}_{\b\theta^\*G'}\o(Ff)^\*_\b\theta\^{A'}\o\kappa\^f
\end{align*}
and, from the universality of the cartesian lifts, we obtain
\begin{align*}
&\theta^\*\lambda\^f\o\theta\^A_{f^\*_\s}=(Ff)^\*_\b\theta\^{A'}\o\kappa\^f
\end{align*}
which is precisely the compatibility between $\theta'$ and the distributors.
\par Conversely, consider a $1$-morphism $(F,F',\kappa)\colon(\X_\s,\I_\s)\to(\X_\b,\I_\b)$ of indexed categories and let us define the corresponding $1$-morphism of fibrations. Let us start with the base functor, which is just the base functor $F\colon\X_\s\to\X_\b$. To define the total functor $F'\colon\X_\s'\to\X_\b'$, let us consider an object $(A,E)\in\El(\X_\s,\I_\s)$ and let us define
\begin{align*}
&F'(A,E)\=(FA,F\^AE)
\end{align*}
where $F\^A\colon\X_\s\^A\to\X_\b\^{FA}$, so $(FA,F\^AE)$ is indeed an object of $\El(\X_\b,\I_\b)$. Now, consider a morphism $(f,\xi\^f)\colon(A,E)\to(A',E')$ of $\El(\X_\s,\I_\s)$. Let us define $F'(f,\xi\^f)$ as the morphism $(Ff,\xi\^{Ff})$, where
\begin{align*}
&\xi\^{Ff}\colon F\^AE\xrightarrow{F\^A\xi\^f}F\^Af^\*_\s E'\xrightarrow{\kappa\^f}(Ff)^\*_\b F\^{A'}E'
\end{align*}
where we employed the distributor $\kappa\^f\colon F\^A\o f^\*_\s\Rightarrow (Ff)^\*_\b\o F\^{A'}$. To prove that $F\^A$ is functorial, we need to employ the compatibilities of the distributors with the unitors and the compositors. Let us start by showing that $F'\id_{(A,E)}=\id_{F'(A,E)}$:
\begin{eqnarray*}
& &F'\left(\id_{(A,E)}\right)\\
&=&F'\left(\id_A,\Unit_\s\colon E\to\id_{A\s}^\*E\right)\\
&=&\left(F\id_A,F\^AE\xrightarrow{F\^A\Unit_\s}F\^A\id_{A\s}^\*E\xrightarrow{\kappa\^{\id_A}}(F\id_A)^\*_\b F\^AE\right)\\
&=&\left(\id_{FA},\Unit_{\b F\^A}\colon F\^AE\to(\id_{FA})^\*_\b F\^AE\right)\\
&=&\id_{F'(A,E)}
\end{eqnarray*}
Consider now two morphisms $(f,\xi\^f)\colon(A,E)\to(A',E')$ and $(g,\xi\^g)\colon(A',E')\to(A'',E'')$ of $\El(\X_\s,\I_\s)$. First, recall that $F'(f,\xi\^f)=(Ff,\xi\^{Ff})$ and $F'(g,\xi\^g)=(Fg,\xi\^{Fg})$, where:
\begin{align*}
&\xi\^{Ff}\colon F\^AE\xrightarrow{F\^A\xi\^f}F\^Af^\*_\s E'\xrightarrow{\kappa\^f}(Ff)^\*_\b F\^{A'}E'\\
&\xi\^{Fg}\colon F\^{A'}E'\xrightarrow{F\^{A'}\xi\^g}F\^{A'}g^\*_\s E''\xrightarrow{\kappa\^g}(Fg)^\*_\b F\^{A''}E''
\end{align*}
Moreover:
\begin{align*}
&\left(g,\xi\^g\right)\o\left(f,\xi\^f\right)=\left(g\o f,E\xrightarrow{\xi\^f}f^\*_\s E'\xrightarrow{f^\*_\s\xi\^g}f^\*_\s g^\*_\s E''\xrightarrow{\Comp_\s}(g\o f)^\*_\s E''\right)
\end{align*}
Therefore, $F'\left(\left(g,\xi\^g\right)\o\left(f,\xi\^f\right)\right)=\left(F(g\o f),\xi\^{F(g\o f)}\right)$, where:
\begin{align*}
&\xi\^{F(g\o f)}\colon F\^AE\xrightarrow{F\^A\xi\^f}F\^Af^\*_\s E'\xrightarrow{F\^Af^\*_\s\xi\^g}F\^Af^\*_\s g^\*_\s E''\xrightarrow{F\^A\Comp_\s}F\^A(g\o f)^\*_\s E''\xrightarrow{\kappa\^{g\o f}}(F(g\o f))^\*_\b F\^{A''}E''
\end{align*}
However, the following diagram commutes:
\begin{equation*}
% https://q.uiver.app/#q=WzAsOSxbMCwwLCJGXFxeQUUiXSxbMSwwLCJGXFxeQWZeXFwqX1xcbyBFJyJdLFsyLDAsIihGZileXFwqX1xcYiBGXFxee0EnfUUnIl0sWzMsMCwiKEZmKV5cXCpfXFxiIEZcXF57QSd9Z15cXCpfXFxvIEUnJyJdLFszLDEsIihGZileXFwqX1xcYihGZyleXFwqX1xcYiBGXFxee0EnJ31FJyciXSxbMywzLCIoRihmZykpXlxcKl9cXGIgRlxcXntBJyd9RScnIl0sWzIsMiwiRlxcXkEoZmcpXlxcKl9cXG8gRScnIl0sWzIsMSwiRlxcXkFmXlxcKl9cXG8gZ15cXCpfXFxvIEUnJyJdLFswLDMsIkZcXF5BRSJdLFswLDEsIkZcXF5BXFx4aVxcXmYiXSxbMSwyLCJcXGthcHBhXFxeZiJdLFsyLDMsIihGZileXFwqX1xcYiBGXFxee0EnfVxceGlcXF5nIl0sWzMsNCwiKEZmKV5cXCpfXFxiXFxrYXBwYVxcXmciXSxbNCw1LCJcXENvbXBfXFxiIEZcXF57QScnfSJdLFs2LDUsIlxca2FwcGFcXF57Zmd9IiwxXSxbMSw3LCJGXFxeQWZeXFwqX1xcb1xceGlcXF5nIiwxXSxbNywzLCJcXGthcHBhXFxeZiBnXlxcKl9cXG8iLDFdLFs3LDYsIkZcXF5BXFxDb21wX1xcbyIsMV0sWzAsOCwiIiwyLHsibGV2ZWwiOjIsInN0eWxlIjp7ImhlYWQiOnsibmFtZSI6Im5vbmUifX19XSxbOCw1LCJGXFxeQVxceGlcXF57Zmd9IiwyXV0=
\begin{tikzcd}
{F\^AE} & {F\^Af^\*_\s E'} & {(Ff)^\*_\b F\^{A'}E'} & {(Ff)^\*_\b F\^{A'}g^\*_\s E''} \\
&& {F\^Af^\*_\s g^\*_\s E''} & {(Ff)^\*_\b(Fg)^\*_\b F\^{A''}E''} \\
&& {F\^A(g\o f)^\*_\s E''} \\
{F\^AE} &&& {(F(g\o f))^\*_\b F\^{A''}E''}
\arrow["{F\^A\xi\^f}", from=1-1, to=1-2]
\arrow[Rightarrow, no head, from=1-1, to=4-1]
\arrow["{\kappa\^f}", from=1-2, to=1-3]
\arrow["{F\^Af^\*_\s\xi\^g}"{description}, from=1-2, to=2-3]
\arrow["{(Ff)^\*_\b F\^{A'}\xi\^g}", from=1-3, to=1-4]
\arrow["{(Ff)^\*_\b\kappa\^g}", from=1-4, to=2-4]
\arrow["{\kappa\^f g^\*_\s}"{description}, from=2-3, to=1-4]
\arrow["{F\^A\Comp_\s}"{description}, from=2-3, to=3-3]
\arrow["{\Comp_\b F\^{A''}}", from=2-4, to=4-4]
\arrow["{\kappa\^{g\o f}}"{description}, from=3-3, to=4-4]
\arrow["{F\^A\xi\^{g\o f}}"', from=4-1, to=4-4]
\end{tikzcd}
\end{equation*}
Therefore:
\begin{eqnarray*}
&=&\left(F(g\o f),\xi\^{F(g\o f)}\right)\\
&=&\left(F(g\o f),\kappa\^{g\o f}\o F\^A\Comp_\s\o F\^Af^\*_\s\xi\^g\o F\^A\xi\^f\right)\\
&=&\left(Fg\o Ff,\Comp_{\b F\^{A''}}\o(Ff)^\*_\b\kappa\^g\o(Ff)^\*_\b F\^{A'}\xi\^g\o\kappa\^f\o F\^A\xi\^f\right)\\
&=&\left(Fg,\xi\^{Fg}\right)\o\left(Ff,\xi\^{Ff}\right)\\
&=&F'\left(g,\xi\^g\right)\o F'\left(f,\xi\^f\right)
\end{eqnarray*}
Let us now consider a $2$-morphism $(\theta,\theta')\colon(F,F',\kappa)\Rightarrow(G,G',\lambda)$ of indexed categories between two $1$-morphisms $(F,F',\kappa),(G,G',\lambda)\colon(\X_\s,\I_\s)\to(\X_\b,\I_\b)$ of indexed categories. Let us define the corresponding $2$-morphism of fibrations $\FF(\theta,\theta')\colon\FF(F,F',\kappa)\Rightarrow\FF(G,G',\lambda)$. First, the base natural transformation is given by $\theta\colon F\Rightarrow G$. The total natural transformation $\theta'\colon F'\to G'$ is the morphism so defined:
\begin{align*}
&(\theta'\colon F'(A,E)=(FA,F\^AE)\to(GA,G\^AE)=G'(A,E))=(\theta\colon FA\to GA,\theta\^A\colon F\^A E\to\theta^\*G\^AE)
\end{align*}
Let us prove that $\theta'$ is natural. Consider a morphism $(f,\xi\^f)\colon(A,E)\to(A',E')$ of $\El(\X_\s,\I_\s)$. Then
\begin{equation*}
% https://q.uiver.app/#q=WzAsNyxbMCwwLCJGXFxeQUUiXSxbMSwwLCJGXFxeQWZeXFwqX1xcbyBFJyJdLFswLDEsIlxcdGhldGFeXFwqIEdcXF5BRSJdLFsxLDEsIlxcdGhldGFeXFwqIEdcXF5BZl5cXCpfXFxvIEUnIl0sWzIsMCwiKEZmKV5cXCpfXFxiIEZcXF57QSd9RSciXSxbMywxLCJcXHRoZXRhXlxcKihHZileXFwqX1xcYiBHXFxee0EnfUUnIl0sWzMsMCwiKEZmKV5cXCpcXHRoZXRhXlxcKiBHXFxee0EnfUUnIl0sWzAsMiwiXFx0aGV0YVxcXkEiLDJdLFsyLDMsIlxcdGhldGFeXFwqIEdcXF5BXFx4aVxcXmYiLDJdLFswLDEsIkZcXF5BXFx4aVxcXmYiXSxbMSwzLCJcXHRoZXRhXFxeQWZeXFwqX1xcbyIsMV0sWzEsNCwiXFxrYXBwYVxcXmYiXSxbMyw1LCJcXHRoZXRhXlxcKlxcbGFtYmRhXFxeZiIsMl0sWzYsNSwiIiwwLHsibGV2ZWwiOjIsInN0eWxlIjp7ImhlYWQiOnsibmFtZSI6Im5vbmUifX19XSxbNCw2LCIoRmYpXlxcKl9cXGJcXHRoZXRhXFxee0EnfSJdXQ==
\begin{tikzcd}
{F\^AE} & {F\^Af^\*_\s E'} & {(Ff)^\*_\b F\^{A'}E'} & {(Ff)_\b^\*\theta^\* G\^{A'}E'} \\
{\theta^\* G\^AE} & {\theta^\* G\^Af^\*_\s E'} && {\theta^\*(Gf)^\*_\b G\^{A'}E'}
\arrow["{F\^A\xi\^f}", from=1-1, to=1-2]
\arrow["{\theta\^A}"', from=1-1, to=2-1]
\arrow["{\kappa\^f}", from=1-2, to=1-3]
\arrow["{\theta\^Af^\*_\s}"{description}, from=1-2, to=2-2]
\arrow["{(Ff)^\*_\b\theta\^{A'}}", from=1-3, to=1-4]
\arrow[Rightarrow, no head, from=1-4, to=2-4]
\arrow["{\theta^\* G\^A\xi\^f}"', from=2-1, to=2-2]
\arrow["{\theta^\*\lambda\^f}"', from=2-2, to=2-4]
\end{tikzcd}
\end{equation*}
where we used the compatibility of $\theta'$ with the distributors. This shows that $\II\colon\Fib\to\Indx$ and $\FF\colon\Indx\to\Fib$ are $2$-functors. We leave to the reader to show that they form an equivalence, since this is just an extension to the classical equivalence $\Fib(\X)\simeq\Indx(\X)$.
\end{proof}
\end{proposition}

%__________________________________________________________________________
\subsection{Tangent category theory}
\label{subsection:tangent-cats-background}
Tangent category theory aims to categorically axiomatize some of the fundamental constructions of differential geometry. In particular, a tangent category, first introduced by Rosick\`y in~\cite{rosicky:tangent-cats} and further generalized by Cockett and Cruttwell in~\cite{cockett:tangent-cats}, consists of a category equipped with an endofunctor which can be interpreted as the tangent bundle functor. An object in a tangent category can be interpreted as a geometric space equipped with a tangent bundle which encodes a notion of \textit{local linearity} for this space. Informally, the fibres of the tangent bundle can be interpreted as \textit{tangent spaces} of the given space, i.e., the best \textit{local linear} approximation of the space at a given point.
\par The notion of linearity in a tangent category is not an intrinsic notion, but rather a contextual notion: the tangent structure establishes what it means for a morphism to be linear. In~\cite{cockett:differential-bundles}, Cockett and Cruttwell introduced the concept of differential bundles in a tangent category, which are the analogs of vector bundles in the category of smooth manifolds, as shown by MacAdam in~\cite{macadam:vector-bundles}. Differential bundles represent precisely those fibre bundle-like objects whose fibres carry a linear structure.
\par In this section, we would like to recall the main definitions of a tangent category, of a differential bundle in a tangent category and discuss some examples.
\par We start with the definition of a tangent category. We suggest the reader to consult~\cite{cockett:tangent-cats}. The original definition of Rosick\`y presented in~\cite{rosicky:tangent-cats} is now known with the denomination of a \textit{tangent category with negatives}, sometimes also called a \textit{Rosick\`y tangent category} (see~\cite{cruttwell:algebraic-geometry} or~\cite{ikonicoff:operadic-algebras-tagent-cats}). First, let us recall the definition of an additive bundle.

\begin{definition}
\label{definition:additive-bundle}
An \textbf{additive bundle} of a category $\X$ is a commutative monoid internal to the slice category $\X/B$ for a given object $B$ of $\X$. Concretely, this consists of a morphism $q\colon E\to B$ of $\X$, called the \textbf{projection}, for which the $n$-fold pullback $E_n\to B$ of $q$ along itself exists, together with a morphism $z_q\colon B\to E$ of $q$, called the \textbf{zero morphism}, and a morphism $s_q\colon E_2\to E$, called the \textbf{sum morphism}, satisfying the following properties:
\begin{itemize}
\item $z_q$ is a section of $q$:
\begin{equation*}
% https://q.uiver.app/#q=WzAsMyxbMCwwLCJFIl0sWzEsMSwiQSJdLFswLDEsIkEiXSxbMiwwLCJ6X3EiXSxbMCwxLCJxIl0sWzIsMSwiIiwyLHsibGV2ZWwiOjIsInN0eWxlIjp7ImhlYWQiOnsibmFtZSI6Im5vbmUifX19XV0=
\begin{tikzcd}
E \\
B & B
\arrow["q", from=1-1, to=2-2]
\arrow["{z_q}", from=2-1, to=1-1]
\arrow[Rightarrow, no head, from=2-1, to=2-2]
\end{tikzcd}
\end{equation*}
\item $s_q$ is a bundle morphism:
\begin{equation*}
% https://q.uiver.app/#q=WzAsNCxbMSwwLCJFIl0sWzEsMSwiQSJdLFswLDAsIkVfMiJdLFswLDEsIkUiXSxbMCwxLCJxIl0sWzIsMywiXFxwaV8xIiwyXSxbMywxLCJxIiwyXSxbMiwwLCJzX3EiXV0=
\begin{tikzcd}
{E_2} & E \\
E & B
\arrow["{s_q}", from=1-1, to=1-2]
\arrow["{\pi_1}"', from=1-1, to=2-1]
\arrow["q", from=1-2, to=2-2]
\arrow["q"', from=2-1, to=2-2]
\end{tikzcd}\hfill
% https://q.uiver.app/#q=WzAsNCxbMSwwLCJFIl0sWzEsMSwiQSJdLFswLDAsIkVfMiJdLFswLDEsIkUiXSxbMCwxLCJxIl0sWzIsMywiXFxwaV8yIiwyXSxbMywxLCJxIiwyXSxbMiwwLCJzX3EiXV0=
\begin{tikzcd}
{E_2} & E \\
E & B
\arrow["{s_q}", from=1-1, to=1-2]
\arrow["{\pi_2}"', from=1-1, to=2-1]
\arrow["q", from=1-2, to=2-2]
\arrow["q"', from=2-1, to=2-2]
\end{tikzcd}
\end{equation*}
$\pi_k\colon E_n\to E$ being the $k$-th projection of the $n$-fold pullback;
\item Associativity:
\begin{equation*}
% https://q.uiver.app/#q=WzAsNCxbMCwwLCJFXzMiXSxbMSwwLCJFXzIiXSxbMSwxLCJFIl0sWzAsMSwiRV8yIl0sWzAsMSwiXFxpZF9FXFx0aW1lc19Bc19xIl0sWzEsMiwic19xIl0sWzMsMiwic19xIiwyXSxbMCwzLCJzX3FcXHRpbWVzX0FcXGlkX0UiLDJdXQ==
\begin{tikzcd}
{E_3} & {E_2} \\
{E_2} & E
\arrow["{\id_E\times_As_q}", from=1-1, to=1-2]
\arrow["{s_q\times_A\id_E}"', from=1-1, to=2-1]
\arrow["{s_q}", from=1-2, to=2-2]
\arrow["{s_q}"', from=2-1, to=2-2]
\end{tikzcd}
\end{equation*}
\item Unitality:
\begin{equation*}
% https://q.uiver.app/#q=WzAsMyxbMCwwLCJFXzIiXSxbMCwxLCJFIl0sWzEsMCwiRSJdLFsxLDAsIlxcPHF6X3EsXFxpZF9FXFw+Il0sWzAsMiwic19xIl0sWzEsMiwiIiwyLHsibGV2ZWwiOjIsInN0eWxlIjp7ImhlYWQiOnsibmFtZSI6Im5vbmUifX19XV0=
\begin{tikzcd}
{E_2} & E \\
E
\arrow["{s_q}", from=1-1, to=1-2]
\arrow["{\<qz_q,\id_E\>}", from=2-1, to=1-1]
\arrow[Rightarrow, no head, from=2-1, to=1-2]
\end{tikzcd}
\end{equation*}
\item Commutativity:
\begin{equation*}
% https://q.uiver.app/#q=WzAsNCxbMCwwLCJFXzIiXSxbMCwxLCJFXzIiXSxbMSwwLCJFIl0sWzEsMSwiRSJdLFsyLDMsIiIsMCx7ImxldmVsIjoyLCJzdHlsZSI6eyJoZWFkIjp7Im5hbWUiOiJub25lIn19fV0sWzAsMiwic19xIl0sWzEsMywic19xIiwyXSxbMCwxLCJcXHRhdSIsMl1d
\begin{tikzcd}
{E_2} & E \\
{E_2} & E
\arrow["{s_q}", from=1-1, to=1-2]
\arrow["\tau"', from=1-1, to=2-1]
\arrow[Rightarrow, no head, from=1-2, to=2-2]
\arrow["{s_q}"', from=2-1, to=2-2]
\end{tikzcd}
\end{equation*}
where $\tau\colon E_2\to E_2$ denotes the flip $\<\pi_2,\pi_1\>$.
\end{itemize}
An \textbf{additive bundle morphism} consists of a pair $(f,g)\colon(B,E,q,z_q,s_q)\to(B',E',q',z_q',s_q')$ of morphisms $f\colon B\to B'$ and $g\colon E\to E'$ satisfying the following properties:
\begin{itemize}
\item Compatibility with the projections:
\begin{equation*}
% https://q.uiver.app/#q=WzAsNCxbMCwwLCJFIl0sWzEsMCwiRSciXSxbMCwxLCJBIl0sWzEsMSwiQSciXSxbMCwxLCJnIl0sWzAsMiwicSIsMl0sWzEsMywicSciXSxbMiwzLCJmIiwyXV0=
\begin{tikzcd}
E & {E'} \\
B & {B'}
\arrow["g", from=1-1, to=1-2]
\arrow["q"', from=1-1, to=2-1]
\arrow["{q'}", from=1-2, to=2-2]
\arrow["f"', from=2-1, to=2-2]
\end{tikzcd}
\end{equation*}
\item Compatibility with the zero morphisms:
\begin{equation*}
% https://q.uiver.app/#q=WzAsNCxbMCwwLCJFIl0sWzEsMCwiRSciXSxbMCwxLCJBIl0sWzEsMSwiQSciXSxbMCwxLCJnIl0sWzIsMCwiel9xIl0sWzMsMSwiel9xJyIsMl0sWzIsMywiZiIsMl1d
\begin{tikzcd}
E & {E'} \\
B & {B'}
\arrow["g", from=1-1, to=1-2]
\arrow["{z_q}", from=2-1, to=1-1]
\arrow["f"', from=2-1, to=2-2]
\arrow["{z_q'}"', from=2-2, to=1-2]
\end{tikzcd}
\end{equation*}
\item Additivity:
\begin{equation*}
% https://q.uiver.app/#q=WzAsNCxbMCwxLCJFIl0sWzEsMSwiRSciXSxbMCwwLCJFXzIiXSxbMSwwLCJFJ18yIl0sWzAsMSwiZyIsMl0sWzIsMywiZ1xcdGltZXNfQWciXSxbMiwwLCJzX3EiLDJdLFszLDEsInNfcSciXV0=
\begin{tikzcd}
{E_2} & {E'_2} \\
E & {E'}
\arrow["{g\times_Bg}", from=1-1, to=1-2]
\arrow["{s_q}"', from=1-1, to=2-1]
\arrow["{s_q'}", from=1-2, to=2-2]
\arrow["g"', from=2-1, to=2-2]
\end{tikzcd}
\end{equation*}
\end{itemize}
\end{definition}

\begin{definition}
\label{definition:tangent-category}
A \textbf{tangent category} is a category $\X$ equipped with a \textbf{tangent structure}, which consists of the following data:
\begin{description}
\item[tangent bundle functor] An endofunctor $\T\colon\X\to\X$, called the \textbf{tangent bundle functor};

\item[projection] A natural transformation $p\colon\T\Rightarrow\id_\X$, called the \textbf{projection}, for which the $n$-fold pullbacks $\T_nA\to A$ of $p\colon\T A\to A$ along itself exist and preserved by the functors $\T^m\=\T\o\T\o{\dots}\o\T$, for every positive integer $m$ (for $m=1$, $\T^1\=\T$). The $k$-th projection is denoted by $\pi_k\colon\T_n\Rightarrow\T$;

\item[zero morphism] A natural transformation $z\colon\id_\X\Rightarrow\T$, called the \textbf{zero morphism};

\item[sum morphism] A natural transformation $s\colon\T_2\Rightarrow\T$, called the \textbf{sum morphism};

\item[vertical lift] A natural transformation $l\colon\T\Rightarrow\T^2$, called the \textbf{vertical lift};

\item[canonical flip] A natural transformation $c\colon\T^2\Rightarrow\T^2$, called the \textbf{canonical flip};
\end{description}
satisfying the following properties:
\begin{itemize}
\item Additive structure. For each object $A\in\X$, $(A,\T A,p,z,s)$ is an additive bundle;

\item Additivity of vertical lift. For each object $A\in\X$
\begin{align*}
&(z,l)\colon(A,\T A,p,z,s)\to(\T A,\T^2A,\T p,\T z,\T s)
\end{align*}
is an additive bundle morphism;

\item Additivity of canonical flip. For each object $A\in\X$
\begin{align*}
&(\id_{\T A},c)\colon(\T A,\T^2A,\T p,\T z,\T s)\to(\T A,\T^2A,p_\T,z_\T,s_\T)
\end{align*}
is an additive bundle morphism;

\item Coassociativity of vertical lift:
\begin{equation*}
% https://q.uiver.app/#q=WzAsNCxbMSwxLCJcXFReM0EiXSxbMCwwLCJcXFQgQSJdLFsxLDAsIlxcVF4yQSJdLFswLDEsIlxcVF4yQSJdLFsxLDMsImwiLDJdLFsxLDIsImwiXSxbMiwwLCJcXFQgbCJdLFszLDAsImxcXFQiLDJdXQ==
\begin{tikzcd}
{\T A} & {\T^2A} \\
{\T^2A} & {\T^3A}
\arrow["l", from=1-1, to=1-2]
\arrow["l"', from=1-1, to=2-1]
\arrow["{\T l}", from=1-2, to=2-2]
\arrow["l\T"', from=2-1, to=2-2]
\end{tikzcd}
\end{equation*}

\item The canonical flip is a symmetric braiding:
\begin{equation*}
% https://q.uiver.app/#q=WzAsMyxbMCwwLCJcXFReMkEiXSxbMSwwLCJcXFReMkEiXSxbMCwxLCJcXFReMkEiXSxbMiwwLCJjIl0sWzAsMSwiYyJdLFsyLDEsIiIsMix7ImxldmVsIjoyLCJzdHlsZSI6eyJoZWFkIjp7Im5hbWUiOiJub25lIn19fV1d
\begin{tikzcd}
{\T^2A} & {\T^2A} \\
{\T^2A}
\arrow["c", from=1-1, to=1-2]
\arrow["c", from=2-1, to=1-1]
\arrow[Rightarrow, no head, from=2-1, to=1-2]
\end{tikzcd}\hfill
% https://q.uiver.app/#q=WzAsNixbMCwwLCJcXFReM0EiXSxbMSwwLCJcXFReM0EiXSxbMiwwLCJcXFReM0EiXSxbMiwxLCJcXFReM0EiXSxbMSwxLCJcXFReM0EiXSxbMCwxLCJcXFReM0EiXSxbMCwxLCJcXFQgYyJdLFs0LDMsImNcXFQiLDJdLFs1LDQsIlxcVCBjIiwyXSxbMSwyLCJjXFxUIl0sWzAsNSwiY1xcVCIsMl0sWzIsMywiXFxUIGMiXV0=
\begin{tikzcd}
{\T^3A} & {\T^3A} & {\T^3A} \\
{\T^3A} & {\T^3A} & {\T^3A}
\arrow["{\T c}", from=1-1, to=1-2]
\arrow["c\T"', from=1-1, to=2-1]
\arrow["c\T", from=1-2, to=1-3]
\arrow["{\T c}", from=1-3, to=2-3]
\arrow["{\T c}"', from=2-1, to=2-2]
\arrow["c\T"', from=2-2, to=2-3]
\end{tikzcd}
\end{equation*}

\item Compatibility between canonical flip and vertical lift:
\begin{equation*}
% https://q.uiver.app/#q=WzAsMyxbMCwwLCJcXFReMkEiXSxbMCwxLCJcXFQgQSJdLFsxLDAsIlxcVF4yQSJdLFsxLDAsImwiXSxbMCwyLCJjIl0sWzEsMiwibCIsMl1d
\begin{tikzcd}
{\T^2A} & {\T^2A} \\
{\T A}
\arrow["c", from=1-1, to=1-2]
\arrow["l", from=2-1, to=1-1]
\arrow["l"', from=2-1, to=1-2]
\end{tikzcd}\hfill
% https://q.uiver.app/#q=WzAsNSxbMSwwLCJcXFReM0EiXSxbMiwwLCJcXFReM0EiXSxbMCwwLCJcXFReMkEiXSxbMiwxLCJcXFReM0EiXSxbMCwxLCJcXFReMkEiXSxbMCwxLCJcXFQgYyJdLFsyLDAsImxcXFQiXSxbMSwzLCJjXFxUIl0sWzIsNCwiYyIsMl0sWzQsMywiXFxUIGwiLDJdXQ==
\begin{tikzcd}
{\T^2A} & {\T^3A} & {\T^3A} \\
{\T^2A} && {\T^3A}
\arrow["l\T", from=1-1, to=1-2]
\arrow["c"', from=1-1, to=2-1]
\arrow["{\T c}", from=1-2, to=1-3]
\arrow["c\T", from=1-3, to=2-3]
\arrow["{\T l}"', from=2-1, to=2-3]
\end{tikzcd}
\end{equation*}

\item Universality of the vertical lift. The following diagram:
\begin{equation*}
% https://q.uiver.app/#q=WzAsNCxbMSwxLCJcXFQgQSJdLFswLDAsIlxcVF8yQSJdLFsxLDAsIlxcVF4yQSJdLFswLDEsIkEiXSxbMSwzLCJcXHBpXzFwIiwyXSxbMSwyLCJcXHBoaSJdLFsyLDAsIlxcVCBwIl0sWzMsMCwieiIsMl1d
\begin{tikzcd}
{\T_2A} & {\T^2A} \\
A & {\T A}
\arrow["\phi", from=1-1, to=1-2]
\arrow["{\pi_1p}"', from=1-1, to=2-1]
\arrow["{\T p}", from=1-2, to=2-2]
\arrow["z"', from=2-1, to=2-2]
\end{tikzcd}
\end{equation*}
is a pullback diagram, where:
\begin{align*}
&\phi\colon\T_2A\xrightarrow{l\times_Az_\T}\T\T_2A\xrightarrow{\T s}\T^2A
\end{align*}
\end{itemize}
A tangent category \textbf{has negatives} when it comes equipped with:
\begin{description}
\item[negation] A natural transformation $n\colon\T\Rightarrow\T$, called \textbf{negation}, such that:
\begin{equation*}
% https://q.uiver.app/#q=WzAsNCxbMSwwLCJcXFRfMkEiXSxbMCwwLCJcXFQgQSJdLFsxLDEsIlxcVCBBIl0sWzAsMSwiQSJdLFsxLDAsIlxcPFxcaWQsblxcPiJdLFswLDIsInMiXSxbMSwzLCJwIiwyXSxbMywyLCJ6IiwyXV0=
\begin{tikzcd}
{\T A} & {\T_2A} \\
A & {\T A}
\arrow["{\<\id,n\>}", from=1-1, to=1-2]
\arrow["p"', from=1-1, to=2-1]
\arrow["s", from=1-2, to=2-2]
\arrow["z"', from=2-1, to=2-2]
\end{tikzcd}
\end{equation*}
\end{description}
\end{definition}

The tangent bundle functor sends an object $A$ to its tangent bundle $\T A$, which is interpreted as the space of tangent vectors of $A$ at each point. Moreover, it sends a morphism $f\colon A\to A'$, regarded as a smooth map, to the morphism $\T f$ which sends each tangent vector $u$ at a given point $x$ to the differential of $f$ at $x$ applied to $u$, i.e., $\d f_x(u)$. The functoriality of $\T$ corresponds to the chain rule of the differential.
\par The projection, the zero morphism, and the sum morphism equip each tangent bundle $\T A\to A$ with an additive structure. Informally, this structure axiomatizes the additive structure of the tangent spaces of a tangent bundle of a smooth manifold. By interpreting $\T A$ as the space of vector tangent vectors over each point of $A$, the projection sends each tangent vector to its base point, the zero morphism associates each point with the zero tangent vector, and the sum morphism allows one to sum together two tangent vectors with the same base point.
\par When the tangent category has negatives, the additive structure over the fibres becomes a commutative group. This can be interpreted as inverting the orientation of the tangent vectors.
\par The double tangent bundle $\T^2A$ can be regarded as the space of infinitesimal differentiable homotopies of $A$, i.e., tangent vectors of tangent vectors. Locally, one can represent such an element as a quadruple $(x,u,v,\omega)$ formed by a point $x$ of $A$, two tangent vectors $u,v$ of $A$ at $x$, and a tangent vector $\omega$ of $\T A$ at $(x,u)$. Then, the vertical lift sends a tangent vector $u$ of $A$ at a given point $x$ to $(x,0,0,u)$.
\par The universality of the vertical lift establishes that the vertical bundle, which is the bundle whose fibres are the kernel of the differential of the projection at each point, is trivial. This can be interpreted as a local linearity condition for the tangent bundle, that is, an object in a tangent category has a local linear behaviour.
\par Finally, the canonical flip encodes the symmetry of the Hessian matrix, i.e., the commutativity of partial derivatives for smooth functions.

\begin{notation}
\label{notation:tangent-category}
For a tangent category $(\X,\TT)$, we denote the tangent bundle functor using the same letter of the tangent structure $\TT$, adopting the font $\T\colon\X\to\X$. The $n$-fold pullback of the projection is denoted by $\T_n$ with the corresponding projections denoted by $\pi_k\colon\T_n\Rightarrow\T$. The projection, the zero morphism, the sum morphism, the vertical lift, and the canonical flip are denoted with the letters $p,z,s,l$, and $c$, respectively; when $(\X,\TT)$ has negatives, the negation is denoted by $n$. When the symbol of the tangent structure is decorated with subscripts or subscripts, the same decoration applies to $\T,\T_n,p,z,s,l,c$, and $n$ whenever the tangent structure has negatives.
\end{notation}

\begin{example}
\label{example:trivial-tangent-cat}
Every category $\X$ admits a trivial tangent structure, whose tangent bundle functor and whose natural transformations are the identities.
\end{example}

\begin{example}
\label{example:differential-geometry}
The category of finite-dimensional smooth manifolds forms a tangent category with negatives whose tangent bundle functor is the usual notion of the tangent bundle functor in differential geometry.
\end{example}

\begin{example}
\label{example:commutative-rigs}
The category of commutative and unital rigs (i.e., rings without negatives, a.k.a. semirings) $\cRig$ has a canonical tangent category whose tangent bundle functor sends a rig $R$ to the rig $R[\epsilon]\=R[x]/(x^2)$; the projection $p\colon R[\epsilon]\to R$ sends $a\in R$ to itself; the zero morphism $z\colon R\to R[\epsilon]$ sends $\epsilon$ to $0$; the sum $s\colon R[\epsilon_1,\epsilon_2]\to R[\epsilon]$ sends both $\epsilon_1$ and $\epsilon_2$ to $\epsilon$, where $R[\epsilon_1,\epsilon_2]\=R[x_1,x_2]/(x_ix_j;i,j=1,2)$; the vertical lift $l\colon R[\epsilon]\to R[\epsilon][\epsilon']$ sends $\epsilon$ to $\epsilon'\epsilon$, where $R[\epsilon][\epsilon']\= R[x,x']/(x^2,{x'}^2)$; finally, the canonical flip $c\colon R[\epsilon][\epsilon']\to R[\epsilon][\epsilon']$ sends $\epsilon$ to $\epsilon'$ and vice versa.
\par A similar tangent structure is also defined in the category $\cRing$ of commutative and unital rings. This tangent structure has negatives, whose negation $n\colon R[\epsilon]\to R[\epsilon]$ sends $\epsilon$ to $-\epsilon$.
\end{example}

\begin{example}
\label{example:algebraic-geometry}
The opposite of the category of commutative and unital rigs comes equipped with a tangent structure, which is adjoint to the tangent structure of Example~\ref{example:commutative-rigs}, i.e., the tangent bundle functors form an adjunction. Concretely, the tangent bundle functor of $\cRig^\op$ sends a rig $R$ to the symmetric rig of the module of K\"ahler differentials $\T R\=\Sym_R\Omega^1_{R/\N}$ of $R$, i.e., to the rig generated by all elements $a$ of $R$ together with symbols $\d a$, for each $a$ of $R$, satisfying the following relations:
\begin{align*}
&a\cdot_{\T R}b=a\cdot_Rb\\
&\d(1)=0\\
&\d(ra+sb)=r\d a+s\d b\\
&\d(ab)=b\d a+a\d b
\end{align*}
Regarding the natural transformations are morphisms of $\cRig$, the projection $p\colon R\to\T R$ sends each $a\in R$ to itself; the zero morphism $\T R\to R$ sends each $a$ to itself and each $\d a$ to $0$; the sum morphism $s\colon\T R\to\T_2R$ sends each $a$ to itself and each $\d a$ to $\d a\otimes 1+1\otimes\d a$, where $\T_2R=\T R\otimes_R\T R$; the vertical lift $l\colon\T^2R\to\T R$ sends each $a$ to itself, each $\d a$ and each $\d' a$ to $0$, and $\d'\d a$ to $\d a$, where $\d'x$ is the K\"ahler differential of $x\in\T R$; finally, the canonical flip $c\colon\T^2R\to\T^2R$ sends each $a$ to itself, $\d a$ to $\d'a$, $\d'a$ to $\d a$, and $\d'\d a$ to itself.
\par This tangent structure restricts to the subcategory $\cRing^\op$. Moreover, $\cRing^\op$ has negatives whose negation $n\colon\T R\to\T R$ sends each $a$ to itself and each $\d a$ to $-\d a$.
\end{example}

A morphism of tangent categories consists of a functor between the underlying categories, together with a distributive law between the functor and the respective tangent bundle functors. The distributive law can be lax, colax, an isomorphism, or the identity. Consequently, there are four different flavours of morphisms of tangent categories.

\begin{definition}
\label{definition:morphisms-tangent-cats}
Given two tangent categories $(\X,\TT)$ and $(\X',\TT')$, a \textbf{lax tangent morphism} consists of a pair $(F,\alpha)\colon(\X,\TT)\to(\X',\TT')$ formed by a functor $F\colon\X\to\X'$ together with a natural transformation $\alpha\colon F\o\T\Rightarrow\T'\o F$, called the \textbf{lax distributive law}, compatible with the tangent structures as follows:
\begin{itemize}
\item Additivity. For every $A\in\X$
\begin{align*}
&(\id_{FA},\alpha)\colon(FA,F\T A,Fp,Fz,Fs)\to(FA,\T'FA,p'_F,z'_F,s'_F)
\end{align*}
is an additive bundle morphism;

\item Compatibility with the vertical lifts:
\begin{equation*}
% https://q.uiver.app/#q=WzAsNSxbMCwxLCJGXFxUIEEiXSxbMiwxLCJcXFQnRkEiXSxbMCwwLCJGXFxUXjJBIl0sWzEsMCwiXFxUJ0ZcXFQgQSJdLFsyLDAsIntcXFQnfV4yRkEiXSxbMCwxLCJcXGFscGhhIiwyXSxbMiwzLCJcXGFscGhhXFxUIl0sWzMsNCwiXFxUJ1xcYWxwaGEiXSxbMCwyLCJGbCJdLFsxLDQsImwnRiIsMl1d
\begin{tikzcd}
{F\T^2A} & {\T'F\T A} & {{\T'}^2FA} \\
{F\T A} && {\T'FA}
\arrow["\alpha\T", from=1-1, to=1-2]
\arrow["{\T'\alpha}", from=1-2, to=1-3]
\arrow["Fl", from=2-1, to=1-1]
\arrow["\alpha"', from=2-1, to=2-3]
\arrow["{l'F}"', from=2-3, to=1-3]
\end{tikzcd}
\end{equation*}

\item Compatibility with the canonical flips:
\begin{equation*}
% https://q.uiver.app/#q=WzAsNixbMCwxLCJGXFxUXjJBIl0sWzIsMSwie1xcVCd9XjJGQSJdLFswLDAsIkZcXFReMkEiXSxbMSwwLCJcXFQnRlxcVCBBIl0sWzIsMCwie1xcVCd9XjJGQSJdLFsxLDEsIlxcVCdGXFxUIEEiXSxbMiwzLCJcXGFscGhhXFxUIl0sWzMsNCwiXFxUJ1xcYWxwaGEiXSxbMiwwLCJGYyIsMl0sWzQsMSwiYydGIl0sWzAsNSwiXFxhbHBoYVxcVCIsMl0sWzUsMSwiXFxUJ1xcYWxwaGEiLDJdXQ==
\begin{tikzcd}
{F\T^2A} & {\T'F\T A} & {{\T'}^2FA} \\
{F\T^2A} & {\T'F\T A} & {{\T'}^2FA}
\arrow["\alpha\T", from=1-1, to=1-2]
\arrow["Fc"', from=1-1, to=2-1]
\arrow["{\T'\alpha}", from=1-2, to=1-3]
\arrow["{c'F}", from=1-3, to=2-3]
\arrow["\alpha\T"', from=2-1, to=2-2]
\arrow["{\T'\alpha}"', from=2-2, to=2-3]
\end{tikzcd}
\end{equation*}
\end{itemize}
A \textbf{colax tangent morphism} consists of a pair $(G,\beta)\colon(\X,\TT)\nto(\X',\TT')$ formed by a functor $G\colon\X\to\X'$ together with a natural transformation $\beta\colon\T'\o G\Rightarrow G\o\T$ which fulfills the duals of the compatibility condition of a lax distributive law. Moreover, the functor $G$ preserves each $n$-fold pullback of the projection along itself and the pullback diagram which establishes the universal property of the vertical lift.\\
Finally, a colax tangent morphism whose distributive law is invertible is called a \textbf{strong tangent morphism} and is a \textbf{strict tangent morphism} when the distributive law is the identity. For colax tangent morphisms we employ the following notation $(G,\beta)\colon(\X,\TT)\nto(\X',\TT')$. Moreover, we identify a strict tangent morphism with its underlying functor. We also adopt the convention to say that a functor between two tangent categories \textbf{strongly} (or \textbf{strictly}) \textbf{preserves} the tangent structures, when such a functor comes equipped with a distributive law which makes it into a strong strict) tangent morphism.
\end{definition}

%__________________________________________________________________________
%__________________________________________________________________________

\section{Tangent fibrations and indexed tangent categories}
\label{section:tangent-fibrations}
\par Cockett and Cruttwell in~\cite[Section~5]{cockett:differential-bundles} observed that differential bundles of a tangent category, which are the analogs of vector bundles in the context of differential geometry (see~\cite{macadam:vector-bundles}), can be naturally organized into a fibration. Such a fibration is also compatible with the tangent structure of differential bundles and with the one of the base category. This led to the introduction of the notion of a \textbf{tangent fibration}. Cockett and Cruttwell also showed that the fibres of a tangent fibration carry a canonical tangent structure. 
\par We start this section by recalling this definition and this result. We then discuss a first attempt to extend the Grothendieck construction to tangent fibrations using Cockett and Cruttwell's result. We will show that their procedure of ``splitting'' the total tangent category into its fibres inevitably loses information about the total tangent structure. We end this section by comparing the notion of internal fibrations in the $2$-category of tangent categories, as introduced by Street in~\cite{street:fibrations}, with the notion of tangent fibrations.

%__________________________________________________________________________
\subsection{Tangent fibrations and their fibres}
\label{subsection:tangent-fibrations-definition}
Informally, a tangent fibration is a fibration between two tangent categories, compatible with the tangent structures. To introduce the required compatibility condition, consider first a cloven fibration $\Pi\colon\X'\to\X$ between two categories, each equipped with an endofunctor $\T'\colon\X'\to\X'$ and $\T\colon\X\to\X$, respectively; suppose also that $\Pi\o\T'=\T\o\Pi$. 
\par On one hand, consider now a morphism $f\colon A\to A'$ of $\X$ and an object $E'$ over $A'$. Since $\Pi\T'E'=\T\Pi E'=\T A$, there is a cartesian lift $\varphi\^{\T f}\colon(\T f)^\*(\T' E')\to\T'E'$ of $\T f$ on $\T'E'$. On the other hand, one can also apply the endofunctor $\T'$ to the cartesian lift $\varphi\^f\colon f^\*E'\to E'$ of $f$ on $E'$ and obtain $\T'\varphi\^f\colon\T'(f^\*E')\to\T'E'$.
\par Both $\varphi\^{\T f}$ and $\T'\varphi\^f$ are lifts of $\T f$, since $\Pi(\varphi\^{\T f})=\Pi(\T'\varphi\^f)=f$. By the universality of the cartesian lift $\varphi\^{\T f}$, there is a unique morphism $\kappa\^f\colon\T'(f^\*E')\dashrightarrow(\T f)^\*(\T'E')$ making the following diagram commutes:
\begin{equation*}
% https://q.uiver.app/#q=WzAsMyxbMCwwLCJcXFQnKGZeXFwqRScpIl0sWzAsMSwiKFxcVCBmKV5cXCooXFxUJ0UnKSJdLFsxLDEsIlxcVCdFJyJdLFsxLDIsIlxcdmFycGhpXFxee1xcVCBmfSIsMl0sWzAsMiwiXFxUJ1xcdmFycGhpXFxeZiJdLFswLDEsIlxceGlcXF5mIiwyLHsic3R5bGUiOnsiYm9keSI6eyJuYW1lIjoiZGFzaGVkIn19fV1d
\begin{tikzcd}
{\T'(f^\*E')} \\
{(\T f)^\*(\T'E')} & {\T'E'}
\arrow["{\kappa\^f}"', dashed, from=1-1, to=2-1]
\arrow["{\T'\varphi\^f}", from=1-1, to=2-2]
\arrow["{\varphi\^{\T f}}"', from=2-1, to=2-2]
\end{tikzcd}
\end{equation*}
Requiring the morphism $\kappa\^f$ to be an isomorphism is to ask the tangent bundle functors $(\T,\T')$ to form a $1$-morphism of fibrations, as stated in Definition~\ref{definition:1-morphisms-fibrations}. We follow the convention of referring to the morphisms $\kappa\^f$ as the \textbf{distributors} of the tangent bundle functors.
\par The original definition of a tangent fibration was first introduced by Cockett and Cruttwell in~\cite[Definition~5.2]{cockett:differential-bundles}. Here we specialize this definition for cloven fibrations.

\begin{definition}
\label{definition:tangent-fibration}
A \textbf{cloven tangent fibration} $\Pi\colon(\X',\TT')\to(\X,\TT)$ between two tangent categories $(\X',\TT')$ and $(\X,\TT)$, respectively called the \textbf{total} and the \textbf{base} tangent category, consists of a cloven fibration $\Pi\colon\X'\to\X$ whose underlying functor is a strict tangent morphism and for which the tangent bundle functors $\T'$ and $\T$ preserve the cartesian lifts, i.e., $(\T,\T')\colon\Pi\to\Pi$ is a $1$-morphism of fibrations.
\end{definition}

\begin{remark}
\label{remark:convention-on-cloven}
In the following, we assume every tangent fibration to be a cloven tangent fibration.
\end{remark}

In~\cite[Section~5]{cockett:differential-bundles}, Cockett and Cruttwell constructed two important examples of tangent fibrations: the tangent fibrations of display bundles and the one of differential bundles. In the next examples, we recall these two constructions. First, we recall a technical concept, introduced by Cockett and Cruttwell in the same paper to ensure the existence of the cartesian lifts. MacAdam reformulated this notion in~\cite[Definition~1.2.1]{macadam:vector-bundles}. Here, we adopt MacAdam's version.

\begin{definition}
\label{definition:tangent-display-system}
A \textbf{tangent display system} for a tangent category $(\X,\TT)$ consists of a family $\Dsply$ of morphisms of $\X$ for which the following two conditions hold:
\begin{itemize}
\item Stability under $\T$-pullbacks. The pullback of a morphism in $\Dsply$ along any other morphism exists, its universality is preserved by all functors $\T^m$, and is contained in $\Dsply$;

\item Stability under the tangent bundle functor. The tangent bundle functor sends morphisms of $\Dsply$ to morphisms of $\Dsply$.
\end{itemize}
In the following, when the tangent display system $\Dsply$ is clear from the context, we refer to the morphisms of $\Dsply$ as \textbf{tangent display maps}.
\end{definition}

In the tangent category of smooth manifolds (see Example~\ref{example:differential-geometry}), submersions constitute a canonical choice of a tangent display system; in the tangent category of affine schemes (see Example~\ref{example:algebraic-geometry}), one can consider the family of all morphisms as a tangent display system, since the category is finitely complete and the tangent bundle functor is the right adjoint to the (opposite of the) tangent bundle functor of the tangent category of Example~\ref{example:commutative-rigs}.

\begin{remark}
\label{remark:tangent-display-maps}
Recently, Cruttwell and the author of this paper found a characterization for the morphisms of the maximal tangent display system of any given tangent category. Such a characterization has not been published yet, however, it has already been presented at~\cite{cruttwell:ATCAT2024-connections},~\cite{cruttwell:CT2024-connections}, and~\cite{lanfranchi:MFCS2024-tangent-display-maps}. This notion will appear in a future paper.
\end{remark}

\begin{example}
\label{example:bundle-tangent-fibration}
Consider a tangent category $(\X,\TT)$ equipped with a tangent display system $\Dsply$ and let $\Dsply(\X,\TT)$ denote the category whose objects are tangent display maps and morphisms are commutative squares between them. It is straightforward to see that $\Dsply(\X,\TT)$ inherits a tangent structure whose tangent bundle functor is the tangent bundle functor $\T$ of $(\X,\TT)$ applied to tangent display maps. The functor $\Pi\colon\Dsply(\X,\TT)\to\X$ which sends each tangent display map $q\colon E\to A$ to its codomain $A$ and each commutative square $(f,g)\colon(q\colon E\to A)\to(q'\colon E'\to A')$ to the base morphism $f\colon A\to A'$, defines a tangent fibration $\Pi\colon\Dsply(\X,\TT)\to(\X,\TT)$ (cf.~\cite[Proposition~5.7]{cockett:differential-bundles}).
\end{example}

\begin{example}
\label{example:differential-bundle-tangent-fibration}
Differential bundles, introduced by Cockett and Cruttwell in~\cite{cockett:differential-bundles} to capture the analog of a vector bundle in a tangent category (cf.~\cite{macadam:vector-bundles}), can be organized in a tangent fibration. To see this, consider a tangent category $(\X,\TT)$ equipped with a tangent display system $\Dsply$. In Example~\ref{example:bundle-tangent-fibration} we discussed that the codomain functor $\Pi\colon\Dsply(\X,\TT)\to(\X,\TT)$ defines a tangent fibration. By restricting $\Pi$ to the full subcategory $\DBnd(\X,\TT)$ of $\Dsply(\X,\TT)$ of differential bundles which are also tangent display maps one obtains a new tangent fibration $\Pi\colon\DBnd(\X,\TT)\to(\X,\TT)$. Moreover, by restricting $\DBnd(\X,\TT)$ to only linear morphisms, i.e., morphisms compatible with the vertical lifts (see~\cite[Definition~2.3]{cockett:differential-bundles}), one also obtains another tangent fibration $\Pi\colon\DBnd_\lnr(\X,\TT)\to(\X,\TT)$. We invite the reader to consult~\cite[Section~5]{cockett:differential-bundles} for details on these constructions.
\end{example}

Cockett and Cruttwell proved in~\cite[Theorem~5.3]{cockett:differential-bundles} that the fibres of a cloven tangent fibration inherit a tangent structure from the total tangent category, strongly preserved by the substitution functors.
\par Let us briefly recall this construction; consider a tangent fibration $\Pi\colon(\X',\TT')\to(\X,\TT)$ and an object $A$ of $\X$. The first step is to show that each fibre comes equipped with a tangent structure $\TT\^A$ so defined:

\begin{description}
\item[tangent bundle functor] The tangent bundle functor $\T\^A\colon\Pi^{-1}(A)\to\Pi^{-1}(A)$ is the functor obtained by composing $\T'\colon\Pi^{-1}(A)\to\Pi^{-1}(\T A)$ with the substitution functor $z^\*\colon\Pi\^{-1}(\T A)\to\Pi^{-1}(A)$, induced by the zero morphism $z\colon A\to\T A$;

\item[projection] The projection $p\^A\colon\T\^A\Rightarrow\id_{\Pi^{-1}A}$ is defined by
\begin{align*}
&p\^A\colon z^\*\T'E\xrightarrow{\varphi\^z}\T'E\xrightarrow{p'}E
\end{align*}
where $\varphi\^z$ is the cartesian lift of the zero morphism of $(\X,\TT)$ and $p'$ denotes the projection of $(\X',\TT')$;

\item[zero morphism] The zero morphism $z\^A\colon\id_{\Pi^{-1}A}\Rightarrow\T\^A$ is the unique morphism $z\^A\colon E\to z^\*\T'E$ induced by the universality of the cartesian lift $\varphi\^z$, which makes the following diagram commutes:
\begin{equation*}
\begin{tikzcd}
E && {z^\*\T'E} \\
& {\T'E}
\arrow["{z'}"', from=1-1, to=2-2]
\arrow["{\varphi\^z}", from=1-3, to=2-2]
\arrow["{z\^A}", dashed, from=1-1, to=1-3]
\end{tikzcd}
\end{equation*}

\item[$n$-fold pullback] The $n$-fold pullback $\T\^A_n$ of the projection $p\^A$ along itself is the functor obtained by composing $\T'_n\colon\Pi^{-1}(A)\to\Pi^{-1}(\T_nA)$ and the substitution functor $z_n^\*\colon\Pi^{-1}(\T_nA)\to\Pi^{-1}(A)$, where $z_n\colon A\to\T_nA$ denotes $\<z,\dots,z\>$. The $k$-th projection $\pi_k\^A\colon\T\^A_n\Rightarrow\T\^A$ is the unique morphism induced by the universality of the cartesian lift $\varphi\^z$, which makes the following diagram commutes:
\begin{equation*}
\begin{tikzcd}
{z_n^\*\T'_nE} & {z^\*\T'E} \\
{\T'_nE} & {\T'E}
\arrow["{\varphi\^z}", from=1-2, to=2-2]
\arrow["{\pi_k}"', from=2-1, to=2-2]
\arrow["{\varphi\^{z_n}}"', from=1-1, to=2-1]
\arrow["{\pi_k\^A}", dashed, from=1-1, to=1-2]
\end{tikzcd}
\end{equation*}

\item[sum morphism] The sum morphism $s\^A\colon\T\^A_2\Rightarrow\T\^A$ is the unique morphism $s\^A\colon z^\*_2\T'_2E\to z^\*\T'E$ induced by the universality of the cartesian lift $\varphi\^z$, which makes the following diagram commutes:
\begin{equation*}
\begin{tikzcd}
{z_2^\*\T'_2E} & {z^\*\T'E} \\
{\T'_2E} & {\T'E}
\arrow["{\varphi\^z}", from=1-2, to=2-2]
\arrow["{s'}"', from=2-1, to=2-2]
\arrow["{s\^A}", dashed, from=1-1, to=1-2]
\arrow["{\varphi\^{z_2}}"', from=1-1, to=2-1]
\end{tikzcd}
\end{equation*}

\item[vertical lift] The vertical lift $l\^A\colon\T\^A\Rightarrow{\T\^A}^2$ is the unique morphism induced by the universality of the cartesian lifts $\varphi\^{z_\T}$ and $\T'\varphi\^z$, which makes the following diagram commutes:
\begin{equation*}
\begin{tikzcd}
{z^\*\T'E} & {z_\T^\*\T'(z^\*\T'E)} \\
& {\T'(z^\*\T'E)} \\
{\T'E} & {{\T'}^2E}
\arrow["{\varphi\^{z_\T}}", from=1-2, to=2-2]
\arrow["{\T'\varphi\^z}", from=2-2, to=3-2]
\arrow["{\varphi\^z}"', from=1-1, to=3-1]
\arrow["{l'}"', from=3-1, to=3-2]
\arrow["{l\^A}", dashed, from=1-1, to=1-2]
\end{tikzcd}
\end{equation*}

\item[canonical flip] The canonical flip $c\^A\colon{\T\^A}^2\Rightarrow{\T\^A}^2$ is the unique morphism induced by the universality of the cartesian lifts $\varphi\^{z_\T}$ and $\T'\varphi\^z$, which makes the following diagram commutes:
\begin{equation*}
\begin{tikzcd}
{z_\T^\*\T'(z^\*\T'E)} & {z_\T^\*\T'(z^\*\T'E)} \\
{\T'(z^\*\T'E)} & {\T'(z^\*\T'E)} \\
{{\T'}^2E} & {{\T'}^2E}
\arrow["{\varphi\^{z_\T}}", from=1-2, to=2-2]
\arrow["{\T'\varphi\^z}", from=2-2, to=3-2]
\arrow["{c'}"', from=3-1, to=3-2]
\arrow["{c\^A}", dashed, from=1-1, to=1-2]
\arrow["{\varphi\^{z_\T}}"', from=1-1, to=2-1]
\arrow["{\T'\varphi\^z}"', from=2-1, to=3-1]
\end{tikzcd}
\end{equation*}
\end{description}

When $(\X',\TT')$ has negatives with negation $n'\colon\T'\Rightarrow\T'$, we can also define:

\begin{description}
\item[negation] The negation $n\^A\colon\T\^A\Rightarrow\T\^A$ is the unique morphism induced by the universality of the cartesian lift $\varphi\^z$, which makes the following diagram commutes:
\begin{equation*}
\begin{tikzcd}
{z^\*\T'E} & {z^\*\T'E} \\
{\T'E} & {\T'E}
\arrow["{\varphi\^z}", from=1-2, to=2-2]
\arrow["{n'}"', from=2-1, to=2-2]
\arrow["{\varphi\^z}"', from=1-1, to=2-1]
\arrow["{n\^A}", dashed, from=1-1, to=1-2]
\end{tikzcd}
\end{equation*}
\end{description}

The next step is to equip each substitution functor $f^\*\colon\Pi^{-1}(A')\to\Pi^{-1}(A)$ induced by a morphism $f\colon A\to A'$ of $\X$ with a suitable distributive law which makes $f^\*$ a strong tangent morphism. Consider the morphism
\begin{align}
\label{equation:distributive-law-subfunctors}
\begin{split}
\alpha\^f\colon&\T\^Af^\*K=z^\*\T'(f^\*E')\xrightarrow{z^\*\kappa\^f}z^\*(\T f)^\*\T'K\xrightarrow{\Comp}(\T f\o z)^\*\T'K\mathdash\\
&\xrightarrow{\T f\o z=z\o f}(z\o f)^\*\T'K\xrightarrow{{\Comp}^\-}f^\*z^\*\T'K=f^\*\T\^{A'}K
\end{split}
\end{align}
where we used the compositor $\Comp$, the distributor $\kappa\^f$ of the tangent bundle functors, and the naturality of $z$.

\begin{remark}
\label{remark:compositor}
To define the distributive law $\alpha\^f$ we make explicit use of the compositor, however, we could have simply written $z^\*(\T f)^\*\T'K=f^\*z^\*\T'K$, since we are using the compositor and its inverse with the naturality of $z$.
\end{remark}

\begin{definition}
\label{definition:indexed-tangent-category}
An \textbf{indexed tangent category} $(\X,\T;\I)$ consists of a tangent category $(\X,\TT)$, called the \textbf{base} tangent category, together with a pseudofunctor $\I\colon\X^\op\to\TngCat_\cong$ from the opposite of the category $\X$ to the $2$-category of tangent categories, strong tangent morphisms, and natural transformations compatible with the distributive laws. For each morphism $f\colon A\to A'$ of the base tangent category, the corresponding strong tangent morphism $(f^\*,\alpha\^f)\=\I(f)$ is called the \textbf{substitution tangent morphism} associated with $f$.
\end{definition}

\begin{remark}
\label{remark:indexed-tangent-category-definition}
In Definition~\ref{definition:indexed-tangent-category}, the tangent structure on the base category seems to be an unnecessary extra structure to carry around. However, such a tangent structure will play a role in reconstructing the base tangent category of the associated tangent fibration in Theorem~\ref{theorem:reduced-Grothendieck-construction}. In particular, the $1$-morphisms of indexed tangent categories we will consider also interact with the base tangent structures.
\end{remark}

\begin{remark}
\label{remark:dorette-vooys-tangent-indexing-functor}
Recently, Pronk and Vooys in~\cite{pronk:equivariant-tangent-cats} introduced a notion of \textit{tangent indexing functor}. In our language, such a functor is an indexed tangent category whose base tangent category is the trivial tangent category $(\X,\1)$ over $\X$ (see Example~\ref{example:trivial-tangent-cat}).
\end{remark}

\begin{notation}
\label{notation:indexed-tangent-cat}
Given an indexed tangent category $(\X,\TT;\I)$, the tangent category $\I(A)$ associated with an object $A$ of $\X$ is denoted by $(\X\^A,\TT\^A)$ and the sustitution tangent morphism associated with a morphism $f\colon A\to A'$ of $\X$, by $(f^\*,\alpha\^f)\colon(\X\^{A'},\TT\^{A'})\nto(\X\^A,\TT\^A)$.
\end{notation}

ILet us rewrite~\cite[Theorem~5.3]{cockett:differential-bundles} in our jargon.

\begin{proposition}
\label{proposition:tangent-fibration-to-indexed-tangent-category}
Every cloven tangent fibration $\Pi\colon(\X',\TT')\to(\X,\TT)$ is associated with an indexed tangent category $\II(\Pi)\=(\X,\TT;\I)$ which sends every object $A$ of $\X$ to the tangent category $(\Pi^{-1}A,\TT\^A)$ defined over the fibre of $A$ and each morphism $f\colon A\to A'$ of $\X$ to the strong tangent morphism $(f^\*,\alpha\^f)$ whose distributive law was defined in Equation~\eqref{equation:distributive-law-subfunctors}.
\end{proposition}

\begin{remark}
\label{remark:colax-strong-version}
In the construction of the strong distributive law $\alpha\^f$ in Equation~\eqref{equation:distributive-law-subfunctors} the invertibility of the distributor $\kappa\^f$ required in Definition~\ref{definition:tangent-fibration} is a sufficient and necessary condition for $(f^\*,\alpha\^f)$ to be a strong tangent morphism. It is natural to wonder what happens when we drop this requirement and only consider a cloven fibration $\Pi\colon\X'\to\X$ whose underlying functor is a strict tangent morphism $\Pi\colon(\X',\TT')\to(\X,\TT)$ but for which we do not require the tangent bundle functors to preserve the cartesian lifts. In particular, one might hope that the fibres of this weaker version of a tangent fibration will still be tangent categories and that the substitution functors will be only \textit{colax} tangent morphisms.
\par Unfortunately, in the definition of the vertical lift $l\^A\colon\T\^A\Rightarrow{\T\^A}^2$ and of the canonical flip $c\^A\colon{\T\^A}^2\Rightarrow{\T\^A}^2$ of the tangent structure $\TT\^A$ over the fibre $\Pi\^{-1}(A)$ of an object $A$ of $\X$, we needed $\T'\varphi\^z$ to be cartesian. This is equivalent to assuming the distributor $\kappa\^z\colon\T'(z^\*E')\to(\T z)^\*(\T'E')$ to be an isomorphism.
\par In a previous version of this paper, we considered a weaker version of a tangent fibration in which we only assumed the distributor $\kappa\^z$ associated with the zero morphism to be invertible. In particular, for an arbitrary morphism $f\colon A\to A'$, the corresponding distributor $\kappa\^f$ was not required to be an isomorphism.
\par This minimal assumption allows one to construct a tangent structure on each fibre and to equip each substitution functor with a \textit{colax} distributive law.
\par We also considered a stronger version in which every distributor $\kappa\^f$ was considered to be the identity. This makes the substitution functors into strict tangent morphisms between the fibres. In our investigations, we did not find any condition on the distributors for which the substitution functors would carry a lax, not necessarily invertible, distributive law. We decided to maintain the notation $(f^\*,\alpha\^f)\colon(\X\^{A'},\TT\^{A'})\nto(\X\^A,\TT\^A)$ to stress the ``\textit{colax nature}'' of the substitution tangent morphisms.
\par So far, we have not found any concrete examples which justify these extensions; therefore, as suggested by the anonymous referees, we decided to omit them in the new version of the paper and only consider the original version due to Cockett and Cruttwell.
\end{remark}

\begin{example}
\label{example:slice-tangent-category}
In example~\ref{example:bundle-tangent-fibration}, we showed that the functor $\Pi\colon\Dsply(\X,\TT)\to(\X,\TT)$ is a tangent fibration, therefore, via Proposition~\ref{proposition:tangent-fibration-to-indexed-tangent-category}, we obtain a corresponding indexed tangent category $\II(\Pi)$ which sends each object $A$ of $\X$ to the slice tangent category of $(\X,\TT)$ over $A$. Concretely, the slice tangent category over $A\in\X$ is the category of tangent display maps whose codomain is the object $A$, and whose tangent bundle functor sends a tangent display map $q\colon E\to A$ to its \textit{vertical bundle} $\T\^Aq\colon VE\to A$, that is the pullback of $\T q$ along the zero morphism:
\begin{equation*}
% https://q.uiver.app/#q=WzAsNCxbMCwwLCJWRSJdLFswLDEsIkEiXSxbMSwxLCJcXFQgQSJdLFsxLDAsIlxcVCBFIl0sWzMsMiwiXFxUIHEiXSxbMSwyLCJ6IiwyXSxbMCwxLCJcXFRcXF5BcSIsMl0sWzAsM10sWzAsMiwiIiwxLHsic3R5bGUiOnsibmFtZSI6ImNvcm5lciJ9fV1d
\begin{tikzcd}
VE & {\T E} \\
A & {\T A}
\arrow[from=1-1, to=1-2]
\arrow["{\T\^Aq}"', from=1-1, to=2-1]
\arrow["\lrcorner"{anchor=center, pos=0.125}, draw=none, from=1-1, to=2-2]
\arrow["{\T q}", from=1-2, to=2-2]
\arrow["z"', from=2-1, to=2-2]
\end{tikzcd}
\end{equation*}
An important property of the slice tangent category is that the differential objects on this tangent category are equivalent to display differential bundles of the original category (see~\cite[Proposition~5.12]{cockett:differential-bundles}).
\end{example}

\begin{example}
\label{example:differential-bundles-indexed-tangent-cat}
The corresponding indexed tangent category of the tangent fibration described in Example~\ref{example:differential-bundle-tangent-fibration} $\Pi\colon\DBnd(\X,\TT)\to(\X,\TT)$ sends each object $A$ of $\X$ to the tangent category of display differential bundles over $A$.
\end{example}

The next step is to extend the assignment which sends a tangent fibration to the corresponding indexed tangent category to a functor. We first need to introduce morphisms of tangent fibrations and morphisms of indexed tangent categories.

\begin{definition}
\label{definition:morphism-of-tangent-fibrations}
Consider two tangent fibrations $\Pi_\s\colon(\X_\s',\TT_\s')\to(\X_\s,\TT_\s)$ and $\Pi_\b\colon(\X'_\b,\TT_\b')\to(\X_\b,\TT_\b)$. A \textbf{lax} (\textbf{colax}) \textbf{morphism of tangent fibrations} $(F,\beta,F',\beta')\colon\Pi_\s\to\Pi_\b$ consists of a pair of lax (colax) tangent morphisms
\begin{align*}
&(F,\beta)\colon(\X_\s,\TT_\s)\to(\X_\b,\TT_\b)\\
&(F',\beta')\colon(\X_\s',\TT_\s')\to(\X_\b',\TT_\b')
\end{align*}
for which the following diagram commutes:
\begin{equation*}
\begin{tikzcd}
{\X_\s'} & {\X_\b'} \\
{\X_\s} & {\X_\b}
\arrow["{\Pi_\s}"', from=1-1, to=2-1]
\arrow["{\Pi_\b}", from=1-2, to=2-2]
\arrow["{(F,\beta)}"', from=2-1, to=2-2]
\arrow["{(F',\beta')}", from=1-1, to=1-2]
\end{tikzcd}
\end{equation*}
in $\TngCat$. Furthermore, $(F,F')\colon\Pi_\s\to\Pi_\b$ is required to be a morphism of fibrations as in Definition~\ref{definition:1-morphisms-fibrations}.\\
Moreover, a morphism of tangent fibration $(F,\beta,F',\beta')$ is:
\begin{itemize}
\item \textbf{strong} if $(F,\beta)$ and $(F',\beta')$ are strong tangent morphisms;
\item \textbf{strict} if $(F,\beta)$ and $(F',\beta')$ are strict;
\item \textbf{strict on the base} whenever the base tangent morphism $(F,\beta)$ is a strict tangent morphism.
\end{itemize}
\end{definition}

\begin{notation}
\label{notation:category-TNGFIB}
Tangent fibrations together with their morphisms can be organized into categories. In particular, $\TngFib$ denotes the category of tangent fibrations and lax morphisms of tangent fibrations; $\TngFib_\co$ denotes the category of tangent fibrations and colax morphisms of tangent fibrations; $\TngFib_\cong$ is the subcategory of $\TngFib$ of tangent fibrations and strong morphisms; finally, $\TngFib_=$ is the subcategory of $\TngFib$ in which only strict morphisms are considered. We also denote by $\TngFib_\strb$, $\TngFib_{\co,\strb}$, and $\TngFib_{\cong,\strb}$ the subcategories of $\TngFib$, $\TngFib_\co$, and $\TngFib_\cong$ in which the morphisms are strict on the base.
\end{notation}

\begin{definition}
\label{definition:morphisms-of-indexed-tangent-categories}
Consider two indexed tangent categories $\I_\s\colon(\X_\s,\TT_\s)\to\TngCat$ and $\I_\b\colon(\X_\b,\TT_\b)\to\TngCat$ a \textbf{lax} (\textbf{colax}) \textbf{morphism of indexed tangent categories} $(F,\beta,F',\beta',\kappa)\colon\I_\s\to\I_\b$ consists of the following data:

\begin{description}
\item[base tangent morphism] A lax (colax) tangent morphism $(F,\beta)\colon(\X_\s,\TT_\s)\to(\X_\b,\TT_\b)$, called the \textbf{base} tangent morphism;

\item[indexed tangent morphism] A collection $(F',\beta')$ of lax (colax) tangent morphisms
\begin{align*}
&(F\^A,\beta\^A)\colon(\X_\s\^A,\TT_\s\^A)\to(\X_\b\^{FA},\TT_\b\^{FA})
\end{align*}
from $\I_\s(A)=(\X_\s\^A,\TT_\s\^A)$ to $\I_\b(FA)=(\X_\b\^{FA},\TT_\b\^{FA})$, indexed by the objects of $\X_\s$;

\item[distributors] A collection $\kappa$ of natural transformations, called \textbf{distributors}
\begin{equation*}
\begin{tikzcd}
{\X_\s\^{A'}} & {\X_\s\^A} \\
{\X_\b\^{FA'}} & {\X_\b\^{FA}}
\arrow["{F\^{A'}}"', from=1-1, to=2-1]
\arrow["{F\^A}", from=1-2, to=2-2]
\arrow["{f^\*}", from=1-1, to=1-2]
\arrow["{(Ff)^\*}"', from=2-1, to=2-2]
\arrow["{\kappa\^f}"{description}, Rightarrow, from=1-2, to=2-1]
\end{tikzcd}
\end{equation*}
indexed by the morphisms $f\colon A\to A'$ of $\X_\s$, and making the following coherence diagram commutes
\begin{equation*}
\begin{tikzcd}
{F\^A\T_\s\^Af^\*} & {F\^Af^\*\T_\s\^{A'}} & {(Ff)^\*F\^{A'}\T_\s\^{A'}}\\
{\T_\b\^{FA}F\^Af^\*} & {\T_\b\^{FA}(Ff)^\*F\^{A'}} & {(Ff)^\*\T_\b\^{FA'}F\^{A'}}
\arrow["{\T_\b\^{FA}\kappa\^f}"', from=2-1, to=2-2]
\arrow["{(\alpha_\b)_{F\^{A'}}\^{Ff}}"', from=2-2, to=2-3]
\arrow["{F\^A\alpha_\s\^f}", from=1-1, to=1-2]
\arrow["{\kappa\^f_{\T_\s\^{A'}}}", from=1-2, to=1-3]
\arrow["{\beta\^Af^\*}"', from=1-1, to=2-1]
\arrow["{(Ff)^\*\beta\^{A'}}", from=1-3, to=2-3]
\end{tikzcd}
\end{equation*}
where $\alpha_\s\^f$ and $\alpha_\b\^{Ff}$ are the distributive laws of $\I_\s(f)$ and $\I_\b(Ff)$, respectively. For colax indexed tangent morphisms the compatibility is the same but the direction of $\beta\^{-}$ is reversed.
\end{description}
Moreover, a morphism of indexed tangent categories is:
\begin{itemize}
\item \textbf{strong} when the base tangent morphism and each indexed tangent morphism $(F\^A,\beta\^A)$ are strong;

\item \textbf{strict} when the base tangent morphism and each indexed tangent morphism $(F\^A,\beta\^A)$ are strict;

\item \textbf{strict on the base} when the base tangent morphism is strict.
\end{itemize}
\end{definition}

\begin{notation}
\label{notation:category-INDXTNG}
Indexed tangent categories together with their morphisms can be organized into categories. In particular, $\IndxTng$ denotes the category of indexed tangent categories and lax morphisms of indexed tangent categories; $\IndxTng_\co$ denotes the category of indexed tangent categories and colax morphisms of indexed tangent categories; $\IndxTng_\cong$ is the subcategory of $\IndxTng$ of indexed tangent categories and strong morphisms; finally, $\IndxTng_=$ is the subcategory of $\IndxTng$ in which only strict morphisms are considered. We also denote by $\IndxTng_\strb$, $\IndxTng_{\co,\strb}$, and $\IndxTng_{\cong,\strb}$ the subcategories of $\IndxTng$, $\IndxTng_\co$, and $\IndxTng_\cong$ whose morphisms are strict on the base.
\end{notation}

Consider two tangent fibrations $\Pi_\s\colon(\X_\s',\TT_\s')\to(\X_\s,\TT_\s)$ and $\Pi_\b\colon(\X_\b',\TT_\b')\to(\X_\b,\TT_\b)$ and a lax morphism of tangent fibrations $(F,F',\beta')\colon\Pi_\s\to\Pi_\b$ and suppose that $(F,F',\beta')$ is strict on the base. Thanks to Proposition~\ref{proposition:tangent-fibration-to-indexed-tangent-category}, $\Pi_\s$ and $\Pi_\b$ define two indexed tangent categories $\I_\s\=\II(\Pi_\s)\colon(\X_\s,\TT_\s)\to\TngCat$ and $\I_\b\=\II(\Pi_\b)\colon(\X_\b,\TT_\b)\to\TngCat$, respectively. The goal is to show that the morphism $(F,F',\beta')$ defines a lax morphism $\I(F,F',\beta')\colon\I_\s\to\I_\b$ of indexed tangent categories which is also strict on the base. Let us start with the base tangent morphism:

\begin{description}
\item[base tangent morphism] The base tangent morphism is $F\colon(\X_\s,\TT_\s)\to(\X_\b,\TT_\b)$;
\end{description}

To define the indexed tangent morphism, let us start by noticing that $F'\colon\X_\s'\to\X_\b'$ maps an object $E$ of the fibre over $A$ of $\Pi_\s$ to $F'E$. However, we have $\Pi_\b F'E=F\Pi_\s E=FA$, thus $F'E$ lives in the fibre over $FA$ of $\Pi_\b$. So, we can restrict $F'$ to the fibres over $A$ and obtain a functor $F\^A\colon\X_\s\^A\to\X_\b\^{FA}$, where $\X_\s\^A=\Pi_\s^{-1}(A)$ and $\X_\b\^{FA}=\Pi_\b^{-1}(FA)$. Let us extend each $F\^A$ to a lax tangent morphism. Consider the following morphism:
\begin{align*}
&\beta\^A\colon F\^A\T_\s\^AE=F'z_\s^\*\T_\s'E\xrightarrow{\kappa\^z}(Fz_\b)^\*F'\T_\s'E\xrightarrow{Fz_\b={z_\b}_F}({z_\b}_F)^\*F'\T_\s'E\xrightarrow{({z_\b}_F)^\*\beta'}({z_\b}_F)^\*\T_\b'F'E=\T_\b\^{FA}F\^AE
\end{align*}

\begin{description}
\item[indexed tangent morphism] The indexed tangent morphism is the collection of lax tangent morphisms $(F\^A,\beta\^A)\colon(\X_\s\^A,\TT_\s\^A)\to(\X_\b\^{FA},\TT_\b\^{FA})$;

\item[distributors] The distributors are precisely the distributors of $(F,F',\beta')$, as in Definition~\ref{definition:1-morphisms-fibrations}.
\end{description}

\begin{remark}
\label{remark:strict-on-base-explanation}
To define the lax distributive law $\beta\^A$ we used that $F$ is a strict tangent morphism when we replace $Fz_\b$ with $(z_\b)_F$. This greatly simplifies the constructions we are going to study in the next section. However, this technical assumption will be dropped in Section~\ref{section:full-Grothendieck} where we will show the most general case.
\end{remark}

Notice that the compatibility between the distributors and the distributive laws is a consequence of $\Pi_\b\beta'=\id_{\X_\b}$ and of the universality of the cartesian lifts. Concretely, this compatibility consists of the commutativity of the following diagram
\begin{equation*}
\begin{tikzcd}
{F'z_\s^\*\T_\s'f^\*K} & {F'z_\s^\*(\T_\s f)^\*\T_\s'K} & {F'f^\*z_\s^\*\T_\s'K} & {(Ff)^\*F'z_\s^\*\T_\s K} \\
\\
{({z_\b}_F)^\*F'\T_\s'f^\*K} & {({z_\b}_F)^\*F'(\T_\s f)^\*\T_\s'K} & {(F\T_\s f)^\*z_\s^\*\T_\s'K} & {(Ff)^\*({z_\b}_F)^\*F'\T_\s K} \\
\\
{({z_\b}_F)^\*\T_\b'F'f^\*K} & {z_\b^\*\T_\b'(Ff)^\*F'K} & {z_\b^\*(\T_\b Ff)^\*\T_\b'F'K} & {(Ff)^\*({z_\b}_F)^\*\T_\b'F'K}
\arrow["{({z_\b}_F)^\*\T_\b'\lambda\^f}"', from=5-1, to=5-2]
\arrow["{({z_\b}_F)^\*\kappa_\b\^{Ff}F'}"', from=5-2, to=5-3]
\arrow[Rightarrow, no head, from=5-3, to=5-4]
\arrow["{(Ff)^\*z_\b^\*\beta'}", from=3-4, to=5-4]
\arrow["{(Ff)^\*\lambda\^z\T_\s'}", from=1-4, to=3-4]
\arrow["{\lambda\^fz_\s^\*\T_\s'}", from=1-3, to=1-4]
\arrow[Rightarrow, no head, from=1-2, to=1-3]
\arrow["{F'z_\s^\*\kappa_\s\^f}", from=1-1, to=1-2]
\arrow["{\lambda\^z\T_\s'f^\*}"', from=1-1, to=3-1]
\arrow["{({z_\b}_F)^\*\beta'f^\*}"', from=3-1, to=5-1]
\arrow["{\lambda\^z(\T_\s f)^\*\T_\s'}", from=1-2, to=3-2]
\arrow["{({z_\b}_F)^\*F'\kappa_\s\^f}", from=3-1, to=3-2]
\arrow["{\lambda\^{\T_\s f}\T_\s'}", from=3-2, to=3-3]
\arrow[Rightarrow, no head, from=3-3, to=3-4]
\end{tikzcd}
\end{equation*}
where $\lambda$ denotes the distributors of the $1$-morphism $(F,F')$ and $\kappa$ denotes the distributors of the tangent bundle functors.
\par To prove the commutativity of the bottom rectangular diagram, let us post-compose the two paths of morphisms by $\varphi\^{Ff}\varphi_\b\^{z_\b}$, as follows:
\begin{equation*}
\adjustbox{scale=.9,center}{
% https://q.uiver.app/#q=WzAsOSxbMCwwLCIoel9cXGIpX0ZeXFwqRidcXFRfXFxvJ2ZeXFwqSyJdLFsxLDAsIih6X1xcYilfRl5cXCpGJyhcXFRfXFxvIGYpXlxcKlxcVF9cXG8nSyJdLFsyLDAsIih6X1xcYilfRl5cXCooRlxcVF9cXG8gZileXFwqRidcXFRfXFxvJ0siXSxbMywwLCIoRmYpXlxcKih6X1xcYilfRl5cXCpGJ1xcVF9cXG8nSyJdLFs0LDAsIihGZileXFwqKHpfXFxiKV9GXlxcKlxcVF9cXGInRidLIl0sWzAsMSwiRidcXFRfXFxvJ2ZeXFwqSyJdLFs0LDEsIih6X1xcYilfRl5cXCpcXFRfXFxiJ0YnSyJdLFs0LDIsIlxcVF9cXGInRidLIl0sWzAsMiwiRidcXFRfXFxvJ0siXSxbMCwxLCIoel9cXGIpX0ZeXFwqRidcXGthcHBhX1xcb1xcXmYiXSxbMSwyLCIoel9cXGIpX0ZeXFwqXFxsYW1iZGFcXF5mXFxUJ19cXG8iXSxbMiwzLCIiLDAseyJsZXZlbCI6Miwic3R5bGUiOnsiaGVhZCI6eyJuYW1lIjoibm9uZSJ9fX1dLFszLDQsIihGZileXFwqKHpfXFxiKV9GXlxcKlxcYmV0YSciXSxbOCw3LCJcXGJldGEnIiwyXSxbMCw1LCJcXHZhcnBoaV9cXGJeXFwqIiwyXSxbNSw4LCJGJ1xcVF9cXG8nXFx2YXJwaGlcXF5mIiwyXSxbNCw2LCJcXHZhcnBoaVxcXntGZn0iXSxbNiw3LCJcXHZhcnBoaVxcXnpfXFxiIl1d
\begin{tikzcd}
{({z_\b}_F)^\*F'\T_\s'f^\*K} & {({z_\b}_F)^\*F'(\T_\s f)^\*\T_\s'K} & {({z_\b}_F)^\*(F\T_\s f)^\*F'\T_\s'K} & {(Ff)^\*({z_\b}_F)^\*F'\T_\s'K} & {(Ff)^\*({z_\b}_F)^\*\T_\b'F'K} \\
{F'\T_\s'f^\*K} &&&& {({z_\b}_F)^\*\T_\b'F'K} \\
{F'\T_\s'K} &&&& {\T_\b'F'K}
\arrow["{({z_\b}_F)^\*F'\kappa_\s\^f}", from=1-1, to=1-2]
\arrow["{\varphi_\b^\*}"', from=1-1, to=2-1]
\arrow["{({z_\b}_F)^\*\lambda\^f\T'_\s}", from=1-2, to=1-3]
\arrow[Rightarrow, no head, from=1-3, to=1-4]
\arrow["{(Ff)^\*({z_\b}_F)^\*\beta'}", from=1-4, to=1-5]
\arrow["{\varphi\^{Ff}}", from=1-5, to=2-5]
\arrow["{F'\T_\s'\varphi\^f}"', from=2-1, to=3-1]
\arrow["{\varphi\^z_\b}", from=2-5, to=3-5]
\arrow["{\beta'}"', from=3-1, to=3-5]
\end{tikzcd}
}
\end{equation*}
\begin{equation*}
\adjustbox{scale=.9,center}{
\begin{tikzcd}
{({z_\b}_F)^\*F'\T_\s'f^\*K} && {({z_\b}_F)^\*\T_\b'F'f^\*K} & {z_\b^\*\T_\b'(Ff)^\*F'K} & {z_\b^\*(\T_\b Ff)^\*\T_\b'F'K} & {(Ff)^\*({z_\b}_F)^\*\T_\b'F'K} \\
\\
{F'\T_\s'f^\*K} && {\T_\b'F'f^\*K} & {\T_\b'(Ff)^\*F'K} & {(\T_\b Ff)^\*\T_\b'F'K} & {({z_\b}_F)^\*\T'_\b F'K} \\
\\
{F'\T_\s'K} &&& {\T_\b'F'K} && {\T'_\b F'K}
\arrow["{({z_\b}_F)^\*\T_\b'\lambda\^f}", from=1-3, to=1-4]
\arrow["{({z_\b}_F)^\*\kappa_\b\^{Ff}F'}", from=1-4, to=1-5]
\arrow[Rightarrow, no head, from=1-5, to=1-6]
\arrow["{({z_\b}_F)^\*\beta'f^\*}", from=1-1, to=1-3]
\arrow["{\varphi\^{Ff}}", from=1-6, to=3-6]
\arrow["{\varphi_\b\^{z_\b}}"', from=1-1, to=3-1]
\arrow["{F'\T_\s'\varphi\^f}"', from=3-1, to=5-1]
\arrow["{\beta'}"', from=5-1, to=5-4]
\arrow["{\varphi_\b\^{z_\b}}"', from=1-5, to=3-5]
\arrow["{\varphi\^{\T_\b Ff}}"', from=3-5, to=5-4]
\arrow["{\varphi_\b\^{z_\b}}"', from=1-4, to=3-4]
\arrow["{\T_\b'\varphi\^{Ff}}"', from=3-4, to=5-4]
\arrow["{\varphi_\b\^{z_\b}}"', from=1-3, to=3-3]
\arrow["{\T_\b'F'\varphi\^f}"', from=3-3, to=5-4]
\arrow["{\varphi_\b\^{z_\b}}", from=3-6, to=5-6]
\arrow[Rightarrow, no head, from=5-4, to=5-6]
\end{tikzcd}
}
\end{equation*}
Since $\Pi_\b\beta'$ is the identity, from the cartesian property of $\varphi\^{Ff}\varphi_\b\^{z_\b}$ we conclude that these two paths of morphisms are identical.

\begin{proposition}
\label{proposition:functor-from-tangent-fibrations-to-indexed-tangent-categories}
The assignment which sends a cloven tangent fibration $\Pi$ to the corresponding indexed tangent category $\II(\Pi)$ extends to a pair of functors. In particular, each lax (colax) morphism of tangent fibrations strict on the base is associated with a lax (colax) morphism of indexed tangent categories, also strict on the base:
\begin{align*}
&\II\colon\TngFib_\strb\to\IndxTng_\strb\\
&\II\colon\TngFib_{\co,\strb}\to\IndxTng_{\co,\strb}
\end{align*}
These functors restrict to the subcategory of strong and strict morphisms of tangent fibrations which are strict on the base:
\begin{align*}
&\II\colon\TngFib_{\cong,\strb}\to\IndxTng_{\cong,\strb}\\
&\II\colon\TngFib_{=,\strb}\to\IndxTng_{=,\strb}
\end{align*}
\end{proposition}

%__________________________________________________________________________
\subsection{The reduced Grothendieck construction for tangent fibrations}
\label{subsection:reduced-Grothendieck}
In the previous section, we extended Cockett and Cruttwell's result to a functorial assignment which sends a tangent fibration to a corresponding indexed tangent category. Crucially, for a tangent fibration $\Pi\colon(\X',\TT')\to(\X,\TT)$ the tangent bundle functor $\T'$ of the total tangent category sends an object $E$ of the fibre over $A\in\X$, to $\T'E$, which lives in the fibre over $\T A$. So, in order to define the tangent bundle functor $\T\^A\colon\Pi^{-1}A\to\Pi^{-1}A$ we post-composed the total tangent bundle functor $\T'$ with the substitution functor $z^\*$ induced by the zero morphism of $(\X,\TT)$.
\par This process inevitably destroys part of the information of the total tangent bundle functor. To convince the reader of this phenomenon, consider a tangent category $(\X,\TT)$ equipped with a display system as in Example~\ref{example:bundle-tangent-fibration}. This produces a tangent fibration $\Pi\colon\Dsply(\X,\TT)\to(\X,\TT)$ whose total tangent category is the tangent category of tangent display maps $E\to A$ of $\X$.
\par As shown in Example~\ref{example:slice-tangent-category}, the corresponding indexed tangent category $\II(\Pi)$ sends each object $A$ of $\X$ to the slice tangent category $(\X,\TT)/A$, whose tangent bundle functor $\TT\^A$ sends a tangent display map $q\colon E\to A$ to the vertical bundle $\T\^Aq\colon VE\to A$, defined by the pullback diagram:
\begin{equation*}
\begin{tikzcd}
{VE} & {\T E} \\
A & {\T A}
\arrow["z"', from=2-1, to=2-2]
\arrow["{\T q}", from=1-2, to=2-2]
\arrow[from=1-1, to=1-2]
\arrow["{\T\^Aq}"', from=1-1, to=2-1]
\arrow["\lrcorner"{anchor=center, pos=0.125}, draw=none, from=1-1, to=2-2]
\end{tikzcd}
\end{equation*}
In general, $VE$ is only strictly included in $\T E$. For instance, in the tangent category of Example~\ref{example:algebraic-geometry}, consider a ring $R$ and let us consider its tangent bundle $p\colon\T R\to R$, where $\T R$ is the ring $R[\epsilon]$, with $\epsilon^2=0$, and $p$ sends $\epsilon$ to $0$. From the universality of the vertical lift, one finds out that the vertical bundle of $p$ consists of the pullback $\T_2R=R[\epsilon_1,\epsilon_2]$ of $p$ along itself, which is the ring of polynomials with coefficients in $R$ in two variables $\epsilon_1$ and $\epsilon_2$, such that $\epsilon_i\epsilon_j=0$, for $i,j=1,2$.
\par Contrarily, $\T(\T R)=\T^2R$ is the ring $R[\epsilon][\epsilon']$, generated by two variables $\epsilon$ and $\epsilon'$, such that $\epsilon^2=\epsilon'^2=0$, but $\epsilon'\epsilon\neq 0$. The map $R[\epsilon_1,\epsilon_2]\to R[\epsilon][\epsilon']$ which sends $\epsilon_1$ to $\epsilon$ and $\epsilon_2$ to $\epsilon'$ is precisely the strict inclusion of the vertical bundle into the tangent bundle.\newline
\par In this section, we partially reconstruct the original tangent fibration starting from the associated indexed tangent category, by constructing a functor of type $\FF\colon\IndxTng\to\TngFib$. We show that this functor together with the functor of Proposition~\ref{proposition:functor-from-tangent-fibrations-to-indexed-tangent-categories} forms an adjunction, which, however, is not an equivalence of categories.
\par The main idea is to equip the category of elements $\El(\X,\I)$ (see Section~\ref{subsection:fibrations-background}) of the underlying indexed category $\I\colon\X^\op\to\TngCat_\cong\to\Cat$ with a suitable tangent structure, obtained by ``gluing'' together the tangent structures of the fibres.

\begin{description}
\item[tangent bundle functor] The tangent bundle functor $\hat\T\colon\El(\X,\I)\to\El(\X,\I)$ sends an object $(A,E)$ to $(\T A,p^\*\T\^AE)$ and a morphism $(f,\xi\^f)\colon(A,E)\to(A',E')$ to the morphism
\begin{align*}
(\T f,\xi\^{\T f})\colon(\T A,p^\*\T\^AE)\to(\T A',p^\*\T\^{A'}E')
\end{align*}
defined by $\T f\colon\T A\to\T A'$ together with the morphism
\begin{align*}
\xi\^{\T f}&\colon p^\*\T\^AE\xrightarrow{p^\*\T\^A\xi\^f}p^\*\T\^Af^\*E'\xrightarrow{p^\*\alpha\^f}p^\*f^\*\T\^{A'}E'\xrightarrow{\Comp}(f\o p)^\*\T\^{A'}E'\mathdash\\
&\xrightarrow{f\o p=p\o \T f}(p\o \T f)^\*\T\^{A'}E'\xrightarrow{{\Comp}^{-1}}(\T f)^\*p^\*\T\^{A'}E'
\end{align*}
where $\Comp$ denotes the compositor;

\item[projection] The projection $\hat p\colon\hat\T\Rightarrow\id_{\El(\X,\I)}$ is the pair
\begin{align*}
&(p,\xi\^p)\colon(\T A,p^\*\T\^AE)\to(A,E)
\end{align*}
formed by $p\colon\T A\to A$ and by the morphism:
\begin{align*}
&\xi\^p\colon p^\*\T\^AE\xrightarrow{p^\*p\^A}p^\*E
\end{align*}

\item[zero morphism] The zero morphism $\hat z\colon\id_{\El(\X,\I)}\Rightarrow\hat\T$ is the pair
\begin{align*}
&(z,\xi\^z)\colon(A,E)\to(\T A,p^\*\T\^AE)
\end{align*}
formed by $z\colon A\to\T A$ and by the morphism:
\begin{align*}
&\xi\^z\colon E\xrightarrow{z\^A}\T\^AE\xrightarrow{\Unit}\id_A^\*\T\^AE\xrightarrow{p\o z=\id_A}(p\o z)^\*\T\^AE\xrightarrow{\Comp}z^\*p^\*\T\^AE
\end{align*}

\item[$n$-fold pullback] The $n$-fold pullback along the projection $p\^A$ is given by:
\begin{align*}
&\hat\T_n(A,E)=(\T_nA,\pi_1^\*p^\*\T\^A_nE)
\end{align*}
where the $k$-th projection $\hat\pi_k\colon\hat\T_n\Rightarrow\hat\T$ is given by the pair
\begin{align*}
&(\pi_k,\xi\^{\pi_k})\colon(\T_nA,\pi_1^\*p^\*\T\^A_nE)\to(\T A,p^\*\T\^AE)
\end{align*}
where $\pi_k\colon\T_nA\to\T A$ and $\xi\^{\pi_k}$ is the morphism:
\begin{align*}
\xi\^{\pi_k}&\colon\pi_1^\*p^\*\T_n\^AE\xrightarrow{\pi_1^\*p^\*\pi_k\^A}\pi_1^\*p^\*\T\^AE\xrightarrow{{\Comp}^{-1}}(p\o\pi_1)^\*\T\^AE\mathdash\\
&\xrightarrow{p\o\pi_1=p\o\pi_k}(p\o\pi_k)^\*\T\^AE\xrightarrow{\Comp}\pi_k^\*p^\*\T\^AE
\end{align*}

\item[sum morphism] The sum morphism $\hat s\colon\hat\T_2\Rightarrow\hat\T$ is the pair
\begin{align*}
&(s,\xi\^s)\colon(\T_2A,\pi_1^\*\T\^A_2E)\to(\T A,p^\*\T\^AE)
\end{align*}
formed by $s\colon\T_2A\to\T A$ and by the morphism:
\begin{align*}
\xi\^s\colon&\pi_1^\*p^\*\T\^A_2E\xrightarrow{\pi_1^\*p^\*s\^A}\pi_1^\*p^\*\T\^AE\xrightarrow{{\Comp}^{-1}}(p\o \pi_1)^\*\T\^AE\mathdash\\
&\xrightarrow{p\o s=p\o\pi_1}(p\o s)^\*\T\^AE\xrightarrow{\Comp}s^\*p^\*\T\^AE
\end{align*}

\item[vertical lift] The vertical lift $\hat l\colon\hat\T\Rightarrow{\hat\T}^2\=\hat\T\o\hat\T$ is the pair
\begin{align*}
&(l,\xi\^l)\colon(\T A,p^\*\T\^AE)\to(\T^2A,p_\T^\*\T\^{\T A}p^\*\T\^AE)
\end{align*}
formed by $l\colon\T A\to\T^2A$ and by the morphism:
\begin{align*}
\xi\^l&\colon p^\*\T\^AE\xrightarrow{p^\*l\^A}p^\*{\T\^A}^2E\xrightarrow{p=p\o z\o p=p\o p_\T\o l}(p\o p_\T\o l)^\*{\T\^A}^2E\mathdash\\
&\xrightarrow{\Comp}l^\*p_\T^\*p^\*{\T\^A}^2E\xrightarrow{l^\*p^\*_\T\alpha\^p\T\^A}l^\*p^\*_\T\T\^{\T A}p^\*\T\^AE
\end{align*}

\item[canonical flip] The canonical flip $\hat c\colon{\hat\T}^2\Rightarrow{\hat\T}^2$ is the pair
\begin{align*}
&(c,\xi\^c)\colon(\T^2A,p_\T^\*\T\^{\T A}p^\*\T\^AE)\to(\T^2A,p_\T^\*\T\^{\T A}p^\*\T\^AE)
\end{align*}
formed by $c\colon\T^2A\to\T^2A$ and by the morphism:
\begin{align*}
\xi\^c&\colon p^\*_\T\T\^{\T A}p^\*\T\^AE\xrightarrow{p^\*_\T{\alpha\^p}^{-1}\T\^A}p^\*_\T p^\*{\T\^A}^2E\xrightarrow{p^\*_\T p^\*c\^A}p^\*_\T p^\*{\T\^A}^2E\mathdash\\
&\xrightarrow{{\Comp}^{-1}}(p\o p_\T)^\*{\T\^A}^2E\xrightarrow{p\o p_\T=p\o\T p\o c=p\o p_\T\o c}(p\o p_\T\o c)^\*{\T\^A}^2E\mathdash\\
&\xrightarrow{\Comp}c^\*p_\T^\*p^\*{\T\^A}^2E\xrightarrow{c^\*p^\*_\T\alpha\^p{\T\^A}^2}c^\*p^\*_\T\T\^{\T A}p^\*\T\^AE
\end{align*}
\end{description}
Moreover, if $(\X,\TT)$ has negatives with negation $n\colon\T\Rightarrow\T$ and each $(\X\^A,\TT\^A)$ has negatives with negation $n\^A\colon\T\^A\Rightarrow\T\^A$, then also the tangent category of elements has negatives whose negation is defined as follows:
\begin{description}
\item[negation] The negation $\hat n\colon\hat\T\Rightarrow\hat\T$ is the pair
\begin{align*}
&(n,\xi\^n)\colon(\T A,p^\*\T\^AE)\to(\T A,p^\*\T\^AE)
\end{align*}
formed by $n\colon\T A\to\T A$ and by the morphism:
\begin{align*}
\xi\^n&\colon p^\*\T\^AE\xrightarrow{p^\*n\^A}p^\*\T\^AE\xrightarrow{p=p\o n}(p\o n)^\*\T\^AE\xrightarrow{\Comp}n^\*p^\*\T\^AE
\end{align*}
\end{description}

\begin{lemma}
\label{lemma:tangent-category-of-elements}
The category $\El(\X,\TT,\I)$ of elements of an indexed tangent category $(\X,\TT,\I)$ comes equipped with a tangent structure:
\begin{align*}
&\hat\TT\=\left(\hat\T,\hat p,\hat z,\hat s,\hat l,\hat c\right)
\end{align*}
Furthermore, if $(\X,\TT)$ and each tangent category $(\X\^A,\TT\^A)$ admit negatives for each $A$ of $\X$, so does $\El(\X,\TT,\I)$.
\begin{proof}
The equational axioms required for $\hat\TT$ to define a tangent structure are a direct, but tedious, consequence of $(\X,\TT)$ and each $(\X\^A,\TT\^A)$ being tangent categories. As an example, let us show that the zero morphism $\hat z$ is a section of the projection $\hat p$, that is:
\begin{align*}
&(p,\xi\^p)\o(z,\xi\^z)=(\id_A,\xi\^{\id_A})\colon(A,E)\to(A,E)
\end{align*}
For starters, notice that the equation holds for the base components of the morphisms:
\begin{align*}
&p\o z=\id_A
\end{align*}
Let us focus on the fibre components. First, let us recall that
\begin{align*}
&(p,\xi\^p)\o(z,\xi\^z)=(p\o z,\xi\^{p\o z})
\end{align*}
where:
\begin{align*}
&\xi\^{p\o z}\colon E\xrightarrow{\xi\^z}z^\*p^\*\T\^AE\xrightarrow{z^\*\xi\^p}z^\*p^\*E\xrightarrow{\Comp}(p\o z)^\*E=\id_A^\*E
\end{align*}
However, by definition:
\begin{align*}
&\xi\^z\colon E\xrightarrow{z\^A}\T\^AE\xrightarrow{\Unit}\id_A^\*\T\^AE=(p\o z)^\*\T\^AE\xrightarrow{{\Comp}^{-1}}z^\*p^\*\T\^AE\\
&\xi\^p\colon p^\*\T\^AE\xrightarrow{p^\*p\^A}p^\*E
\end{align*}
Thus, using that $p\^A\o z\^A=\id_{E}$, it follows straightforwardly that
\begin{align*}
&\xi\^{p\o z}=\Unit
\end{align*}
proving the desired equation. One can prove all the other equational axioms in a similar way. The remaining axioms to be proven are the existence of the $n$-fold pullbacks of the projection along itself, together with the universality of the vertical lift. Let us focus on the existence of the $n$-fold pullback of the projection along itself.
\par let us start by considering an object $(A',E')$ of $\El(\X,\TT,\I)$ and a collection morphisms of $\El(\X,\TT,\I)$
\begin{align*}
&(f_k,\xi\^f_k)\colon(A',E')\to\hat\T(A,E)
\end{align*}
such that, for each $i,j=1\,n$:
\begin{align*}
&\hat p\o(f_i,\xi\^f_i)=\hat p\o(f_j,\xi\^f_j)
\end{align*}
From the universality of the $n$-fold pullback of $p$ along itself, we obtain a morphism $\<f_1\,f_n\>\colon A'\dashrightarrow\T_nA$ of $\X$:
\begin{equation*}
% https://q.uiver.app/#q=WzAsNSxbMCwwLCJBJyJdLFsxLDEsIlxcVF9uQSJdLFsxLDIsIlxcVCBBIl0sWzIsMSwiXFxUIEEiXSxbMiwyLCJBIl0sWzIsNCwicCIsMl0sWzMsNCwicCJdLFsxLDIsIlxccGlfMSIsMl0sWzEsMywiXFxwaV9uIl0sWzIsMywiXFxkb3RzIiwzLHsib2Zmc2V0IjotMywic3R5bGUiOnsiYm9keSI6eyJuYW1lIjoibm9uZSJ9LCJoZWFkIjp7Im5hbWUiOiJub25lIn19fV0sWzAsMiwiZl8xIiwyLHsiY3VydmUiOjN9XSxbMCwzLCJmXzIiLDAseyJjdXJ2ZSI6LTN9XSxbMCwxLCJcXDxmXzFcXCxmX25cXD4iLDEseyJzdHlsZSI6eyJib2R5Ijp7Im5hbWUiOiJkYXNoZWQifX19XV0=
\begin{tikzcd}
{A'} \\
& {\T_nA} & {\T A} \\
& {\T A} & A
\arrow["{\<f_1\,f_n\>}"{description}, dashed, from=1-1, to=2-2]
\arrow["{f_n}", bend left, from=1-1, to=2-3]
\arrow["{f_1}"', bend right, from=1-1, to=3-2]
\arrow["{\pi_n}", from=2-2, to=2-3]
\arrow["{\pi_1}"', from=2-2, to=3-2]
\arrow["p", from=2-3, to=3-3]
\arrow["\dots"{marking, allow upside down}, shift left=3, draw=none, from=3-2, to=2-3]
\arrow["p"', from=3-2, to=3-3]
\end{tikzcd}
\end{equation*}
let us now consider the following diagram in $\X\^{A'}$:
\begin{equation*}
% https://q.uiver.app/#q=WzAsNyxbMSwxLCIocFxcbyBmXzEpXlxcKlxcVF9uXFxeQUUiXSxbMSwyLCIocFxcbyBmXzEpXlxcKlxcVFxcXkFFIl0sWzIsMSwiKHBcXG8gZl8xKV5cXCpcXFRcXF5BRSJdLFsyLDIsIihwXFxvIGZfMSleXFwqRSJdLFswLDIsImZfMV5cXCpwXlxcKlxcVF9uXFxeQUUiXSxbMCwwLCJFJyJdLFsyLDAsImZfbl5cXCpwXlxcKlxcVF9uXFxeQUUiXSxbMCwxLCIocFxcbyBmXzEpXlxcKlxceGlcXF57XFxwaV8xfSIsMl0sWzAsMiwiKHBcXG8gZl8xKV5cXCpcXHhpXFxee1xccGlfbn0iXSxbMSwzLCIocFxcbyBmXzEpXlxcKnBcXF5BIiwyXSxbMiwzLCIocFxcbyBmXzEpXlxcKnBcXF5BIl0sWzAsMywiIiwxLHsic3R5bGUiOnsibmFtZSI6ImNvcm5lciJ9fV0sWzEsMiwiXFxkb3RzIiwzLHsib2Zmc2V0IjotMywic3R5bGUiOnsiYm9keSI6eyJuYW1lIjoibm9uZSJ9LCJoZWFkIjp7Im5hbWUiOiJub25lIn19fV0sWzUsNCwiXFx4aVxcXmZfMSIsMl0sWzQsMSwiXFxDb21wIiwyXSxbNSw2LCJcXHhpXFxeZl9uIl0sWzYsMiwiXFxDb21wIl0sWzUsMCwiXFx4aVxcXntcXDxmXzFcXCxmX25cXD59IiwxLHsic3R5bGUiOnsiYm9keSI6eyJuYW1lIjoiZGFzaGVkIn19fV1d
\begin{tikzcd}
{E'} && {f_n^\*p^\*\T_n\^AE} \\
& {(p\o f_1)^\*\T_n\^AE} & {(p\o f_1)^\*\T\^AE} \\
{f_1^\*p^\*\T_n\^AE} & {(p\o f_1)^\*\T\^AE} & {(p\o f_1)^\*E}
\arrow["{\xi\^f_n}", from=1-1, to=1-3]
\arrow["{\xi\^f_1}"', from=1-1, to=3-1]
\arrow["\Comp", from=1-3, to=2-3]
\arrow["{(p\o f_1)^\*\xi\^{\pi_n}}", from=2-2, to=2-3]
\arrow["{(p\o f_1)^\*\xi\^{\pi_1}}"', from=2-2, to=3-2]
\arrow["\lrcorner"{anchor=center, pos=0.125}, draw=none, from=2-2, to=3-3]
\arrow["{(p\o f_1)^\*p\^A}", from=2-3, to=3-3]
\arrow["\Comp"', from=3-1, to=3-2]
\arrow["\dots"{marking, allow upside down}, shift left=3, draw=none, from=3-2, to=2-3]
\arrow["{(p\o f_1)^\*p\^A}"', from=3-2, to=3-3]
\end{tikzcd}
\end{equation*}
Since the substitution tangent morphisms are strong, the functor $(p\o f_1)^\*$ preserves the $n$-fold pullbacks of the projection $p\^A$ with itself. Therefore, we obtain a unique morphism $E'\dashrightarrow(p\o f_1)^\*\T_n\^AE$. let us define:
\begin{align*}
&\xi\^{\<f_1\,f_n\>}\colon E'\dashrightarrow(p\o f_1)^\*\T_n\^AE\xrightarrow{f_1=\pi_1\o\<f_1\,f_n\>}(p\o \pi_1\o\<f_1\,f_n\>)^\*\T_n\^AE\xrightarrow{\Comp}\<f_1\,f_n\>^\*(p\o \pi_1)^\*\T_n\^AE
\end{align*}
Thus, we obtain a morphism:
\begin{align*}
&\left\<(f_1,\xi\^f_1)\,(f_n,\xi\^f_n)\right\>\=\left(\<f_1\,f_n\>,\xi\^{\<f_1\,f_n\>}\right)\colon(A',E')\to\hat\T_n(A,E)
\end{align*}
It is easy to prove that this is the unique morphism satisfying the equation
\begin{align*}
&\hat\pi_k\o\left\<(f_1,\xi\^f_1)\,(f_n,\xi\^f_n)\right\>=(f_k,\xi\^f_k)
\end{align*}
for each $k=1\,n$. Finally, to prove the universal property of the vertical lift, first notice that since $p^\*$ is a strong tangent morphism, it preserves the universal property of the vertical lift $l\^A\colon\T\^A\to{\T\^A}^2$ of $(\X\^A,\TT\^A)$. From this observation and the universal property of the vertical lift of the base tangent category $(\X,\T)$ we prove the desired universal property.
\end{proof}
\end{lemma}

\begin{definition}
\label{definition:tangent-category-of-elements}
The tangent category defined in Lemma~\ref{lemma:tangent-category-of-elements}, is called the \textbf{tangent category of elements} of an indexed tangent category $\I$.
\end{definition}

It is not hard to see that the functor $\II(\Pi)\colon\El(\X,\TT;\I)\to(\X,\TT)$ which sends each object $(A,E)$ of $\El(\X,\TT;\I)$ to $A$ defines a cloven tangent fibration. In particular, by construction, $\Pi$ strictly preserves the tangent structures. To prove that the tangent bundle functors preserve the cartesian lifts, notice that the cartesian lift of a morphism $f\colon A\to A'$ of $\X$ along $\Pi$ is the morphism
\begin{align*}
&\varphi\^f\=(f,\xi\^f)\colon(A,f^\*E')\to(A',E')
\end{align*}
where:
\begin{align*}
&\xi^f\colon f^\*E'\xrightarrow{f^\*\id_{E'}}f^\*E'
\end{align*}
From this, one can see that $\hat\T(\varphi\^f)$ is the morphism $(\T f,\alpha\^f)$, where $\alpha\^f$ is the distributive law of $f^\*$. Since $f^\*$ is a strong tangent morphism, it follows directly that the distributor, which coincides with $\alpha\^f$, is an isomorphism.\newline
\par The next step is to prove that this assignment extends to a functor $\FF\colon\IndxTng_\strb\to\TngFib_\strb$ which forms an adjunction with $\II$ of Proposition~\ref{proposition:functor-from-tangent-fibrations-to-indexed-tangent-categories}.
\par Let us start by considering two indexed tangent categories $\I_\s\colon(\X_\s,\TT_\s)\to\TngCat$ and $\I_\b\colon(\X_\b,\TT_\b)\to\TngCat$ and a lax (colax) morphism of indexed tangent categories $(F,F',\beta',\kappa)\colon\I_\s\to\I_\b$ which is strict on the base. The goal is to define a lax (colax) morphism of tangent fibrations strict on the base, between the corresponding tangent fibrations $\Pi(\I_\s)$ and $\Pi(\I_\b)$:
\begin{description}
\item[base tangent morphism] The base tangent morphism is given by $F\colon(\X_\s,\TT_\s)\to(\X_\b,\TT_\b)$;

\item[total tangent morphism] The lax (colax) tangent morphism between the total categories is given by the functor $\El(F,F')\colon\El(\X_\s,\I_\s)\to\El(\X_\b,\I_b)$ which sends a pair $(A,E$ to $(FA,F'E)$ and a morphism $(f,\xi\^f)\colon(A,E)\to(A',E')$ to the morphism $(Ff,\xi\^{Ff})\colon(FA,F'E)\to(FA',F'E')$ where $\xi\^{Ff}$ is the morphism:
\begin{align*}
&\xi\^{Ff}\colon F'E\xrightarrow{F'\xi\^f}F'f^\*E\xrightarrow{\kappa\^f}(Ff)^\*F'K
\end{align*}
The lax distributive law $\beta'\colon\El(F,F')\o\hat\T_\s\Rightarrow\hat\T_\b\o\El(F,F')$ is the morphism $(\beta,\xi\^\beta)\colon(F\T_\s A,F'p_\s^\*\T_\s\^AE)\to(\T_\b FA,(p_\b)_F^\*\T_\b\^{FA}F'E)$ given by $F\T_\s A=\T_\b FA$ and by the morphism
\begin{align*}
&\xi\^\beta\colon F'p_\s^\*\T_\s\^AE\xrightarrow{\kappa\^{p_\s}\T_\s\^A}(F p_\s)^\*F'\T_\s\^AE\xrightarrow{(Fp_\s)^\*\beta\^A}(Fp_\s)^\*\T_\b\^{FA}F'E=(p_\b)_F^\*\T_\b\^{FA}F'E
\end{align*}
where we used that $Fp_\s=(p_\b)_F$. When the morphism of tangent fibrations is colax, the corresponding colax distributive law of $\El(F,F')$ is the morphism $(\beta,\xi\^\beta)\colon(\T_\b FA,(p_\b)_F^\*\T_\b\^{FA}F'E)\nto(F\T_\s A,F'p_\s^\*\T_\s\^AE)$, given by $\T_\b FA=F\T_\s A$ and by the morphism:
\begin{align*}
&\xi\^\beta\colon(p_\b)_F^\*\T_\b\^{FA}F'E=(Fp_\s)^\*\T_\b\^{FA}F'E\xrightarrow{(Fp_\s)^\*\beta\^A}(F p_\s)^\*F'\T_\s\^AE\xrightarrow{{\kappa\^{p_\s}}^{-1}\T_\s\^A}F'p_\s^\*\T_\s\^AE
\end{align*}
\end{description}

\begin{proposition}
\label{proposition:functor-from-indexed-tangent-categories-to-tangent-fibrations}
The assignment which sends an indexed tangent category $\I$ to the corresponding cloven tangent fibration $\Pi(\I)$ extends to a pair of functors. In particular, each lax (colax) morphism of indexed tangent categories strict on the base is associated with a lax (colax) morphism of tangent fibrations, also strict on the base:
\begin{align*}
&\FF\colon\IndxTng_\strb\to\TngFib_\strb\\
&\FF\colon\IndxTng_{\co,\strb}\to\TngFib_{\co,\strb}
\end{align*}
These functors restrict to the subcategory of strong and strict morphisms of indexed tangent categories which are strict on the base:
\begin{align*}
&\FF\colon\IndxTng_{\cong,\strb}\to\TngFib_{\co,\strb}\\
&\FF\colon\IndxTng_{=,\strb}\to\TngFib_{\co,\strb}
\end{align*}
\end{proposition}

We can finally state the main theorem of this section.

\begin{theorem}[reduced Grothendieck construction for tangent categories]
\label{theorem:reduced-Grothendieck-construction}
The functors of Propositions~\ref{proposition:functor-from-tangent-fibrations-to-indexed-tangent-categories} and~\ref{proposition:functor-from-indexed-tangent-categories-to-tangent-fibrations} form the following adjunctions:
\begin{equation*}
% https://q.uiver.app/#q=WzAsMixbMCwwLCJcXEluZHhUbmdfXFxzdHJiIl0sWzEsMCwiXFxUbmdGaWJfXFxzdHJiIl0sWzAsMSwiXFxGRiIsMCx7ImN1cnZlIjotM31dLFsxLDAsIlxcSUkiLDAseyJjdXJ2ZSI6LTN9XSxbMiwzLCJcXGJvdCIsMSx7InNob3J0ZW4iOnsic291cmNlIjoyMCwidGFyZ2V0IjoyMH0sImxldmVsIjoxLCJzdHlsZSI6eyJib2R5Ijp7Im5hbWUiOiJub25lIn0sImhlYWQiOnsibmFtZSI6Im5vbmUifX19XV0=
\begin{tikzcd}
{\IndxTng_\strb} & {\TngFib_\strb}
\arrow[""{name=0, anchor=center, inner sep=0}, "\FF", bend left=50, from=1-1, to=1-2]
\arrow[""{name=1, anchor=center, inner sep=0}, "\II", bend left=50, from=1-2, to=1-1]
\arrow["\bot"{description}, draw=none, from=0, to=1]
\end{tikzcd}\hfill
% https://q.uiver.app/#q=WzAsMixbMCwwLCJcXEluZHhUbmdfe1xcY28sXFxzdHJifSJdLFsxLDAsIlxcVG5nRmliX3tcXGNvLFxcc3RyYn0iXSxbMCwxLCJcXEZGIiwwLHsiY3VydmUiOi0zfV0sWzEsMCwiXFxJSSIsMCx7ImN1cnZlIjotM31dLFsyLDMsIlxcYm90IiwxLHsic2hvcnRlbiI6eyJzb3VyY2UiOjIwLCJ0YXJnZXQiOjIwfSwibGV2ZWwiOjEsInN0eWxlIjp7ImJvZHkiOnsibmFtZSI6Im5vbmUifSwiaGVhZCI6eyJuYW1lIjoibm9uZSJ9fX1dXQ==
\begin{tikzcd}
{\IndxTng_{\co,\strb}} & {\TngFib_{\co,\strb}}
\arrow[""{name=0, anchor=center, inner sep=0}, "\FF", bend left=50, from=1-1, to=1-2]
\arrow[""{name=1, anchor=center, inner sep=0}, "\II", bend left=50, from=1-2, to=1-1]
\arrow["\bot"{description}, draw=none, from=0, to=1]
\end{tikzcd}
\end{equation*}
Moreover, these adjunctions restrict to strong and strict morphisms of indexed tangent categories by forming two further adjunctions:
\begin{equation*}
% https://q.uiver.app/#q=WzAsMixbMCwwLCJcXEluZHhUbmdfe1xcY29uZyxcXHN0cmJ9Il0sWzEsMCwiXFxUbmdGaWJfe1xcY29uZyxcXHN0cmJ9Il0sWzAsMSwiXFxGRiIsMCx7ImN1cnZlIjotM31dLFsxLDAsIlxcSUkiLDAseyJjdXJ2ZSI6LTN9XSxbMiwzLCJcXGJvdCIsMSx7InNob3J0ZW4iOnsic291cmNlIjoyMCwidGFyZ2V0IjoyMH0sImxldmVsIjoxLCJzdHlsZSI6eyJib2R5Ijp7Im5hbWUiOiJub25lIn0sImhlYWQiOnsibmFtZSI6Im5vbmUifX19XV0=
\begin{tikzcd}
{\IndxTng_{\cong,\strb}} & {\TngFib_{\cong,\strb}}
\arrow[""{name=0, anchor=center, inner sep=0}, "\FF", bend left=50, from=1-1, to=1-2]
\arrow[""{name=1, anchor=center, inner sep=0}, "\II", bend left=50, from=1-2, to=1-1]
\arrow["\bot"{description}, draw=none, from=0, to=1]
\end{tikzcd}\hfill
% https://q.uiver.app/#q=WzAsMixbMCwwLCJcXEluZHhUbmdfez0sXFxzdHJifSJdLFsxLDAsIlxcVG5nRmliX3s9LFxcc3RyYn0iXSxbMCwxLCJcXEZGIiwwLHsiY3VydmUiOi0zfV0sWzEsMCwiXFxJSSIsMCx7ImN1cnZlIjotM31dLFsyLDMsIlxcYm90IiwxLHsic2hvcnRlbiI6eyJzb3VyY2UiOjIwLCJ0YXJnZXQiOjIwfSwibGV2ZWwiOjEsInN0eWxlIjp7ImJvZHkiOnsibmFtZSI6Im5vbmUifSwiaGVhZCI6eyJuYW1lIjoibm9uZSJ9fX1dXQ==
\begin{tikzcd}
{\IndxTng_{=,\strb}} & {\TngFib_{=,\strb}}
\arrow[""{name=0, anchor=center, inner sep=0}, "\FF", bend left=50, from=1-1, to=1-2]
\arrow[""{name=1, anchor=center, inner sep=0}, "\II", bend left=50, from=1-2, to=1-1]
\arrow["\bot"{description}, draw=none, from=0, to=1]
\end{tikzcd}
\end{equation*}
Finally, in each of these adjunctions, the counit is an isomorphism.
\begin{proof}
First, recall that by the Grothendieck construction, given a tangent fibration $\Pi\colon(\X',\TT')\to(\X,\TT)$ there is an isomorphism of fibrations between the underlying fibration of $\Pi$ and the one of $\FF(\II(\Pi))$. Similarly, for an indexed tangent category $(\X,\TT;\I)$, the underlying indexed category is isomorphic to the one underlying $\II(\FF(\I))$. So, we only need to compare the tangent structure of $\El(\X,\TT;\II(\Pi))$ with $\TT'$, and the indexed tangent structure of $\II(\FF(\I))$ with the one of $\I$.
\par Let us start with an indexed tangent category $(\X,\TT;\I)$. It is not hard to see that for an object $A\in\X$ the tangent category associated with $A$ via $\II(\FF(\I))$ is defined as follows. The category is precisely $\X\^A=\I(A)$, since $\FF(\I)^{-1}A=\X\^A$. The tangent bundle functor sends an object $E$ of $\X\^A$ to $z^\*p^\*\T\^AE$, where $\T\^A$ is the tangent bundle functor of $(\X\^A,\TT\^A)=\I(A)$. However, notice that $z^\*p^\*\T\^AE$ is isomorphic to $\T\^AE$ via the natural isomorphism:
\begin{align*}
&z^\*p^\*\T\^AE\xrightarrow{{\Comp}^\-}(p\o z)^\*\T\^AE\xrightarrow{p\o z=\id_\X}\id^\*\T\^AE\xrightarrow{{\Unit}^\-}\T\^AE
\end{align*}
Similarly, a morphism $\psi\colon E\to E'$ corresponds to a morphism $(\id_A,\psi)\colon(A,E)\to(A,E')$ in the tangent category of elements. Thus, the tangent bundle functor of $\II(\FF(\I))(A)$ sends $\psi$ to $z^\*p^\*\psi$. Therefore, the tangent categories $\I(A)$ and $\II(\FF(\I))(A)$ are isomorphic via a natural isomorphism
\begin{align*}
&\II\o\FF\Rightarrow\id_{\IndxTng_\strb}
\end{align*}
which constitutes the counit of the adjunction.
\par The next step is to define the unit. Let us consider a tangent fibration $\Pi\colon(\X',\TT')\to(\X,\TT)$ and let us compare $(\X',\TT')$ with $\El(\X,\TT;\II(\Pi))$. The objects of the latter are pairs $(A,E)$ formed by $A\in\X$ and $E\in\Pi^{-1}A$, morphisms are pairs $(f,\xi\^f)\colon(A,E)\to(A',E')$ formed by a morphism $f\colon A\to A'$ of $\X$ together with a morphism $\xi\^f\colon E\to f^\*E'$ of $\Pi^{-1}A$, that is $\Pi(\xi\^f)=\id_A$. The tangent bundle functor $\hat\T$ sends a pair $(A,E)$ to $(\T A,p^\*z^\*\T'E)$ and a morphism $(f,\xi\^f)\colon(A,E)\to(A',E')$ to $(\T f,\xi\^{\T f})$ where:
\begin{align*}
&\xi\^{\T f}\colon p^\*z^\*\T'E\xrightarrow{p^\*z^\*\T'\xi\^f}p^\*z^\*\T'f^\*E'\xrightarrow{p^\*z^\*\kappa\^f}p^\*z^\*(\T f)^\*\T'E'\cong p^\*f^\*z^\*\T'E'\cong(\T f)^\*p^\*z^\*\T'E'
\end{align*}
However, since $p$ is not invertible, because $z\o p$ is not the identity, the tangent bundle functor $\hat\T$ is not isomorphic to $\T'$. We can still define a natural transformation $\T'\Rightarrow\hat\T$, induced by the universality of the cartesian lift $\varphi\^p$ and $\varphi\^z$, as the unique morphism which makes the following diagram commutes:
\begin{equation*}
\begin{tikzcd}
{\T'E} & {p^\*z^\*\T'E} \\
E & {z^\*\T'E} \\
{\T'E} & {\T'E}
\arrow["{\varphi\^pz^\*\T'}", from=1-2, to=2-2]
\arrow["{\varphi\^z\T'}", from=2-2, to=3-2]
\arrow["p'"', from=1-1, to=2-1]
\arrow["z'"', from=2-1, to=3-1]
\arrow[Rightarrow, no head, from=3-1, to=3-2]
\arrow[dashed, from=1-1, to=1-2]
\end{tikzcd}
\end{equation*}
This natural transformation together with the isomorphism of categories between the category of elements and $\X'$ defines the unit of the adjunction.
\end{proof}
\end{theorem}

Since the base tangent category is preserved by the functors of the adjunctions of Theorem~\ref{theorem:reduced-Grothendieck-construction}, these adjunction restrict to tangent fibrations and indexed tangent categories on a fixed base tangent category.

\begin{corollary}
\label{corollary:reduced-grothendieck-over-fixed-base}
Given a fixed tangent category $(\X,\TT)$, the adjunctions of Theorem~\ref{theorem:reduced-Grothendieck-construction} restricts to the $2$-categories of tangent fibrations and of indexed tangent categories on $(\X,\TT)$:
\begin{equation*}
% https://q.uiver.app/#q=WzAsMixbMCwwLCJcXEluZHhUbmdfXFxzdHJiKFxcWCxcXFRUKSJdLFsxLDAsIlxcVG5nRmliX1xcc3RyYihcXFgsXFxUVCkiXSxbMCwxLCJcXEZGIiwwLHsiY3VydmUiOi0zfV0sWzEsMCwiXFxJSSIsMCx7ImN1cnZlIjotM31dLFsyLDMsIlxcYm90IiwxLHsic2hvcnRlbiI6eyJzb3VyY2UiOjIwLCJ0YXJnZXQiOjIwfSwibGV2ZWwiOjEsInN0eWxlIjp7ImJvZHkiOnsibmFtZSI6Im5vbmUifSwiaGVhZCI6eyJuYW1lIjoibm9uZSJ9fX1dXQ==
\begin{tikzcd}
{\IndxTng_\strb(\X,\TT)} & {\TngFib_\strb(\X,\TT)}
\arrow[""{name=0, anchor=center, inner sep=0}, "\FF", bend left=50, from=1-1, to=1-2]
\arrow[""{name=1, anchor=center, inner sep=0}, "\II", bend left=50, from=1-2, to=1-1]
\arrow["\bot"{description}, draw=none, from=0, to=1]
\end{tikzcd}\hfill
% https://q.uiver.app/#q=WzAsMixbMCwwLCJcXEluZHhUbmdfe1xcY28sXFxzdHJifShcXFgsXFxUVCkiXSxbMSwwLCJcXFRuZ0ZpYl97XFxjbyxcXHN0cmJ9KFxcWCxcXFRUKSJdLFswLDEsIlxcRkYiLDAseyJjdXJ2ZSI6LTN9XSxbMSwwLCJcXElJIiwwLHsiY3VydmUiOi0zfV0sWzIsMywiXFxib3QiLDEseyJzaG9ydGVuIjp7InNvdXJjZSI6MjAsInRhcmdldCI6MjB9LCJsZXZlbCI6MSwic3R5bGUiOnsiYm9keSI6eyJuYW1lIjoibm9uZSJ9LCJoZWFkIjp7Im5hbWUiOiJub25lIn19fV1d
\begin{tikzcd}
{\IndxTng_{\co,\strb}(\X,\TT)} & {\TngFib_{\co,\strb}(\X,\TT)}
\arrow[""{name=0, anchor=center, inner sep=0}, "\FF", bend left=50, from=1-1, to=1-2]
\arrow[""{name=1, anchor=center, inner sep=0}, "\II", bend left=50, from=1-2, to=1-1]
\arrow["\bot"{description}, draw=none, from=0, to=1]
\end{tikzcd}
\end{equation*}
\begin{equation*}
% https://q.uiver.app/#q=WzAsMixbMCwwLCJcXEluZHhUbmdfe1xcY29uZyxcXHN0cmJ9KFxcWCxcXFRUKSJdLFsxLDAsIlxcVG5nRmliX3tcXGNvbmcsXFxzdHJifShcXFgsXFxUVCkiXSxbMCwxLCJcXEZGIiwwLHsiY3VydmUiOi0zfV0sWzEsMCwiXFxJSSIsMCx7ImN1cnZlIjotM31dLFsyLDMsIlxcYm90IiwxLHsic2hvcnRlbiI6eyJzb3VyY2UiOjIwLCJ0YXJnZXQiOjIwfSwibGV2ZWwiOjEsInN0eWxlIjp7ImJvZHkiOnsibmFtZSI6Im5vbmUifSwiaGVhZCI6eyJuYW1lIjoibm9uZSJ9fX1dXQ==
\begin{tikzcd}
{\IndxTng_{\cong,\strb}(\X,\TT)} & {\TngFib_{\cong,\strb}(\X,\TT)}
\arrow[""{name=0, anchor=center, inner sep=0}, "\FF", bend left=50, from=1-1, to=1-2]
\arrow[""{name=1, anchor=center, inner sep=0}, "\II", bend left=50, from=1-2, to=1-1]
\arrow["\bot"{description}, draw=none, from=0, to=1]
\end{tikzcd}\hfill
% https://q.uiver.app/#q=WzAsMixbMCwwLCJcXEluZHhUbmdfez0sXFxzdHJifShcXFgsXFxUVCkiXSxbMSwwLCJcXFRuZ0ZpYl97PSxcXHN0cmJ9KFxcWCxcXFRUKSJdLFswLDEsIlxcRkYiLDAseyJjdXJ2ZSI6LTN9XSxbMSwwLCJcXElJIiwwLHsiY3VydmUiOi0zfV0sWzIsMywiXFxib3QiLDEseyJzaG9ydGVuIjp7InNvdXJjZSI6MjAsInRhcmdldCI6MjB9LCJsZXZlbCI6MSwic3R5bGUiOnsiYm9keSI6eyJuYW1lIjoibm9uZSJ9LCJoZWFkIjp7Im5hbWUiOiJub25lIn19fV1d
\begin{tikzcd}
{\IndxTng_{=,\strb}(\X,\TT)} & {\TngFib_{=,\strb}(\X,\TT)}
\arrow[""{name=0, anchor=center, inner sep=0}, "\FF", bend left=50, from=1-1, to=1-2]
\arrow[""{name=1, anchor=center, inner sep=0}, "\II", bend left=50, from=1-2, to=1-1]
\arrow["\bot"{description}, draw=none, from=0, to=1]
\end{tikzcd}
\end{equation*}
\end{corollary}

\begin{corollary}
\label{corollary:reduced-grothendieck-base-trivial}
When the base tangent category $(\X,\TT)$ has a trivial tangent structure, the adjunctions of Corollary~\ref{corollary:reduced-grothendieck-over-fixed-base} become equivalences.
\begin{proof}
In the proof of Theorem~\ref{theorem:reduced-Grothendieck-construction}, we showed that the reason why the unit of the adjunction is not an isomorphism is that, in general, $z\o p$ is not the identity, where $z$ and $p$ are the zero morphism and the projection of the base tangent category of a tangent fibration. However, when the base tangent structure is trivial, this condition holds, since $p$ and $z$ are the identities.
\end{proof}
\end{corollary}

%__________________________________________________________________________
\subsection{Tangent fibrations as pseudoalgebras}
\label{subsection:pseudoalgebras}
In~\cite{street:fibrations}, Street presented (cloven) fibrations as pseudoalgebras of a $2$-monad and he employed this characterization to extend the notion of a fibration in the context of $2$-category theory. Since tangent categories form a $2$-category, it is natural to wonder whether or not Street's definition of fibrations applied to the context of tangent categories provides the same notion of a tangent fibration of Definition~\ref{definition:tangent-fibration}.
\par This section is dedicated to exploring this question. Let us start by recalling Street's construction. Consider the $2$-category $\Cat$ of categories and let $\X$ be a fixed category. On the slice $2$-category $\Cat/\X$, one can define a $2$-functor $M$ which sends a functor $\Pi\colon\X'\to\X$ to $M(\Pi)\colon\X\downarrow\Pi\to\X$, where $\X\downarrow\Pi$ denotes the comma category of $\Pi$ and the functor $M(\Pi)$ sends each $(f\colon A\to\Pi(E))\in\X\downarrow\Pi$ to $A\in\X$.
\par $M\colon\Cat/\X\to\Cat/\X$ comes equipped with a $2$-monad structure and cloven fibrations are precisely pseudoalgebras of $M$. Let us recall the notion of a pseudoalgebra for a $2$-monad.

\begin{definition}
\label{definition:pseudoalgebra}
A \textbf{pseudoalgebra} of a $2$-monad $M\colon\CC\to\CC$ consists of an object $A$ of $\CC$ together with a $1$-morphism $c\colon MA\to A$ of $\CC$ and two $2$-morphisms
\begin{equation*}
% https://q.uiver.app/#q=WzAsNCxbMCwwLCJBIl0sWzEsMSwiQSJdLFswLDEsIk1BIl0sWzEsMCwiQSJdLFswLDIsIlxcZXRhIiwyXSxbMiwxLCJjIiwyXSxbMCwzLCIiLDAseyJsZXZlbCI6Miwic3R5bGUiOnsiaGVhZCI6eyJuYW1lIjoibm9uZSJ9fX1dLFszLDEsIiIsMCx7ImxldmVsIjoyLCJzdHlsZSI6eyJoZWFkIjp7Im5hbWUiOiJub25lIn19fV0sWzMsMiwiXFx6ZXRhIiwwLHsic2hvcnRlbiI6eyJzb3VyY2UiOjIwfSwibGV2ZWwiOjJ9XV0=
\begin{tikzcd}
A & A \\
MA & A
\arrow[Rightarrow, no head, from=1-1, to=1-2]
\arrow["\eta"', from=1-1, to=2-1]
\arrow["\zeta", shorten <=4pt, Rightarrow, from=1-2, to=2-1]
\arrow[Rightarrow, no head, from=1-2, to=2-2]
\arrow["c"', from=2-1, to=2-2]
\end{tikzcd}\hfill
% https://q.uiver.app/#q=WzAsNCxbMCwwLCJNXjJBIl0sWzEsMCwiTUEiXSxbMCwxLCJNQSJdLFsxLDEsIkEiXSxbMCwxLCJNYyJdLFsxLDMsImMiXSxbMCwyLCJcXGdhbW1hIiwyXSxbMiwzLCJjIiwyXSxbMSwyLCJcXHRoZXRhIiwwLHsibGV2ZWwiOjJ9XV0=
\begin{tikzcd}
{M^2A} & MA \\
MA & A
\arrow["Mc", from=1-1, to=1-2]
\arrow["\gamma"', from=1-1, to=2-1]
\arrow["\theta", Rightarrow, from=1-2, to=2-1]
\arrow["c", from=1-2, to=2-2]
\arrow["c"', from=2-1, to=2-2]
\end{tikzcd}
\end{equation*}
where $\eta\colon\id_\CC\Rightarrow M$ and $\gamma\colon M^2\Rightarrow M$ are the unit and the multiplication of the $2$-monad $M$, respectively. Moreover, $\zeta$ and $\theta$ satisfy some coherence conditions with the monad structure (see~\cite[Section~2]{street:fibrations}).
\end{definition}

Concretely, a pseudoalgebra of the $2$-monad $M\colon\Cat/\X\to\Cat/\X$ is a functor $\Pi\colon\X'\to\X$ together with a $2$-functor $\trans\colon\X\downarrow\Pi\to\X'$ which sends each $f\colon A\to\Pi(E)$ to an object $f^\*E$ of $\X'$. Moreover, the $2$-functor $\trans$ sends the square diagram
\begin{equation*}
% https://q.uiver.app/#q=WzAsNCxbMCwwLCJBIl0sWzEsMCwiXFxQaShFKSJdLFswLDEsIlxcUGkoRSkiXSxbMSwxLCJcXFBpKEUpIl0sWzAsMSwiZiJdLFswLDIsImYiLDJdLFsyLDMsIiIsMix7ImxldmVsIjoyLCJzdHlsZSI6eyJoZWFkIjp7Im5hbWUiOiJub25lIn19fV0sWzEsMywiXFxQaShcXGlkX0UpIl1d
\begin{tikzcd}
A & {\Pi(E)} \\
{\Pi(E)} & {\Pi(E)}
\arrow["f", from=1-1, to=1-2]
\arrow["f"', from=1-1, to=2-1]
\arrow["{\Pi(\id_E)}", from=1-2, to=2-2]
\arrow[Rightarrow, no head, from=2-1, to=2-2]
\end{tikzcd}
\end{equation*}
to a morphism $f^\*E\to E$ of $\X'$, which corresponds to the cartesian lift of $f\colon A\to\Pi(E)$ on $E$. Finally, the $2$-morphisms $\zeta$ and $\theta$ provide the unitor and the compositor.
\par For an arbitrary $2$-category $\CC$, the $2$-monad $M$ is replaced with the $2$-monad which sends each object $\Pi\colon\X'\to\X$ of the slice $2$-category $\CC/\X$ to $M(\Pi)\colon\X\downarrow\Pi\to\X$, where $\X\downarrow\Pi$ denotes the comma object of $\Pi$, provided the existence of such comma objects.\newline
\par To see whether or not pseudoalgebras of the $2$-monad $M$ in the $2$-category of tangent categories correspond to the notion of a tangent fibration of Definition~\ref{definition:tangent-fibration}, first we need to establish if such a $2$-category has comma objects. In particular, we need to see whether or not the comma category $\X\downarrow\Pi$ carries a tangent structure when $(\Pi,\alpha)$ is a tangent morphism. In general, this is not the case, however, by considering a subclass of tangent morphisms we obtain the following characterization.

\begin{lemma}
\label{lemma:comma-tangent-category}
Let $(\Pi,\alpha)\colon(\X',\TT')\to(\X,\TT)$ be a strong tangent morphism. Then, the comma category $\X\downarrow\Pi$ carries a tangent structure whose tangent bundle functor sends each $f\colon A\to\Pi(E)$ to:
\begin{align*}
&\T A\xrightarrow{\T f}\T\Pi(E)\xrightarrow{\alpha^{-1}}\Pi\T'(E)
\end{align*}
Moreover, this tangent category is the comma object of the $2$-category $\TngCat_\cong$ of tangent categories and strong tangent morphisms which preserve the $n$-fold pullbacks of the projection and the universality of the vertical lift.
\begin{proof}
Let us briefly sketch the proof. To show that the comma category $\X\downarrow\Pi$ comes equipped with a tangent structure, one uses the compatibilities of the strong distributive law $\alpha$ with the tangent structures and that the functor $\Pi$ preserves the $n$-fold pullbacks of the projection along itself and the universality of the vertical lift.
\par To prove that such a construction defines comma objects, let us consider another tangent category $(\X'',\TT'')$ together with two strong tangent morphisms $(G,\beta)\colon(\X'',\TT'')\to(\X,\TT)$ and $(H,\gamma)\colon(\X'',\TT'')\to(\X',\TT')$ and a tangent natural transformation:
\begin{equation*}
% https://q.uiver.app/#q=WzAsNCxbMCwwLCIoXFxYJycsXFxUVCcnKSJdLFswLDEsIihcXFgsXFxUVCkiXSxbMSwxLCIoXFxYLFxcVFQpIl0sWzEsMCwiKFxcWCcsXFxUVCcpIl0sWzMsMiwiKFxcUGksXFxhbHBoYSkiXSxbMSwyLCIiLDIseyJsZXZlbCI6Miwic3R5bGUiOnsiaGVhZCI6eyJuYW1lIjoibm9uZSJ9fX1dLFswLDEsIihHLFxcYmV0YSkiLDJdLFswLDMsIihILFxcZ2FtbWEpIl0sWzMsMSwiXFx2YXJwaGkiLDEseyJsZXZlbCI6Mn1dXQ==
\begin{tikzcd}
{(\X'',\TT'')} & {(\X',\TT')} \\
{(\X,\TT)} & {(\X,\TT)}
\arrow["{(H,\gamma)}", from=1-1, to=1-2]
\arrow["{(G,\beta)}"', from=1-1, to=2-1]
\arrow["\varphi"{description}, Rightarrow, from=1-2, to=2-1]
\arrow["{(\Pi,\alpha)}", from=1-2, to=2-2]
\arrow[Rightarrow, no head, from=2-1, to=2-2]
\end{tikzcd}
\end{equation*}
Let us consider the functor
\begin{align*}
K\colon&\X''\to\X\downarrow\Pi\\
&A\mapsto(\varphi_A\colon GA\to\Pi(HA))\\
&(f\colon A\to B)\mapsto((\Pi(Hf),Gf)\colon\varphi_A\to\varphi_B)
\end{align*}
together with the distributive law $\delta\colon K\o\T''\Rightarrow\bar\T\o K$, so defined:
\begin{equation*}
% https://q.uiver.app/#q=WzAsNSxbMCwwLCJHXFxUJydBIl0sWzAsMiwiXFxQaSBIXFxUJydBIl0sWzEsMCwiXFxUIEdBIl0sWzEsMSwiXFxUXFxQaSBIQSJdLFsxLDIsIlxcUGlcXFQnSEEiXSxbMCwxLCJcXHZhcnBoaV97XFxUJydBfSIsMl0sWzIsMywiXFxUXFx2YXJwaGlfQSJdLFszLDQsIlxcYWxwaGFeXFwtSCJdLFsxLDQsIlxcUGlcXGdhbW1hIiwyXSxbMCwyLCJcXGJldGEiXSxbMCwxLCJLXFxUJydBIiwyLHsib2Zmc2V0Ijo0LCJjdXJ2ZSI6M31dLFsyLDQsIlxcVCBLQSIsMCx7Im9mZnNldCI6LTQsImN1cnZlIjotM31dLFs1LDMsIlxcZGVsdGEiLDAseyJzaG9ydGVuIjp7InNvdXJjZSI6MTB9fV1d
\begin{tikzcd}
{G\T''A} & {\T GA} \\
& {\T\Pi HA} \\
{\Pi H\T''A} & {\Pi\T'HA}
\arrow["\beta", from=1-1, to=1-2]
\arrow[""{name=0, anchor=center, inner sep=0}, "{\varphi_{\T''A}}"', from=1-1, to=3-1]
\arrow["{K\T''A}"', shift right=4, bend right=50, from=1-1, to=3-1]
\arrow["{\T\varphi_A}", from=1-2, to=2-2]
\arrow["{\T KA}", shift left=4, bend left=50, from=1-2, to=3-2]
\arrow["{\alpha^\-H}", from=2-2, to=3-2]
\arrow["{\Pi\gamma}"', from=3-1, to=3-2]
\arrow["\delta", shorten <=3pt, Rightarrow, from=0, to=2-2]
\end{tikzcd}
\end{equation*}
One can prove that this strong tangent morphism preserves the $n$-fold pullbacks of the projection along itself, the universality of the vertical lift, and that it is the unique strong tangent morphism such that
\begin{equation*}
% https://q.uiver.app/#q=WzAsMyxbMCwwLCIoXFxYJycsXFxUVCcnKSJdLFsxLDEsIihcXFgsXFxUVCkiXSxbMSwwLCIoXFxYLFxcVFQpXFxkb3duYXJyb3coXFxQaSxcXGFscGhhKSJdLFswLDEsIihHLFxcYmV0YSkiLDJdLFsyLDEsIlxccGlfMSJdLFswLDIsIihLLFxcZGVsdGEpIiwwLHsic3R5bGUiOnsiYm9keSI6eyJuYW1lIjoiZGFzaGVkIn19fV1d
\begin{tikzcd}
{(\X'',\TT'')} & {(\X,\TT)\downarrow(\Pi,\alpha)} \\
& {(\X,\TT)}
\arrow["{(K,\delta)}", dashed, from=1-1, to=1-2]
\arrow["{(G,\beta)}"', from=1-1, to=2-2]
\arrow["{\pi_1}", from=1-2, to=2-2]
\end{tikzcd}\hfill
% https://q.uiver.app/#q=WzAsMyxbMCwwLCIoXFxYJycsXFxUVCcnKSJdLFsxLDEsIihcXFgnLFxcVFQnKSJdLFsxLDAsIihcXFgsXFxUVClcXGRvd25hcnJvdyhcXFBpLFxcYWxwaGEpIl0sWzAsMSwiKEgsXFxnYW1tYSkiLDJdLFsyLDEsIlxccGlfMiJdLFswLDIsIihLLFxcZGVsdGEpIiwwLHsic3R5bGUiOnsiYm9keSI6eyJuYW1lIjoiZGFzaGVkIn19fV1d
\begin{tikzcd}
{(\X'',\TT'')} & {(\X,\TT)\downarrow(\Pi,\alpha)} \\
& {(\X',\TT')}
\arrow["{(K,\delta)}", dashed, from=1-1, to=1-2]
\arrow["{(H,\gamma)}"', from=1-1, to=2-2]
\arrow["{\pi_2}", from=1-2, to=2-2]
\end{tikzcd}
\end{equation*}
where $\pi_1$ and $\pi_2$ denotes the projections of the comma object $\X\downarrow\Pi$, which are strict tangent morphisms, and such that:
\begin{align*}
&\varphi=\psi_{(K,\delta)}
\end{align*}
where:
\begin{equation*}
% https://q.uiver.app/#q=WzAsNCxbMCwxLCIoXFxYLFxcVFQpIl0sWzAsMCwiKFxcWCxcXFRUKVxcZG93bmFycm93KFxcUGksXFxhbHBoYSkiXSxbMSwwLCIoXFxYJyxcXFRUJykiXSxbMSwxLCIoXFxYLFxcVFQpIl0sWzEsMCwiXFxwaV8xIiwyXSxbMSwyLCJcXHBpXzIiXSxbMiwwLCJcXHBzaSIsMSx7ImxldmVsIjoyfV0sWzIsMywiKFxcUGksXFxhbHBoYSkiXSxbMCwzLCIiLDIseyJsZXZlbCI6Miwic3R5bGUiOnsiaGVhZCI6eyJuYW1lIjoibm9uZSJ9fX1dXQ==
\begin{tikzcd}
{(\X,\TT)\downarrow(\Pi,\alpha)} & {(\X',\TT')} \\
{(\X,\TT)} & {(\X,\TT)}
\arrow["{\pi_2}", from=1-1, to=1-2]
\arrow["{\pi_1}"', from=1-1, to=2-1]
\arrow["\psi"{description}, Rightarrow, from=1-2, to=2-1]
\arrow["{(\Pi,\alpha)}", from=1-2, to=2-2]
\arrow[Rightarrow, no head, from=2-1, to=2-2]
\end{tikzcd}
\end{equation*}
is the tangent natural transformation which picks out the object of the comma category.
\end{proof}
\end{lemma}

Thanks to Lemma~\ref{lemma:comma-tangent-category}, we can define the $2$-monad $M\colon\TngCat_\cong/(\X,\TT)\to\TngCat_\cong/(\X,\TT)$ which sends a strong tangent morphism $(\Pi,\alpha)\colon(\X',\TT')\to(\X,\TT)$ to $M(\Pi,\alpha)\colon(\X,\TT)\downarrow(\Pi,\alpha)\to(\X,\TT)$.

\begin{theorem}
\label{theorem:tangent-fibrations-are-pseudoalgebras}
Tangent fibrations on a given tangent category $(\X,\TT)$ are precisely pseudoalgebras of the $2$-monad $M\colon\TngCat_\cong/(\X,\TT)\to\TngCat_\cong/(\X,\TT)$ whose underlying tangent morphism is strict.
\begin{proof}
A pseudoalgebra of $M$ consists of a strong tangent morphism $(\Pi,\alpha)\colon(\X',\TT')\to(\X,\TT)$ together with another strong tangent morphism $(\trans,\kappa)\colon(\X,\TT)\downarrow(\Pi,\alpha)\to(\X',\TT')$ and two natural transformations $\zeta$ and $\theta$. In particular, following the same argument used above for when the ambient $2$-category is $\Cat$, $\Pi\colon\X'\to\X$ together with the functor $\trans\colon\X\downarrow\Pi\to\X'$ defines a cloven fibration. Moreover, the strong distributive law $\kappa\colon\trans\o\Bar\T\Rightarrow\T'\o\trans$ defines the distributors of the tangent bundle functors, where $\Bar\T$ denotes the tangent bundle functor of $(\X,\TT)\downarrow(\Pi,\alpha)$.
\end{proof}
\end{theorem}

There is only one important difference between tangent fibrations and pseudoalgebras of $M$: in the latter, the underlying tangent morphism $(\Pi,\alpha)\colon(\X',\TT)'\to(\X,\TT)$ is strong but not necessarily strict. One could have tried to consider the $2$-category $\TngCat_=$ of tangent categories and strict tangent morphism instead of $\TngCat_\cong$. However, in this context, the corresponding pseudoalgebras are instead fibrations between tangent categories whose tangent bundle functors strictly (not strongly) preserve the cartesian lifts.

\begin{remark}
\label{remark:weaker-tangent-fibrations}
The notion of an internal fibration of the $2$-category $\TngCat_\cong$ suggests to investigate a weaker version of the notion of a tangent fibration, whose underlying functor only preserves strongly (and not strictly) the tangent structures. In future work, we are interested in exploring this definition and looking for examples.
\end{remark}

%__________________________________________________________________________
%__________________________________________________________________________

\section{Tangent objects: a formal approach to tangent categories}
\label{section:tangent-objects}
Theorem~\ref{theorem:reduced-Grothendieck-construction} provides an adjunction between tangent fibrations and indexed tangent categories; however, the unit of this adjunction is not, in general, an isomorphism. The issue is due to the loss of information caused by pulling back along the zero morphism in order to reposition the tangent bundle functor into the right fibre.
\par In order to prove a genuine Grothendieck equivalence, we should avoid pulling back along the zero morphism; this would yield a collection of functors of type $\T\^A\colon\Pi^{-1}(A)\to\Pi^{-1}(\T A)$ indexed by the objects of the base category $\X$. This intuition suggests looking at a ``global tangent structure'' over the entire indexed category, instead of a ``local tangent structure'' on each fibre.
\par In this section, we want to make precise the idea of a \textit{tangent structure} for indexed categories. Our approach is far more general: we aim to introduce the notion of a tangent structure on \emph{any} object of a given strict $2$-category. This approach is inspired by the formal theory of Street presented in~\cite{street:formal-theory-monads}, in which the formal theory of monads on a given $2$-category was explored. Here, we focus our efforts on introducing the main concept; in future work, we intend to explore further this idea.\newline

\par Leung in his Ph.D. thesis~\cite{leung:weil-algebras} proposed a simple and effective classification of the tangent structures for a given category $\X$. In particular, he showed that tangent structures $\TT$ for $\X$ are in one-to-one correspondence with strict monoidal functors:
\begin{align*}
&\Leung[\TT]\colon\Weil\to\End(\X)
\end{align*}
from (a subcategory of) the monoidal category of Weil algebras to the monoidal category $\End(\X)\=\CC(\X,\X)$ of endofunctors over the category $\X$, satisfying extra conditions. A Weil algebra is a commutative and unital $\N$-algebra, obtained by quotienting the $\N$-algebra $\N[x_1,\dots,x_n]$ of $\N$-polynomials in $n$ variables by an ideal generated by monomials of order $2$. More precisely, $\Weil$ denotes the monoidal category generated by the Weil algebras $W^n\=\N[x_1\,x_n]/(x_ix_j,i\leq j)$, for positive integers $n$. As shown by Leung, the category $\Weil$ comes equipped with a tangent structure whose tangent bundle functor sends a Weil algebra $A$ to $A\otimes W$, where $W\=W^1$ and whose natural transformations are induced by the following morphisms:

\begin{description}
\item[projection] The projection $p\colon W\to\N$, which sends the generator $x$ of $W$ to $0$;

\item[zero morphism] The zero morphism $z\colon\N\to W$, which sends an integer to itself;

\item[sum morphism] The sum morphism $s\colon W^2\to W$, which sends the two generators $x_1$ and $x_2$ to the unique generator $x$;

\item[vertical lift] The vertical lift $l\colon W\to W\x W$, which sends the generator $x$ to $x\x y$;

\item[canonical flip] The canonical flip $c\colon W\x W\to W\x W$, which sends the generator $x$ of the left $W$ to the generator $y$ of the right $W$, and vice versa, i.e., $y$ to $x$.
\end{description}
Leung's classification establishes that a tangent structure $\TT$ over a category $\X$ is precisely given by a strict monoidal functor $\Leung[\TT]$ which sends the generators $W^n$ of $\Weil$ to the functors $\T_n$, and the morphisms listed above to the homonymous natural transformations of $\TT$. In particular, the tangent bundle functor is $\T=\Leung[\TT](W)$, the double tangent bundle functor is $\T^2=\Leung[\TT](W\x W)$, the projection is $p=\Leung[\TT](p)\colon\T=\Leung[\TT](W)\to \Leung[\TT](\N)=\id_\X$, etcetera.
\par In this section we explore a generalization of this classification which leads to a simple but important tool for our discussion: the concept of a tangent object.\newline
\par Let $\CC$ be a fixed strict $2$-category, that is a category enriched over $\Cat$. In future work, we would like to explore weaker versions of this concept, but for now, let us focus on the strict case.
\par Before defining a tangent object, we first need to introduce an important technical definition.

\begin{definition}
\label{definition:pointwise-limits}
Given a strict $2$-category $\CC$ and two objects $\X$ and $\Y$ of $\CC$, a limit in the category $\CC(\X,\Y)$ is \textbf{pointwise} when it is preserved by all functors $\CC(f,\Y)\colon\CC(\X,\Y)\to\CC(\X',\Y)$ for every $1$-morphism $f\colon\X'\to\X$ in $\CC$.
\end{definition}

\begin{remark}
\label{remark:pointwise-limits}
We would like to thank Rory Lucyshyn-Wright for suggesting this assumption and pointing out its importance for tangent objects in an informal discussion with the author. This aspect was missing in the original definition provided by the author.
\end{remark}

When the $2$-category $\CC$ is the $2$-category $\Cat$ of categories, pointwise limits of $\CC(\X,\Y)$ are those limit diagrams in the category of functors of type $F\colon\X\to\Y$ that are preserved by the evaluation functor. Concretely, this means that, for an object $X$ of $\X$, and a diagram $D\colon\X_0\to\CC(\X,\Y)$, the functor $\lim D\colon\X\to\Y$ evaluated at $X$ is isomorphic to the object $\lim D(X)$ of $\Y$, where $D(X)$ represents the diagram $\X_0\to\Y$ obtained by evaluating each functor $D_A\colon\X\to\Y$, corresponding to each $A$ of $\X_0$, at $X$.
\par When the target category $\Y$ has all finite limits, then so does the category of functors $\CC(\X,\Y)$ and each limit is pointwise. However, when the target category is not known to be finitely complete, there is no guarantee that the limits of $\CC(\X,\Y)$ will be pointwise. A counterexample can be found in~\cite[Section~3.3]{kelly:enriched-cats}.
\par Unfortunately, in tangent category theory, the requirement of the existence of limits is a subtle matter since in differential geometry not every pair of morphisms admits a pullback. In particular, a tangent category cannot be required to be finitely complete since this would rule out one of the main examples of a tangent category. Consequently, in order to make our definition of tangent objects compatible with the usual notion of a tangent category when $\CC$ is assumed to be the $2$-category of categories, we need to require the limits involved in the definition of a tangent object to be pointwise.

\begin{definition}
\label{definition:tangent-object}
A \textbf{tangent object} in a $2$-category $\CC$ is an object $\X$ of $\CC$ equipped with a \textbf{tangent structure} $\TT$, which consists of a \textbf{Leung monoidal functor}, which is a strict monoidal functor $\Leung[\TT]\colon\Weil\to\End(\X)$ from the monoidal category of Weil algebras to the monoidal category of endomorphisms over $\X$ in $\CC$, satisfying the following two universal conditions:
\begin{itemize}
\item $\Leung[\TT]$ preserves the \textbf{foundational pullbacks}, which are pullbacks of the form:
\begin{equation*}
\begin{tikzcd}
{A\x(B\times C)} & {A\x C} \\
{A\x B} & A
\arrow["{A\x\pi_2}", from=1-1, to=1-2]
\arrow["{A\x\pi_1}"', from=1-1, to=2-1]
\arrow["{A\x p}"', from=2-1, to=2-2]
\arrow["{A\x p}", from=1-2, to=2-2]
\arrow["\lrcorner"{anchor=center, pos=0.125}, draw=none, from=1-1, to=2-2]
\end{tikzcd}
\end{equation*}
for all $A,B,C\in\Weil$ (cf.~\cite[Definition~3.17]{leung:weil-algebras}), where, for the sake of simplicity, we omitted the unitors in the diagram. Moreover, these pullbacks are pointwise limits;
\item $\Leung[\TT]$ preserves the universality of the vertical lift, i.e., the pullback diagram:
\begin{equation*}
% https://q.uiver.app/#q=WzAsNCxbMCwwLCJXXjIiXSxbMSwwLCJXXFx4IFciXSxbMSwxLCJXIl0sWzAsMSwiXFxOIl0sWzAsMSwiXFxudSJdLFsxLDIsIldcXHggcCJdLFszLDIsInoiLDJdLFswLDMsIlxccGlfMXAiLDJdLFswLDIsIiIsMSx7InN0eWxlIjp7Im5hbWUiOiJjb3JuZXIifX1dXQ==
\begin{tikzcd}
{W^2} & {W\x W} \\
\N & W
\arrow["\phi", from=1-1, to=1-2]
\arrow["{W\x p}", from=1-2, to=2-2]
\arrow["z"', from=2-1, to=2-2]
\arrow["{\pi_1p}"', from=1-1, to=2-1]
\arrow["\lrcorner"{anchor=center, pos=0.125}, draw=none, from=1-1, to=2-2]
\end{tikzcd}
\end{equation*}
where $\phi\=\<z\x W,l\>(W\x s)$ and $\pi_1\colon W^2\to W$ sends $x_1$ to $x$ and $x_2$ to zero, where, for the sake of simplicity, we omitted the unitors in the diagram. Moreover, this pullback is a pointwise limit.
\end{itemize}
\end{definition}

\begin{remark}
\label{remark:universality-vertical-lift-weil}
In Leung's original result, the universality of the vertical lift of Definition~\ref{definition:tangent-object}, is replaced with the universality of an equalizer. However, Cockett and Cruttwell proved in~\cite[Lemma~2.12]{cockett:tangent-cats} that the universality of the pullback diagram of Definition~\ref{definition:tangent-object} is equivalent to the universality of the equalizer diagram proposed by Leung. To stay consistent with the rest of the paper, we adopted the pullback version.
\end{remark}

\begin{remark}
\label{remark:negatives-weil}
To classify tangent structures with negatives one can replace the rig $\N$ with the ring $\Z$ in the definition of a Weil algebra and then introduce the negation as follows:
\begin{description}
\item[negation] The negation $n\colon W\to W$ sends the generator $x$ to $-x$.
\end{description}
Thus, Leung's classification extends as follows: tangent structures with negatives $\TT$ over a category $\X$ are in one-to-one correspondence with strict monoidal functors $\Leung^-[\TT]\colon\Weil^-\to\End(\X)$ preserving foundational pullbacks and the universality of the vertical lift, where $\Weil^-$ is the category of Weil algebras over the ring $\Z$.
\end{remark}

Thanks to Remark~\ref{remark:negatives-weil}, we can also define a tangent object with negatives as follows.

\begin{definition}
\label{definition:tangent-object-with-negatives}
A tangent object \textbf{with negatives} in a $2$-category $\CC$ is an object $\X$ of $\CC$ together with a tangent structure with negatives $\TT$, which consists of a \textbf{Leung monoidal functor with negatives}, which is a strict monoidal functor $\Leung^-[\TT]\colon\Weil^-\to\End(\X)$, preserving foundational pullbacks and the universality of the vertical lift as in Definition~\ref{definition:tangent-object}.
\end{definition}

Using a similar strategy to the one used by Leung to classify tangent structure on a given category, we can unwrap Definitions~\ref{definition:tangent-object} and~\ref{definition:tangent-object-with-negatives} to have a more concrete understanding of these notions. Let us start by introducing a useful concept.

\begin{definition}
\label{definition:additive-bundle-object}
An \textbf{additive bundle object} over an object $\X$ of a $2$-category $\CC$ is an additive bundle in the category $\End(\X)$ of endomorphisms of $\X$. Concretely, it consists of an object $\X\in\CC$, two $1$-endomorphisms $B\colon\X\to\X$ and $E\colon\X\to\X$, together with a $2$-morphism $q\colon E\Rightarrow B$, called the \textbf{projection}, a $2$-morphism $z_q\colon B\Rightarrow E$ called the \textbf{zero morphism}, and a $2$-morphism $s_q\colon E_2\Rightarrow E$, called the \textbf{sum morphism}, satisfying the following properties:
\begin{itemize}
\item $z_q$ is a section of $q$:
\begin{equation*}
% https://q.uiver.app/#q=WzAsMyxbMCwwLCJFIl0sWzEsMSwiQSJdLFswLDEsIkEiXSxbMiwwLCJ6X3EiXSxbMCwxLCJxIl0sWzIsMSwiIiwyLHsibGV2ZWwiOjIsInN0eWxlIjp7ImhlYWQiOnsibmFtZSI6Im5vbmUifX19XV0=
\begin{tikzcd}
E \\
B & B
\arrow["q", from=1-1, to=2-2]
\arrow["{z_q}", from=2-1, to=1-1]
\arrow[Rightarrow, no head, from=2-1, to=2-2]
\end{tikzcd}
\end{equation*}

\item $n$-fold pullbacks: for any positive integer $n$, the $n$-fold pullback of the projection $q$ along itself exists in the category $\End(\X)$ of endomorphisms of $\X$, is a pointwise limit, and is preserved by each $E^m\=E\o\dots\o E$, for every positive integer $m$. The $k$-th projection $\pi_k\colon E_n\Rightarrow E$ is denoted by $\pi_k$;

\item $s_q$ is a bundle morphism:
\begin{equation*}
% https://q.uiver.app/#q=WzAsNCxbMSwwLCJFIl0sWzEsMSwiQSJdLFswLDAsIkVfMiJdLFswLDEsIkUiXSxbMCwxLCJxIl0sWzIsMywiXFxwaV8xIiwyXSxbMywxLCJxIiwyXSxbMiwwLCJzX3EiXV0=
\begin{tikzcd}
{E_2} & E \\
E & B
\arrow["{s_q}", from=1-1, to=1-2]
\arrow["{\pi_1}"', from=1-1, to=2-1]
\arrow["q", from=1-2, to=2-2]
\arrow["q"', from=2-1, to=2-2]
\end{tikzcd}\hfill
% https://q.uiver.app/#q=WzAsNCxbMSwwLCJFIl0sWzEsMSwiQSJdLFswLDAsIkVfMiJdLFswLDEsIkUiXSxbMCwxLCJxIl0sWzIsMywiXFxwaV8yIiwyXSxbMywxLCJxIiwyXSxbMiwwLCJzX3EiXV0=
\begin{tikzcd}
{E_2} & E \\
E & B
\arrow["{s_q}", from=1-1, to=1-2]
\arrow["{\pi_2}"', from=1-1, to=2-1]
\arrow["q", from=1-2, to=2-2]
\arrow["q"', from=2-1, to=2-2]
\end{tikzcd}
\end{equation*}

\item Associativity:
\begin{equation*}
% https://q.uiver.app/#q=WzAsNCxbMCwwLCJFXzMiXSxbMSwwLCJFXzIiXSxbMSwxLCJFIl0sWzAsMSwiRV8yIl0sWzAsMSwiXFxpZF9FXFx0aW1lc19Bc19xIl0sWzEsMiwic19xIl0sWzMsMiwic19xIiwyXSxbMCwzLCJzX3FcXHRpbWVzX0FcXGlkX0UiLDJdXQ==
\begin{tikzcd}
{E_3} & {E_2} \\
{E_2} & E
\arrow["{\id_E\times_As_q}", from=1-1, to=1-2]
\arrow["{s_q\times_A\id_E}"', from=1-1, to=2-1]
\arrow["{s_q}", from=1-2, to=2-2]
\arrow["{s_q}"', from=2-1, to=2-2]
\end{tikzcd}
\end{equation*}

\item Unitality:
\begin{equation*}
% https://q.uiver.app/#q=WzAsMyxbMCwwLCJFXzIiXSxbMCwxLCJFIl0sWzEsMCwiRSJdLFsxLDAsIlxcPHF6X3EsXFxpZF9FXFw+Il0sWzAsMiwic19xIl0sWzEsMiwiIiwyLHsibGV2ZWwiOjIsInN0eWxlIjp7ImhlYWQiOnsibmFtZSI6Im5vbmUifX19XV0=
\begin{tikzcd}
{E_2} & E \\
E
\arrow["{s_q}", from=1-1, to=1-2]
\arrow["{\<qz_q,\id_E\>}", from=2-1, to=1-1]
\arrow[Rightarrow, no head, from=2-1, to=1-2]
\end{tikzcd}
\end{equation*}

\item Commutativity:
\begin{equation*}
% https://q.uiver.app/#q=WzAsNCxbMCwwLCJFXzIiXSxbMCwxLCJFXzIiXSxbMSwwLCJFIl0sWzEsMSwiRSJdLFsyLDMsIiIsMCx7ImxldmVsIjoyLCJzdHlsZSI6eyJoZWFkIjp7Im5hbWUiOiJub25lIn19fV0sWzAsMiwic19xIl0sWzEsMywic19xIiwyXSxbMCwxLCJcXHRhdSIsMl1d
\begin{tikzcd}
{E_2} & E \\
{E_2} & E
\arrow["{s_q}", from=1-1, to=1-2]
\arrow["\tau"', from=1-1, to=2-1]
\arrow[Rightarrow, no head, from=1-2, to=2-2]
\arrow["{s_q}"', from=2-1, to=2-2]
\end{tikzcd}
\end{equation*}
where $\tau\colon E_2\to E_2$ denotes the flip $\<\pi_2,\pi_1\>$.
\end{itemize}
A $1$-\textbf{morphism of additive bundle objects} $(\psi,\varphi)\colon(B,E,q,z_q,s_q)\to(B',E',q',z_q',s_q')$ over $\X$ consists of two $2$-morphisms $\varphi\colon E\Rightarrow E'$ and $\psi\colon B\Rightarrow B'$, satisfying the following properties:
\begin{itemize}
\item Compatibility with the projections:
\begin{equation*}
% https://q.uiver.app/#q=WzAsNCxbMCwwLCJFIl0sWzEsMCwiRSciXSxbMCwxLCJBIl0sWzEsMSwiQSciXSxbMCwxLCJnIl0sWzAsMiwicSIsMl0sWzEsMywicSciXSxbMiwzLCJmIiwyXV0=
\begin{tikzcd}
E & {E'} \\
B & {B'}
\arrow["g", from=1-1, to=1-2]
\arrow["q"', from=1-1, to=2-1]
\arrow["{q'}", from=1-2, to=2-2]
\arrow["f"', from=2-1, to=2-2]
\end{tikzcd}
\end{equation*}

\item Compatibility with the zero morphisms:
\begin{equation*}
% https://q.uiver.app/#q=WzAsNCxbMCwwLCJFIl0sWzEsMCwiRSciXSxbMCwxLCJBIl0sWzEsMSwiQSciXSxbMCwxLCJnIl0sWzIsMCwiel9xIl0sWzMsMSwiel9xJyIsMl0sWzIsMywiZiIsMl1d
\begin{tikzcd}
E & {E'} \\
B & {B'}
\arrow["g", from=1-1, to=1-2]
\arrow["{z_q}", from=2-1, to=1-1]
\arrow["f"', from=2-1, to=2-2]
\arrow["{z_q'}"', from=2-2, to=1-2]
\end{tikzcd}
\end{equation*}

\item Additivity:
\begin{equation*}
% https://q.uiver.app/#q=WzAsNCxbMCwxLCJFIl0sWzEsMSwiRSciXSxbMCwwLCJFXzIiXSxbMSwwLCJFJ18yIl0sWzAsMSwiZyIsMl0sWzIsMywiZ1xcdGltZXNfQWciXSxbMiwwLCJzX3EiLDJdLFszLDEsInNfcSciXV0=
\begin{tikzcd}
{E_2} & {E'_2} \\
E & {E'}
\arrow["{g\times_Bg}", from=1-1, to=1-2]
\arrow["{s_q}"', from=1-1, to=2-1]
\arrow["{s_q'}", from=1-2, to=2-2]
\arrow["g"', from=2-1, to=2-2]
\end{tikzcd}
\end{equation*}
\end{itemize}
\end{definition}

\begin{notation}
\label{notation:whiskering}
In the following, given two $1$-morphisms $\T\colon\X\to\X$, $\T'\colon\X'\to\X'$ and two $2$-morphisms
\begin{equation*}
\begin{tikzcd}
\X & {\X'} \\
{\X''} & {\X'}
\arrow[from=1-1, to=1-2]
\arrow[Rightarrow, no head, from=1-2, to=2-2]
\arrow[from=2-1, to=2-2]
\arrow[from=1-1, to=2-1]
\arrow["f"{description}, Rightarrow, from=1-2, to=2-1]
\end{tikzcd}\hfill
\begin{tikzcd}
\X & {\X'} \\
\X & {\X''}
\arrow[from=1-1, to=1-2]
\arrow[from=1-2, to=2-2]
\arrow[from=2-1, to=2-2]
\arrow[Rightarrow, no head, from=1-1, to=2-1]
\arrow["g"{description}, Rightarrow, from=1-2, to=2-1]
\end{tikzcd}
\end{equation*}
we write $\T'f$ for
\begin{equation*}
\begin{tikzcd}
\X & {\X'} & {\X'} \\
{\X''} & {\X'} & {\X'}
\arrow[from=1-1, to=1-2]
\arrow[from=1-1, to=2-1]
\arrow["{\T'}", from=1-2, to=1-3]
\arrow["f"{description}, Rightarrow, from=1-2, to=2-1]
\arrow[Rightarrow, no head, from=1-2, to=2-2]
\arrow[Rightarrow, no head, from=1-3, to=2-3]
\arrow[from=2-1, to=2-2]
\arrow["{\T'}"', from=2-2, to=2-3]
\end{tikzcd}
\end{equation*}
and $g_\T$ for
\begin{equation*}
% https://q.uiver.app/#q=WzAsNixbMSwwLCJcXFgiXSxbMiwwLCJcXFgnIl0sWzEsMSwiXFxYIl0sWzIsMSwiXFxYJyciXSxbMCwwLCJcXFgiXSxbMCwxLCJcXFgiXSxbMCwxXSxbMiwzXSxbMCwyLCIiLDIseyJsZXZlbCI6Miwic3R5bGUiOnsiaGVhZCI6eyJuYW1lIjoibm9uZSJ9fX1dLFsxLDNdLFsxLDIsImciLDEseyJsZXZlbCI6Mn1dLFs0LDAsIlxcVCJdLFs1LDIsIlxcVCIsMl0sWzQsNSwiIiwyLHsibGV2ZWwiOjIsInN0eWxlIjp7ImhlYWQiOnsibmFtZSI6Im5vbmUifX19XV0=
\begin{tikzcd}
\X & \X & {\X'} \\
\X & \X & {\X''}
\arrow["\T", from=1-1, to=1-2]
\arrow[Rightarrow, no head, from=1-1, to=2-1]
\arrow[from=1-2, to=1-3]
\arrow[Rightarrow, no head, from=1-2, to=2-2]
\arrow["g"{description}, Rightarrow, from=1-3, to=2-2]
\arrow[from=1-3, to=2-3]
\arrow["\T"', from=2-1, to=2-2]
\arrow[from=2-2, to=2-3]
\end{tikzcd}
\end{equation*}
\end{notation}

A tangent object $(\X,\TT)$ in $\CC$ is an object $\X$ of $\CC$ equipped with the following data:
\begin{description}
\item[tangent $1$-morphism] A $1$-morphism $\T\colon\X\to\X$;

\item[projection] A $2$-morphism $p\colon\T\Rightarrow\id_\X$;

\item[zero $2$-morphism] A $2$-morphism $z\colon\id_\X\Rightarrow\T$;

\item[sum $2$-morphism] A $2$-morphism $s\colon\T_2\Rightarrow\T$, where $\T_2$ denotes the pullback of $p$ along iself;
\item[vertical lift] A $2$-morphism $l\colon\T\Rightarrow\T^2$;

\item[canonical flip] A $2$-morphism $c\colon\T^2\Rightarrow\T^2$;
\end{description}

satisfying the following properties:
\begin{itemize}
\item Additive structure. $(\id_\X,\T,p,z,s)$ is an additive bundle object of $\CC$;

\item Additivity of vertical lift:
\begin{align*}
&(z,l)\colon(\id_\X,\T,p,z,s)\to(\T,\T^2,\T p,\T z,\T s)
\end{align*}
is a morphism of additive bundle objects;

\item Additivity of canonical flip:
\begin{align*}
&(\id_\T,c)\colon(\T,\T^2,\T p,\T z,\T s)\to(\T,\T^2,p_\T,z_\T,s_\T)
\end{align*}
is a morphism of additive bundle objects;

\item Coassociativity of vertical lift:
\begin{equation*}
% https://q.uiver.app/#q=WzAsNCxbMSwxLCJcXFReM0EiXSxbMCwwLCJcXFQgQSJdLFsxLDAsIlxcVF4yQSJdLFswLDEsIlxcVF4yQSJdLFsxLDMsImwiLDJdLFsxLDIsImwiXSxbMiwwLCJcXFQgbCJdLFszLDAsImxcXFQiLDJdXQ==
\begin{tikzcd}
{\T} & {\T^2} \\
{\T^2} & {\T^3}
\arrow["l", from=1-1, to=1-2]
\arrow["l"', from=1-1, to=2-1]
\arrow["{\T l}", from=1-2, to=2-2]
\arrow["l\T"', from=2-1, to=2-2]
\end{tikzcd}
\end{equation*}

\item The canonical flip is a symmetric braiding:
\begin{equation*}
% https://q.uiver.app/#q=WzAsMyxbMCwwLCJcXFReMkEiXSxbMSwwLCJcXFReMkEiXSxbMCwxLCJcXFReMkEiXSxbMiwwLCJjIl0sWzAsMSwiYyJdLFsyLDEsIiIsMix7ImxldmVsIjoyLCJzdHlsZSI6eyJoZWFkIjp7Im5hbWUiOiJub25lIn19fV1d
\begin{tikzcd}
{\T^2} & {\T^2} \\
{\T^2}
\arrow["c", from=1-1, to=1-2]
\arrow["c", from=2-1, to=1-1]
\arrow[Rightarrow, no head, from=2-1, to=1-2]
\end{tikzcd}\hfill
% https://q.uiver.app/#q=WzAsNixbMCwwLCJcXFReM0EiXSxbMSwwLCJcXFReM0EiXSxbMiwwLCJcXFReM0EiXSxbMiwxLCJcXFReM0EiXSxbMSwxLCJcXFReM0EiXSxbMCwxLCJcXFReM0EiXSxbMCwxLCJcXFQgYyJdLFs0LDMsImNcXFQiLDJdLFs1LDQsIlxcVCBjIiwyXSxbMSwyLCJjXFxUIl0sWzAsNSwiY1xcVCIsMl0sWzIsMywiXFxUIGMiXV0=
\begin{tikzcd}
{\T^3} & {\T^3} & {\T^3} \\
{\T^3} & {\T^3} & {\T^3}
\arrow["{\T c}", from=1-1, to=1-2]
\arrow["c\T"', from=1-1, to=2-1]
\arrow["c\T", from=1-2, to=1-3]
\arrow["{\T c}", from=1-3, to=2-3]
\arrow["{\T c}"', from=2-1, to=2-2]
\arrow["c\T"', from=2-2, to=2-3]
\end{tikzcd}
\end{equation*}

\item Compatibility between canonical flip and vertical lift:
\begin{equation*}
% https://q.uiver.app/#q=WzAsMyxbMCwwLCJcXFReMkEiXSxbMCwxLCJcXFQgQSJdLFsxLDAsIlxcVF4yQSJdLFsxLDAsImwiXSxbMCwyLCJjIl0sWzEsMiwibCIsMl1d
\begin{tikzcd}
{\T^2} & {\T^2} \\
{\T}
\arrow["c", from=1-1, to=1-2]
\arrow["l", from=2-1, to=1-1]
\arrow["l"', from=2-1, to=1-2]
\end{tikzcd}\hfill
% https://q.uiver.app/#q=WzAsNSxbMSwwLCJcXFReM0EiXSxbMiwwLCJcXFReM0EiXSxbMCwwLCJcXFReMkEiXSxbMiwxLCJcXFReM0EiXSxbMCwxLCJcXFReMkEiXSxbMCwxLCJcXFQgYyJdLFsyLDAsImxcXFQiXSxbMSwzLCJjXFxUIl0sWzIsNCwiYyIsMl0sWzQsMywiXFxUIGwiLDJdXQ==
\begin{tikzcd}
{\T^2} & {\T^3} & {\T^3} \\
{\T^2} && {\T^3}
\arrow["l\T", from=1-1, to=1-2]
\arrow["c"', from=1-1, to=2-1]
\arrow["{\T c}", from=1-2, to=1-3]
\arrow["c\T", from=1-3, to=2-3]
\arrow["{\T l}"', from=2-1, to=2-3]
\end{tikzcd}
\end{equation*}

\item Universality of the vertical lift. The following diagram
\begin{equation*}
% https://q.uiver.app/#q=WzAsNCxbMSwxLCJcXFQgQSJdLFswLDAsIlxcVF8yQSJdLFsxLDAsIlxcVF4yQSJdLFswLDEsIkEiXSxbMSwzLCJcXHBpXzFwIiwyXSxbMSwyLCJcXHBoaSJdLFsyLDAsIlxcVCBwIl0sWzMsMCwieiIsMl1d
\begin{tikzcd}
{\T_2} & {\T^2} \\
A & {\T}
\arrow["\phi", from=1-1, to=1-2]
\arrow["{\pi_1p}"', from=1-1, to=2-1]
\arrow["{\T p}", from=1-2, to=2-2]
\arrow["z"', from=2-1, to=2-2]
\end{tikzcd}
\end{equation*}
is a pointwise pullback diagram, where:
\begin{align*}
&\phi\colon\T_2A\xrightarrow{l\times_Az_\T}\T\T_2A\xrightarrow{\T s}\T^2A
\end{align*}
\end{itemize}
We refer to the tuple $\T\=(\T,p,z,s,l,c)$ as a \textbf{tangent structure} over $\X$. Finally, a tangent object with negatives is a tangent object equipped with an extra structure:
\begin{description}
\item[negation] A $2$-morphism $n\colon\T\Rightarrow\T$, called \textbf{negation}, such that:
\begin{equation*}
% https://q.uiver.app/#q=WzAsNCxbMSwwLCJcXFRfMkEiXSxbMCwwLCJcXFQgQSJdLFsxLDEsIlxcVCBBIl0sWzAsMSwiQSJdLFsxLDAsIlxcPFxcaWQsblxcPiJdLFswLDIsInMiXSxbMSwzLCJwIiwyXSxbMywyLCJ6IiwyXV0=
\begin{tikzcd}
{\T} & {\T_2} \\
\id_\X & {\T}
\arrow["{\<\id,n\>}", from=1-1, to=1-2]
\arrow["p"', from=1-1, to=2-1]
\arrow["s", from=1-2, to=2-2]
\arrow["z"', from=2-1, to=2-2]
\end{tikzcd}
\end{equation*}
\end{description}

We introduce the following naming convention.

\begin{notation}
\label{notation:name-tangent-objects}
Given a $2$-category $\CC$ whose objects are called with a name "\textit{x}", we refer to a tangent object of $\CC$ as a \textbf{tangent \textit{x}}.
\end{notation}

The next example shows that our naming convention is consistent with the notion of tangent category, that is, tangent categories are tangent objects in the $2$-category of categories.

\begin{example}
\label{example:tangent-cats-vs-tangent-objects}
The obvious example of tangent objects is given by tangent categories. Thanks to Leung's classification theorem, a tangent category is a category $\X$ equipped with a strict monoidal functor $\Leung[\TT]\colon\Weil\to\End(\X)$ satisfying some universality conditions. So, by taking the $2$-category $\Cat$ of (small) categories, functors and natural transformations, we see that a tangent object of $\Cat$ is precisely a tangent category.
\par Notice that, as pointed out by Lucyshyn-Wright (see Remark~\ref{remark:pointwise-limits}), for tangent objects of $\Cat$ to be tangent categories, it is important that the limit diagrams involved in Definition~\ref{definition:tangent-object} are pointwise.
\end{example}

\begin{example}
\label{example:tangent-monads-vs-tangent-objets}
Let $\CC$ be a $2$-category and consider the $2$-category $\Mnd(\CC)$ whose objects are pairs $(\X,S)$ formed by an object $\X$ of $\CC$ and a monad $S$ of $\X$. Recall that a monad in a $2$-category over an object $\X$ consists of a monoid in the monoidal category $\End(\X)$ of endomorphisms of $\X$. Concretely, a monad consists of an endomorphism $S\colon\X\to\X$ together with two $2$-morphisms $\eta\colon\id_\X\Rightarrow S$ and $\gamma\colon S^2\Rightarrow S$, where $S^2\= S\o S$, satisfying associativity and unitality conditions. A morphism of monads $(F,\alpha)\colon(\X,S)\to(\X',S')$ consists of a $1$-morphism $F\colon\X\to\X'$ together with a $2$-morphism:
\begin{equation*}
\begin{tikzcd}
\X & {\X'} \\
\X & {\X'}
\arrow["S"', from=1-1, to=2-1]
\arrow["{S'}", from=1-2, to=2-2]
\arrow["F", from=1-1, to=1-2]
\arrow["F"', from=2-1, to=2-2]
\arrow["\alpha"{description}, Rightarrow, from=1-2, to=2-1]
\end{tikzcd}
\end{equation*}
satisfying some compatibilities with the units $\eta$ and $\eta'$ and the multiplications $\gamma$ and $\gamma'$ of the monads $S$ and $S'$ (see~\cite{street:formal-theory-monads} for details). A $2$-morphism $\theta\colon(F,\alpha)\Rightarrow(G,\beta)$ of $\Mnd(\CC)$ between two $1$-morphisms $(F,\alpha),(G,\beta)\colon(\X,\TT)\to(\X',\TT')$ consists of a $2$-morphism $\theta\colon F\Rightarrow G$, compatible with the distributive laws.
\par By spelling out the details one finds out that a tangent object of $\Mnd(\CC)$ consists of a tangent object $(\X,\TT)$ of $\CC$ together with a monad $S$ of $\X$ equipped with a $2$-morphism $\alpha\colon S\o\T\Rightarrow\T\o S$ compatible with the tangent structure $\TT$ of $\X$. We refer to $(\X,\TT;S,\alpha)$ as a \textbf{tangent monad} in $\CC$. When $\CC$ is the $2$-category $\Cat$, tangent monads in $\CC$ are precisely tangent monads, as introduced by Cockett, Lemay, and Lucyshyn-Wright in~\cite{cockett:tangent-monads}.
\end{example}

\begin{example}
\label{example:tangent-monoidal-categories}
Let us consider the $2$-category $\MonCat$ whose objects are monoidal categories $(\X,\x,\1,\alpha,\lambda,\rho)$ with associator $\alpha$ and left and right unitors $\lambda$ and $\rho$, respectively, $1$-morphisms are strong monoidal functors $(F,\epsilon,\mu)\colon(\X,\x,\1,\alpha,\lambda,\rho)\to(\X',\x',\1',\alpha',\lambda',\rho')$, i.e., functors $F\colon\X\to\X'$ together with an isomorphism $\epsilon\colon\1'\to F\1$ and a natural isomorphism $\mu_{X,Y}\colon F(X)\x'F(Y)\to F(X\x Y)$, compatible with the associators and the unitors, and $2$-morphisms are natural transformations compatible with the morphisms $\epsilon$ and $\mu$ of the strong monoidal functors.
\par A tangent object of $\MonCat$ consists of a monoidal category $(\X,\x,\1,\alpha,\lambda,\rho)$ equipped with a tangent structure, so that $(\X,\TT)$ is also a tangent category and with an isomorphism $\1\to\T(\1)$, that we call \textbf{tangent unitor}, and a natural isomorphism $\T M\x\T N\to\T(M\x N)$ that we call \textbf{tangent distributor}, compatible with the associator and the unitors. Employing the Convention~\ref{notation:name-tangent-objects}, we call the tangent objects of $\MonCat$, \textbf{tangent monoidal categories}.
\par Notice that the $2$-category $\TngCat$ of tangent categories admits products. This allows one to wonder what are pseudomonoids of $\TngCat$. Recall that a pseudomonoid in a $2$-category $\CC$ with products consists of an object $\X$, together with two $1$-morphisms $\x\colon\X\times\X\to\X$ and $\eta\colon\1\to\X$, $\1$ being terminal in $\CC$, and three $2$-isomorphisms:
\begin{equation*}
% https://q.uiver.app/#q=WzAsNCxbMCwwLCJcXFhcXHRpbWVzXFxYXFx0aW1lc1xcWCJdLFsxLDAsIlxcWFxcdGltZXNcXFgiXSxbMSwxLCJcXFgiXSxbMCwxLCJcXFhcXHRpbWVzXFxYIl0sWzAsMSwiXFx4XFx0aW1lc1xcaWRfXFxYIl0sWzAsMywiXFxpZF9cXFhcXHRpbWVzXFx4IiwyXSxbMSwyLCJcXHgiXSxbMywyLCJcXHgiLDJdLFsxLDMsIlxcYWxwaGEiLDEseyJsZXZlbCI6Mn1dXQ==
\begin{tikzcd}
\X\times\X\times\X & \X\times\X \\
\X\times\X & \X
\arrow["{\x\times\id_\X}", from=1-1, to=1-2]
\arrow["{\id_\X\times\x}"', from=1-1, to=2-1]
\arrow["\x", from=1-2, to=2-2]
\arrow["\x"', from=2-1, to=2-2]
\arrow["\alpha"{description}, Rightarrow, from=1-2, to=2-1]
\end{tikzcd}\hfill
% https://q.uiver.app/#q=WzAsNSxbMSwwLCJcXDFcXHRpbWVzXFxYIl0sWzIsMCwiXFxYXFx0aW1lc1xcWCJdLFsyLDEsIlxcWCJdLFswLDAsIlxcWCJdLFswLDEsIlxcWCJdLFswLDEsIlxcZXRhXFx0aW1lc1xcaWRfXFxYIl0sWzEsMiwiXFx4Il0sWzMsMCwiXFw8LCEsXFxpZF9cXFhcXD4iXSxbNCwyLCIiLDIseyJsZXZlbCI6Miwic3R5bGUiOnsiaGVhZCI6eyJuYW1lIjoibm9uZSJ9fX1dLFszLDQsIiIsMix7ImxldmVsIjoyLCJzdHlsZSI6eyJoZWFkIjp7Im5hbWUiOiJub25lIn19fV0sWzEsNCwiXFxsYW1iZGEiLDEseyJsZXZlbCI6Mn1dXQ==
\begin{tikzcd}
\X & \1\times\X & \X\times\X \\
\X && \X
\arrow["{\eta\times\id_\X}", from=1-2, to=1-3]
\arrow["\x", from=1-3, to=2-3]
\arrow["{\<!,\id_\X\>}", from=1-1, to=1-2]
\arrow[Rightarrow, no head, from=2-1, to=2-3]
\arrow[Rightarrow, no head, from=1-1, to=2-1]
\arrow["\lambda"{description}, Rightarrow, from=1-3, to=2-1]
\end{tikzcd}
\end{equation*}
\begin{equation*}
% https://q.uiver.app/#q=WzAsNSxbMSwwLCJYXFx0aW1lc1xcMSJdLFsyLDAsIlxcWFxcdGltZXNcXFgiXSxbMiwxLCJcXFgiXSxbMCwwLCJcXFgiXSxbMCwxLCJcXFgiXSxbMCwxLCJcXGlkX1xcWFxcdGltZXNcXGV0YSJdLFsxLDIsIlxceCJdLFszLDAsIlxcPFxcaWRfXFxYLCFcXD4iXSxbNCwyLCIiLDIseyJsZXZlbCI6Miwic3R5bGUiOnsiaGVhZCI6eyJuYW1lIjoibm9uZSJ9fX1dLFszLDQsIiIsMix7ImxldmVsIjoyLCJzdHlsZSI6eyJoZWFkIjp7Im5hbWUiOiJub25lIn19fV0sWzEsNCwiXFxyaG8iLDEseyJsZXZlbCI6Mn1dXQ==
\begin{tikzcd}
\X & X\times\1 & \X\times\X \\
\X && \X
\arrow["{\id_\X\times\eta}", from=1-2, to=1-3]
\arrow["\x", from=1-3, to=2-3]
\arrow["{\<\id_\X,!\>}", from=1-1, to=1-2]
\arrow[Rightarrow, no head, from=2-1, to=2-3]
\arrow[Rightarrow, no head, from=1-1, to=2-1]
\arrow["\rho"{description}, Rightarrow, from=1-3, to=2-1]
\end{tikzcd}
\end{equation*}
satisfying the same pentagonal and hexagonal diagrams of the associator and unitors in the definition of a monoidal category.
\end{example}

\par As tangent objects in the $2$-category of monads are equivalent to monads in the $2$-category of tangent objects, tangent monoidal categories are equivalent to pseudomonoids in the $2$-category $\TngCat$ of tangent categories. We refer to this second description as \textbf{monoidal tangent categories} and in the future we use the terms tangent monoidal categories and monoidal tangent categories, interchangeably.

\begin{proposition}
\label{proposition:tangent-monoidal-categories}
There is an equivalence between the category of tangent monoidal categories and the category of monoidal tangent categories, defined as pseudomonoids in the category of tangent categories.
\end{proposition}

\par Example~\ref{example:tangent-monoidal-categories} can be extended to other classes of monoidal categories. For example, one can consider braided monoidal categories, symmetric monoidal categories, or closed symmetric monoidal categories. The corresponding tangent objects are then called \textbf{tangent braided monoidal categories}, \textbf{tangent symmetric monoidal categories}, and \textbf{tangent closed symmetric monoidal categories}.

\begin{example}
\label{example:tangent-enriched-categories}
An enriched category $\Y$ over a monoidal category $(\X,\x,\1,\alpha,\lambda,\rho)$ consists of a collection of objects together, for each pair $M,N$ of objects, an object $\Y(M,N)$ of $\X$, which plays the role of the Hom-Set functor (cf.~\cite{kelly:enriched-cats}). Moreover, an enriched category comes equipped with a collection of morphisms of $\X$:
\begin{align*}
&\o\colon\Y(N,P)\x\Y(M,N)\to\Y(M,P)
\end{align*}
which plays the role of the composition of morphisms, and a collection of morphisms of $\X$:
\begin{align*}
&\id\colon\1\to\Y(M,M)
\end{align*}
which plays the role of the identity morphisms. One can define a $2$-category $\Enrch$ whose objects are pairs $(\X,\Y)$ formed by a monoidal category $\X\=(\X,\x,\1,\alpha,\lambda,\rho)$ together with an enriched category $\Y$ over $\X$. A $1$-morphism $(F,G,\beta)\colon(\X,\Y)\to(\X',\Y')$ of $\Enrch$ consists of a strong monoidal functor $F\=(F,\epsilon,\eta)\colon\X\to\X'$ together with a pair $(G,\beta)\colon\Y\to\Y'$, formed by an assignment $G$ which sends an object $M$ of $\Y$ to an object $GM$ of $\Y'$, together with a collection of isomorphisms of $\X$:
\begin{equation*}
\beta\colon F(\Y(M,N))\to\Y'(GM,GN)
\end{equation*}
compatible with the morphisms $\o$ and $\id$. Finally, $2$-morphisms $(\varphi,\psi)\colon(F,G,\beta)\to(F',G',\beta')$ between two $1$-morphisms $(F,G,\beta),(F',G',\beta')\colon(\X,\Y)\to(\X',\Y')$ consist of a natural transformation of strong monoidal functors $\varphi\colon F\to F'$ together with a collection of morphisms:
\begin{align*}
&\psi\colon\Y'(GM,GN)\to\Y'(G'M,G'N)
\end{align*}
satisfying the following condition:
\begin{equation*}
% https://q.uiver.app/#q=WzAsNCxbMCwwLCJGXFxZKE0sTikiXSxbMSwwLCJGJ1xcWShNLE4pIl0sWzAsMSwiXFxZJyhHTSxHTikiXSxbMSwxLCJcXFknKEcnTSxHJ04pIl0sWzAsMSwiXFx2YXJwaGkiXSxbMiwzLCJcXHBzaSIsMl0sWzAsMiwiXFxiZXRhIiwyXSxbMSwzLCJcXGJldGEnIl1d
\begin{tikzcd}
{F\Y(M,N)} & {F'\Y(M,N)} \\
{\Y'(GM,GN)} & {\Y'(G'M,G'N)}
\arrow["\varphi", from=1-1, to=1-2]
\arrow["\beta"', from=1-1, to=2-1]
\arrow["{\beta'}", from=1-2, to=2-2]
\arrow["\psi"', from=2-1, to=2-2]
\end{tikzcd}
\end{equation*}
Spelling out the details one finds that a tangent object of $\Enrch$ consists of tangent monoidal category $(\X,\TT)$ together with an $\X$-enriched category $\Y$ equipped with an assignment $\T'$ which sends an object $M$ to another object $\T'M$ of $\Y$, with a collection of isomorphisms $\beta\colon\T(\Y(M,N))\to\Y(\T'M,\T'N)$, compatible with $\o$ and $\id$. Moreover, $\Y$ comes equipped with a list of collections of morphisms of $\X$, $p'\colon\Y(\T' M,\T'N)\to\Y(M,N)$, $z'\colon\Y(M,N)\to(\T'M,\T'N)$, $s'\colon\Y_2(\T'M,\T'N)\to\Y(\T'M,\T'N)$, $\Y_2(\T'M,\T'N)$ denoting the pullback of $p'$ along itself, $l'\colon\Y(\T'M,\T'N)\to\Y({\T'}^2M,{\T'}^2N)$, and $c\colon\Y({\T'}^2M,{\T'}^2N)\to\Y({\T'}^2M,{\T'}^2N)$, satisfying the compatibility conditions required by the axioms of a tangent structure of $\X$.
\end{example}

In future work, we intend to explore notions like tangent model categories, tangent internal categories, double tangent categories (that are tangent objects in the $2$-category of tangent categories), tangent double categories (that are tangent objects in the $2$-category of double categories), tangent topoi, tangent sheaves and many others.
\par The next step is to introduce $1$-morphisms of tangent objects.

\begin{definition}
\label{definition:morphism-of-tangent-objects}
A \textbf{lax $1$-morphism of tangent objects} $(F,\alpha)\colon(\X,\TT)\to(\X',\TT')$ between two tangent objects $(\X,\TT)$ and $(\X',\TT')$ in a $2$-category $\CC$ consists of a $1$-morphism $F\colon\X\to\X'$ of $\CC$ together with a $2$-morphism:
\begin{equation*}
\begin{tikzcd}
\X & \X \\
\X' & \X'
\arrow["\T", from=1-1, to=1-2]
\arrow["{\T'}"', from=2-1, to=2-2]
\arrow["F", from=1-2, to=2-2]
\arrow["F"', from=1-1, to=2-1]
\arrow["\alpha"{description}, Rightarrow, from=1-2, to=2-1]
\end{tikzcd}
\end{equation*}
compatible with the tangent structures as follows:
\begin{itemize}
\item Additivity: $(\id_F,\alpha)\colon(F,F\o\T,Fp,Fz,Fs)\to(F,\T'\o F,p'_F,z'_F,s'_F)$ is an additive bundle morphism;

\item Compatibility with the vertical lifts:
\begin{equation*}
% https://q.uiver.app/#q=WzAsNSxbMCwxLCJGXFxUIEEiXSxbMiwxLCJcXFQnRkEiXSxbMCwwLCJGXFxUXjJBIl0sWzEsMCwiXFxUJ0ZcXFQgQSJdLFsyLDAsIntcXFQnfV4yRkEiXSxbMCwxLCJcXGFscGhhIiwyXSxbMiwzLCJcXGFscGhhXFxUIl0sWzMsNCwiXFxUJ1xcYWxwaGEiXSxbMCwyLCJGbCJdLFsxLDQsImwnRiIsMl1d
\begin{tikzcd}
{F\o\T^2} & {\T'\o F\o\T} & {{\T'}^2\o F} \\
{F\o\T} && {\T'\o F}
\arrow["\alpha\T", from=1-1, to=1-2]
\arrow["{\T'\alpha}", from=1-2, to=1-3]
\arrow["Fl", from=2-1, to=1-1]
\arrow["\alpha"', from=2-1, to=2-3]
\arrow["{l'F}"', from=2-3, to=1-3]
\end{tikzcd}
\end{equation*}

\item Compatibility with the canonical flips:
\begin{equation*}
% https://q.uiver.app/#q=WzAsNixbMCwxLCJGXFxUXjJBIl0sWzIsMSwie1xcVCd9XjJGQSJdLFswLDAsIkZcXFReMkEiXSxbMSwwLCJcXFQnRlxcVCBBIl0sWzIsMCwie1xcVCd9XjJGQSJdLFsxLDEsIlxcVCdGXFxUIEEiXSxbMiwzLCJcXGFscGhhXFxUIl0sWzMsNCwiXFxUJ1xcYWxwaGEiXSxbMiwwLCJGYyIsMl0sWzQsMSwiYydGIl0sWzAsNSwiXFxhbHBoYVxcVCIsMl0sWzUsMSwiXFxUJ1xcYWxwaGEiLDJdXQ==
\begin{tikzcd}
{F\o\T^2} & {\T'\o F\o\T} & {{\T'}^2\o F} \\
{F\o\T^2} & {\T'\o F\o\T} & {{\T'}^2\o F}
\arrow["\alpha\T", from=1-1, to=1-2]
\arrow["Fc"', from=1-1, to=2-1]
\arrow["{\T'\alpha}", from=1-2, to=1-3]
\arrow["{c'F}", from=1-3, to=2-3]
\arrow["\alpha\T"', from=2-1, to=2-2]
\arrow["{\T'\alpha}"', from=2-2, to=2-3]
\end{tikzcd}
\end{equation*}
\end{itemize}
Similarly, a \textbf{colax $1$-morphism of tangent objects} $(G,\beta)\colon(\X,\TT)\nto(\X',\TT')$ consists of a $1$-morphism $G\colon\X\to\X'$ together with a $2$-morphism:
\begin{equation*}
\begin{tikzcd}
\X & {\X'} \\
\X & {\X'}
\arrow["{\T'}", from=1-2, to=2-2]
\arrow["\T"', from=1-1, to=2-1]
\arrow["F"', from=2-1, to=2-2]
\arrow["F", from=1-1, to=1-2]
\arrow["\beta"{description}, Rightarrow, from=1-2, to=2-1]
\end{tikzcd}
\end{equation*}
satisfying the dual conditions of $\alpha$. Moreover, the functor:
\begin{align*}
&\CC(\X,G)\colon\CC(\X,\X)\to\CC(\X,\X')
\end{align*}
is required to preserve the universality of the $n$-fold pointwise pullback of the projection along itself and the pointwise pullback of the universal property of the vertical lift.
\par A colax morphism $(F,\alpha)$ of tangent objects is \textbf{strong} if $\alpha$ is invertible and \textbf{strict} if $\alpha$ is the identity.
\end{definition}

We are also interested in defining $2$-morphisms of tangent objects.

\begin{definition}
\label{definition:2-morphisms-of-tangent-objects}
A \textbf{lax $2$-morphism of tangent objects} $\theta\colon(F,\alpha)\to(G,\beta)$ between two $1$-morphisms $(F,\alpha),(G,\beta)\colon(\X,\TT)\to(\X',\TT')$ consists of a $2$-morphism of $\CC$ $\theta\colon F\Rightarrow G$, satisfying the following condition:
\begin{equation*}
% https://q.uiver.app/#q=WzAsNCxbMCwwLCJGXFxvXFxUIl0sWzAsMSwiR1xcb1xcVCJdLFsxLDAsIlxcVCdcXG8gRiJdLFsxLDEsIlxcVCdcXG8gRyJdLFswLDIsIlxcYWxwaGEiXSxbMSwzLCJcXGJldGEiLDJdLFswLDEsIlxcdGhldGFcXFQiLDJdLFsyLDMsIlxcVCdcXHRoZXRhIl1d
\begin{tikzcd}
{F\o\T} & {\T'\o F} \\
{G\o\T} & {\T'\o G}
\arrow["\alpha", from=1-1, to=1-2]
\arrow["{\theta\T}"', from=1-1, to=2-1]
\arrow["{\T'\theta}", from=1-2, to=2-2]
\arrow["\beta"', from=2-1, to=2-2]
\end{tikzcd}
\end{equation*}
Similarly, a \textbf{colax $2$-morphism of tangent objects} $\theta\colon(F,\alpha)\to(G,\beta)$ between two $1$-morphisms $(F,\alpha),\linebreak(G,\beta)\colon(\X,\TT)\nto(\X',\TT')$ consists of a $2$-morphism of $\CC$ $\theta\colon F\Rightarrow G$, satisfying the dual conditions of a lax $2$-morphism. Finally, a \textbf{double morphism of tangent objects}
\begin{equation*}
% https://q.uiver.app/#q=WzAsNCxbMCwwLCIoXFxYX1xccyxcXFRUX1xccykiXSxbMSwwLCIoXFxYX1xccycsXFxUVF9cXHMnKSJdLFswLDEsIihcXFhfXFxiLFxcVFRfXFxiKSJdLFsxLDEsIihcXFhfXFxiJyxcXFRUX1xcYicpIl0sWzAsMiwiKEcsXFxiZXRhKSIsMix7InN0eWxlIjp7ImJvZHkiOnsibmFtZSI6ImJhcnJlZCJ9fX1dLFsxLDMsIihHJyxcXGJldGEnKSIsMCx7InN0eWxlIjp7ImJvZHkiOnsibmFtZSI6ImJhcnJlZCJ9fX1dLFswLDEsIihGX1xccyxcXGFscGhhX1xccykiXSxbMiwzLCIoRl9cXGIsXFxhbHBoYV9cXGIpIiwyXSxbMSwyLCJcXHRoZXRhIiwxLHsibGV2ZWwiOjJ9XV0=
\begin{tikzcd}
{(\X_\s,\TT_\s)} & {(\X_\s',\TT_\s')} \\
{(\X_\b,\TT_\b)} & {(\X_\b',\TT_\b')}
\arrow["{(F_\s,\alpha_\s)}", from=1-1, to=1-2]
\arrow["{(G,\beta)}"', "\shortmid"{marking}, from=1-1, to=2-1]
\arrow["\theta"{description}, Rightarrow, from=1-2, to=2-1]
\arrow["{(G',\beta')}", "\shortmid"{marking}, from=1-2, to=2-2]
\arrow["{(F_\b,\alpha_\b)}"', from=2-1, to=2-2]
\end{tikzcd}
\end{equation*}
for the lax $1$-morphisms $(F_\s,\alpha_\s)$ and $(F_\b,\alpha_\b)$, and the colax $1$-morphisms $(G,\beta)$ and $(G',\beta')$, is a $2$-morphism
\begin{equation*}
\begin{tikzcd}
{\X_\s} & {\X_\s'} \\
{\X_\b} & {\X_\b'}
\arrow["{F_\b}"', from=2-1, to=2-2]
\arrow["{F_\s}", from=1-1, to=1-2]
\arrow["{G'}", from=1-2, to=2-2]
\arrow["G"', from=1-1, to=2-1]
\arrow["\theta"{description}, Rightarrow, from=1-2, to=2-1]
\end{tikzcd}
\end{equation*}
satisfying the following condition:
\begin{equation*}
% https://q.uiver.app/#q=WzAsNixbMCwwLCJHJ1xcbyBGX1xcb1xcb1xcVF9cXG8iXSxbMSwwLCJHJ1xcb1xcVF9cXG8nXFxvIEZfXFxvIl0sWzIsMCwiXFxUX1xcYidcXG8gRydcXG8gRl9cXG8iXSxbMCwxLCJGX1xcYlxcbyBHXFxvXFxUX1xcbyJdLFsxLDEsIkZfXFxiXFxvXFxUX1xcbydcXG8gRyJdLFsyLDEsIlxcVF9cXGInXFxvIEZfXFxiXFxvIEciXSxbMCwxLCJHJ1xcYWxwaGFfXFxvIl0sWzEsMiwiXFxiZXRhJ0ZfXFxvIl0sWzAsMywiXFx0aGV0YVxcVF9cXG8iLDJdLFszLDQsIkZfXFxiXFxiZXRhIiwyXSxbNCw1LCJcXGFscGhhX1xcYiBHIiwyXSxbMiw1LCJcXFRfXFxiJ1xcdGhldGEiXV0=
\begin{tikzcd}
{G'\o F_\s\o\T_\s} & {G'\o\T_\s'\o F_\s} & {\T_\b'\o G'\o F_\s} \\
{F_\b\o G\o\T_\s} & {F_\b\o\T_\s'\o G} & {\T_\b'\o F_\b\o G}
\arrow["{G'\alpha_\s}", from=1-1, to=1-2]
\arrow["{\theta\T_\s}"', from=1-1, to=2-1]
\arrow["{\beta'F_\s}", from=1-2, to=1-3]
\arrow["{\T_\b'\theta}", from=1-3, to=2-3]
\arrow["{F_\b\beta}"', from=2-1, to=2-2]
\arrow["{\alpha_\b G}"', from=2-2, to=2-3]
\end{tikzcd}
\end{equation*}
\end{definition}

\begin{notation}
\label{notation:category-TNG}
Tangent objects of a $2$-category $\CC$ together with lax tangent $1$-morphisms and lax $2$-morphisms form a $2$-category $\Tng(\CC)$. Similarly, tangent objects of $\CC$ together with colax tangent $1$-morphisms and colax $2$-morphisms form a $2$-category $\Tng_\co(\CC)$. The $2$-subcategory of $\Tng(\CC)$ whose $1$-morphisms are strong, that is the distributive law is invertible, is denoted by $\Tng_\cong(\CC)$ and the $2$-subcategory of $\Tng(\CC)$ whose $1$-morphisms are strict, i.e., the distributive law is the identity, is denoted by $\Tng_=(\CC)$. Finally, tangent objects together with lax tangent $1$-morphisms as horizontal morphisms, colax tangent $1$-morphisms as vertical morphisms, and double tangent cells for double cells form also a double category denoted by $\TTng(\CC)$.
\end{notation}

When $\CC$ is the $2$-category $\Cat$ of categories, the double category $\TTng(\Cat)$ is precisely the double category $\TTngCat$ of tangent categories, first described in~\cite[Proposition~2.2]{lanfranchi:differential-bundles-operadic-affine-schemes}.
\par We conclude this section, by showing that the assignment $\Tng$ which sends a $2$-category $\CC$ to the $2$-category $\Tng(\CC)$ of tangent objects of $\CC$ extends to a $2$-functor. For this purpose, we first need to select the correct class of morphisms between $2$-categories. Indeed, if $F\colon\CC\to\CC'$ is an arbitrary $2$-functor and $(\X,\TT)$ a tangent object of $\CC$, in general, there is no reason to believe that $F\X$ is also a tangent object of $\CC'$.
\par Indeed, in order to make $F\X$ into a tangent object, the $2$-functor $F$ must preserve the $n$-fold pullbacks of the projection with itself and the universality of the vertical lift.
\par Recall that a $2$-functor $F\colon\CC\to\CC'$ (notice that here we work with strict $2$-functors) is an assignment which sends objects $M$, $1$-morphisms $f\colon M\to N$, and $2$-morphisms $\theta\colon f\Rightarrow g$ of $\CC$ to objects $F_0M$, $1$-morphisms $F_1f\colon F_0M\to F_0N$, and $2$-morphisms $F_2\theta\colon F_1f\Rightarrow F_2g$, respectively, in a compatible way with the composition and the identities.

\begin{definition}
\label{definition:tangentially-continuous-functors}
A $2$-functor $F\colon\CC\to\CC'$ is \textbf{$2$-pullback preserving} if it preserves pullbacks of type
\begin{equation*}
% https://q.uiver.app/#q=WzAsNCxbMCwwLCJcXGJ1bGxldCJdLFsxLDAsIlxcYnVsbGV0Il0sWzAsMSwiXFxidWxsZXQiXSxbMSwxLCJcXGJ1bGxldCJdLFswLDIsIlxccGlfMSIsMl0sWzAsMSwiXFxwaV8yIl0sWzEsMywicCJdLFsyLDMsInEiLDJdXQ==
\begin{tikzcd}
	\bullet & \bullet \\
	\bullet & \bullet
	\arrow["{\pi_2}", from=1-1, to=1-2]
	\arrow["{\pi_1}"', from=1-1, to=2-1]
	\arrow["p", from=1-2, to=2-2]
	\arrow["q"', from=2-1, to=2-2]
\end{tikzcd}
\end{equation*}
in each category $\End(\X)$, for each object $\X$ of $\CC$.
\end{definition}

If $F$ is a $2$-pullback preserving $2$-functor and $(\X,\TT)$ a tangent object of $\CC$, it is not hard to see that, $F\X\=F_0\X$ comes equipped with a tangent structure so defined:
\begin{description}
\item[tangent bundle morphism] The tangent bundle morphism $F\T\colon F\X\to F\X$ is given by:
\begin{align*}
&F_1\T\colon F_0\X\to F_0\X
\end{align*}
\item[projection] The projection $Fp\colon F\T\Rightarrow\id_{F\X}$ is given by:
\begin{align*}
&F_2p\colon F_1\T\Rightarrow F_1\id_\X=\id_{F_0\X}
\end{align*}
\item[zero morphism] The zero morphism $Fz\colon\id_{F\X}\Rightarrow F\T$ is given by:
\begin{align*}
&F_2z\colon\id_{F_0\X}=F_1\id_\X\Rightarrow F_1\T
\end{align*}
\item[sum morphism] The sum morphism $Fs\colon(F\T)_2\Rightarrow F\T$ is given by:
\begin{align*}
&(F_1\T)_2\cong F_1(\T_2)\xRightarrow{F_2s}F_1\T
\end{align*}
\item[vertical lift] The vertical lift $Fl\colon F\T\Rightarrow(F\T)^2$ is given by:
\begin{align*}
&F_2l\colon F_1\T\Rightarrow F_1(\T^2)=(F_1\T)^2
\end{align*}
\item[canonical flip] The canonical flip $Fc\colon(F\T)^2\Rightarrow(F\T)^2$ is given by:
\begin{align*}
&F_2c\colon(F_1\T)^2=F_1(\T^2)\Rightarrow F_1(\T^2)=(F_1\T)^2
\end{align*}
\end{description}
Moreover, if $(\X,\TT)$ has negatives with negation $n$, then:
\begin{description}
\item[negation] The negation $Fn\colon F\T\Rightarrow F\T$ is given by:
\begin{align*}
&F_2n\colon F_1\T\Rightarrow F_1\T
\end{align*}
\end{description}

\begin{remark}
\label{remark:tangential-continuity}
Note that the pullback preservation property is only a sufficient condition for a $2$-functor to preserve tangent objects. Indeed, one can ask for stricter conditions on $2$-functors. For the sake of simplicity, we decided to adopt the weaker condition expressed by Definition~\ref{definition:tangentially-continuous-functors}.
\end{remark}

\begin{notation}
\label{notation:category-TWOCAT}
$2$-categories, $2$-functors, and $2$-natural transformations form a $2$-category denoted by $\twoCat$. Moreover, the composition of two $2$-pullback preserving $2$-functors is still a $2$-pullback preserving $2$-functor. Thus, also $2$-categories, $2$-pullback preserving $2$-functors and $2$-natural transformations form a $2$-category that will be denoted by $\twoCat_\pp$.
\end{notation}

\begin{proposition}
\label{proposition:Tng-functiorality}
The assignment $\Tng$ which sends a $2$-category $\CC$ to the $2$-category $\Tng(\CC)$ of tangent objects of $\CC$ extends to a $2$-functor $\Tng\colon\twoCat_\pp\to\twoCat$. Similarly, also $\Tng_\co$, $\Tng_\cong$, and $\Tng_=$ extend to $2$-functors $\twoCat_\pp\to\twoCat$.
\end{proposition}

%__________________________________________________________________________
%__________________________________________________________________________

\section{The Grothendieck construction in the context of tangent categories}
\label{section:full-Grothendieck}
In Section~\ref{subsection:reduced-Grothendieck} we presented an adjunction between tangent fibrations and indexed tangent categories. However, this adjunction is not an equivalence of categories.
\par To find the correct notion, we employ the concept of a tangent object, introduced in Section~\ref{section:tangent-objects}. First, we show that tangent objects in the $2$-category $\Fib$ of fibrations are precisely tangent fibrations.

\begin{proposition}
\label{proposition:tangent-fibrations-are-tangent-objects}
The following are equivalences of $2$-categories:
\begin{align*}
\Tng(\Fib)&\simeq\TngFib\\
\Tng_\co(\Fib)&\simeq\TngFib_\co\\
\Tng_\cong(\Fib)&\simeq\TngFib_\cong\\
\Tng_=(\Fib)&\simeq\TngFib_=
\end{align*}
\end{proposition}

\begin{remark}
\label{remark:fixed-base-morphisms}
Notice that we do not require the base tangent morphism of a $1$-morphim in $\TngFib$ to be strict anymore (see Remark~\ref{remark:strict-on-base-explanation}), as required in Definition~\ref{definition:morphism-of-tangent-fibrations}. This was employed to simplify the definition of the functor $\I$ in the reduced Grothendieck construction. Now, there is no more need for this technical assumption. Notice also that, if one prefers base tangent morphisms to be strict, one can simply consider $2$-morphisms $(\theta;\theta')$ in $\Fib$ where the base natural transformation $\theta$ is the identity. Then, $1$-morphisms of tangent objects in this $2$-category are precisely pairs of tangent morphisms whose base morphism is strict.
\end{remark}

The goal is to compare tangent fibrations with some notion of indexed categories. Thanks to Proposition~\ref{proposition:tangent-fibrations-are-tangent-objects}, we now know that tangent fibrations are tangent objects in $\Fib$. Moreover, the Grothendieck construction defines an equivalence between fibrations and indexed categories. Therefore, we can apply the $2$-functor $\Tng$ to this equivalence.

\begin{definition}
\label{definition:tangent-indexed-category}
A \textbf{tangent indexed category} is a tangent object in the $2$-category $\Indx$ of indexed categories. Moreover, a tangent indexed category \textbf{with negatives} is a tangent object with negatives in the $2$-category $\Indx$.
\end{definition}

\begin{remark}
\label{remark:tangent-indexed-vs-indexed-tangent}
A tangent indexed category should not be confused with an indexed tangent category, defined in Definition~\ref{definition:indexed-tangent-category}. We show in a moment that these two notions are related via an adjunction, which however is not an equivalence.
\end{remark}

Let us unpack Definition~\ref{definition:tangent-indexed-category}. A tangent indexed category consists of the following data:
\begin{description}
\item[base tangent category] A base tangent category $(\X,\TT)$;

\item[indexed category] An indexed category $\I\colon\X^\op\to\Cat$;

\item[indexed tangent bundle functor] A collection of functors $\T\^A\colon\X\^A\to\X\^{\T A}$, indexed by the objects $A$ of $\X$;

\item[tangent distributors] A collection of natural transformations $\kappa\^f$, indexed by the morphisms $f\colon A\to A'$ of $\X$:
\begin{equation*}
\begin{tikzcd}
{\X\^{A'}} & {\X\^A} \\
{\X\^{\T A'}} & {\X\^{\T A}}
\arrow["{\T\^{A'}}"', from=1-1, to=2-1]
\arrow["{\T\^A}", from=1-2, to=2-2]
\arrow["{f^\*}", from=1-1, to=1-2]
\arrow["{(\T f)^\*}"', from=2-1, to=2-2]
\arrow["{\kappa\^f}"{description}, Rightarrow, from=1-2, to=2-1]
\end{tikzcd}
\end{equation*}

\item[indexed projection] A collection of natural transformations $p\^A\colon\T\^A\Rightarrow p^\*$, indexed by the objects $A$ of $\X$
\begin{equation*}
\begin{tikzcd}
{\X\^A} & {\X\^{\T A}} \\
{\X\^A} & {\X\^{\T A}}
\arrow["{\T\^A}", from=1-1, to=1-2]
\arrow["{p^\*}"', from=2-1, to=2-2]
\arrow[Rightarrow, no head, from=1-1, to=2-1]
\arrow[Rightarrow, no head, from=1-2, to=2-2]
\arrow["{p\^A}"{description}, Rightarrow, from=1-2, to=2-1]
\end{tikzcd}
\end{equation*}
satisfying the following property:
\begin{equation*}
% https://q.uiver.app/#q=WzAsNCxbMCwwLCJcXFRcXF5BXFxvIGZeXFwqIl0sWzEsMCwicF5cXCpcXG8gZl5cXCoiXSxbMCwxLCIoXFxUIGYpXlxcKlxcb1xcVFxcXntBJ30iXSxbMSwxLCIoXFxUIGYpXlxcKlxcbyBwXlxcKiJdLFswLDEsInBcXF5BIGZeXFwqIl0sWzAsMiwiXFxrYXBwYVxcXmYiLDJdLFsyLDMsIihcXFQgZileXFwqcFxcXntBJ30iLDJdLFsxLDMsIlxcZ2FtbWEiXV0=
\begin{tikzcd}
{\T\^A\o f^\*} & {p^\*\o f^\*} \\
{(\T f)^\*\o\T\^{A'}} & {(\T f)^\*\o p^\*}
\arrow["{p\^A f^\*}", from=1-1, to=1-2]
\arrow["{\kappa\^f}"', from=1-1, to=2-1]
\arrow["\Comp", from=1-2, to=2-2]
\arrow["{(\T f)^\*p\^{A'}}"', from=2-1, to=2-2]
\end{tikzcd}
\end{equation*}

\item[indexed zero] A collection of natural transformations $z\^A\colon\id_{\X\^A}\Rightarrow z^\*\o\T\^A$, indexed by the objects $A$ of $\X$
\begin{equation*}
% https://q.uiver.app/#q=WzAsNCxbMCwwLCJcXFhcXF5BIl0sWzEsMCwiXFxYXFxeQSJdLFswLDEsIlxcWFxcXntcXFQgQX0iXSxbMSwxLCJcXFhcXF5BIl0sWzAsMSwiIiwxLHsibGV2ZWwiOjIsInN0eWxlIjp7ImhlYWQiOnsibmFtZSI6Im5vbmUifX19XSxbMCwyLCJcXFRcXF5BIiwyXSxbMiwzLCJ6XlxcKiIsMl0sWzEsMywiIiwyLHsibGV2ZWwiOjIsInN0eWxlIjp7ImhlYWQiOnsibmFtZSI6Im5vbmUifX19XSxbMSwyLCJ6XFxeQSIsMSx7ImxldmVsIjoyfV1d
\begin{tikzcd}
{\X\^A} & {\X\^A} \\
{\X\^{\T A}} & {\X\^A}
\arrow[Rightarrow, no head, from=1-1, to=1-2]
\arrow["{\T\^A}"', from=1-1, to=2-1]
\arrow["{z^\*}"', from=2-1, to=2-2]
\arrow[Rightarrow, no head, from=1-2, to=2-2]
\arrow["{z\^A}"{description}, Rightarrow, from=1-2, to=2-1]
\end{tikzcd}
\end{equation*}
satisfying the following property:
\begin{equation*}
% https://q.uiver.app/#q=WzAsNCxbMCwwLCJmXlxcKiJdLFsxLDAsInpeXFwqXFxvXFxUXFxeQVxcbyBmXlxcKiJdLFsxLDEsInpeXFwqXFxvKFxcVCBmKV5cXCpcXG9cXFRcXF57QSd9Il0sWzAsMSwiZl5cXCpcXG8gel5cXCpcXG9cXFRcXF57QSd9Il0sWzAsMSwielxcXkEgZl5cXCoiXSxbMSwyLCJ6XlxcKlxca2FwcGFcXF5mIl0sWzAsMywiZl5cXCp6XFxee0EnfSIsMl0sWzIsMywiXFxnYW1tYSJdXQ==
\begin{tikzcd}
{f^\*} & {z^\*\o\T\^A\o f^\*} \\
{f^\*\o z^\*\o\T\^{A'}} & {z^\*\o(\T f)^\*\o\T\^{A'}}
\arrow["{z\^A f^\*}", from=1-1, to=1-2]
\arrow["{f^\*z\^{A'}}"', from=1-1, to=2-1]
\arrow["{z^\*\kappa\^f}", from=1-2, to=2-2]
\arrow["\Comp", from=2-2, to=2-1]
\end{tikzcd}
\end{equation*}

\item[indexed $n$-fold pullback] For every positive integer $n$, a collection of functors $\T\^A_n\colon\X\^A\to\X\^{\T_nA}$, indexed by the objects $A$ of $\X$, together with a collection of natural transformations $\kappa_n\^f\colon\T\^A_n\o f^\*\Rightarrow(\T_nf)^\*\o\T_n\^{A'}$, indexed by the morphism $f\colon A\to A'$ of $\X$
\begin{equation*}
% https://q.uiver.app/#q=WzAsNCxbMCwwLCJcXFhcXF5CIl0sWzEsMCwiXFxYXFxeQSJdLFsxLDEsIlxcWFxcXntcXFRfbkF9Il0sWzAsMSwiXFxYXFxee1xcVF9uQn0iXSxbMCwxLCJmXlxcKiJdLFsxLDIsIlxcVFxcXkFfbiJdLFswLDMsIlxcVFxcXkJfbiIsMl0sWzMsMiwiKFxcVF9uZileXFwqIiwyXSxbMSwzLCJcXHhpX25cXF5mIiwxLHsibGV2ZWwiOjJ9XV0=
\begin{tikzcd}
{\X\^{A'}} & {\X\^A} \\
{\X\^{\T_nA'}} & {\X\^{\T_nA}}
\arrow["{f^\*}", from=1-1, to=1-2]
\arrow["{\T\^A_n}", from=1-2, to=2-2]
\arrow["{\T\^{A'}_n}"', from=1-1, to=2-1]
\arrow["{(\T_nf)^\*}"', from=2-1, to=2-2]
\arrow["{\kappa_n\^f}"{description}, Rightarrow, from=1-2, to=2-1]
\end{tikzcd}
\end{equation*}
and a collection of natural tranformations $\pi\^A_k\colon\T\^A_n\Rightarrow\pi_k\^*\o\T\^A$, indexed by the objects $A$ of $\X$
\begin{equation*}
\begin{tikzcd}
{\X\^A} & {\X\^{\T_nA}} \\
{\X\^{\T A}} & {\X\^{\T_nA}}
\arrow["{\T\^A}"', from=1-1, to=2-1]
\arrow["{\pi_k^\*}"', from=2-1, to=2-2]
\arrow[Rightarrow, no head, from=1-2, to=2-2]
\arrow["{\T\^A_n}", from=1-1, to=1-2]
\arrow["{\pi\^A_k}"{description}, Rightarrow, from=1-2, to=2-1]
\end{tikzcd}
\end{equation*}
satisfying the following property:
\begin{equation*}
% https://q.uiver.app/#q=WzAsNSxbMCwwLCJcXFRfblxcXkFcXG8gZl5cXCoiXSxbMSwwLCJcXHBpX2teXFwqXFxvXFxUXFxeQVxcbyBmXlxcKiJdLFsxLDEsIlxccGlfa15cXCpcXG8oXFxUIGYpXlxcKlxcb1xcVFxcXntBJ30iXSxbMSwyLCIoXFxUX25mKV5cXCpcXG9cXHBpX2teXFwqXFxvXFxUXFxee0EnfSJdLFswLDIsIihcXFRfbmYpXlxcKlxcb1xcVF9uXFxee0EnfSJdLFswLDEsIlxccGlfa1xcXkFmXlxcKiJdLFsxLDIsIlxccGlfa15cXCpcXGthcHBhXFxeZiJdLFsyLDMsIlxcZ2FtbWEiXSxbMCw0LCJcXGthcHBhX25cXF5mZl5cXCoiLDJdLFs0LDMsIihcXFRfbmYpXlxcKlxccGlfa1xcXntBJ30iLDJdXQ==
\begin{tikzcd}
{\T_n\^A\o f^\*} & {\pi_k^\*\o\T\^A\o f^\*} \\
& {\pi_k^\*\o(\T f)^\*\o\T\^{A'}} \\
{(\T_nf)^\*\o\T_n\^{A'}} & {(\T_nf)^\*\o\pi_k^\*\o\T\^{A'}}
\arrow["{\pi_k\^Af^\*}", from=1-1, to=1-2]
\arrow["{\kappa_n\^ff^\*}"', from=1-1, to=3-1]
\arrow["{\pi_k^\*\kappa\^f}", from=1-2, to=2-2]
\arrow["\Comp", from=2-2, to=3-2]
\arrow["{(\T_nf)^\*\pi_k\^{A'}}"', from=3-1, to=3-2]
\end{tikzcd}
\end{equation*}

\item[indexed sum morphism] A collection of natural transformations $s\^A\colon\T_2\^A\Rightarrow s^\*\o\T\^A$, indexed by the objects $A$ of $\X$
\begin{equation*}
% https://q.uiver.app/#q=WzAsNCxbMCwwLCJcXFhcXF5BIl0sWzEsMCwiXFxYXFxee1xcVF8yQX0iXSxbMSwxLCJcXFhcXF57XFxUXzJBfSJdLFswLDEsIlxcWFxcXntcXFQgQX0iXSxbMCwxLCJcXFRcXF5BXzIiXSxbMSwyLCIiLDAseyJsZXZlbCI6Miwic3R5bGUiOnsiaGVhZCI6eyJuYW1lIjoibm9uZSJ9fX1dLFswLDMsIlxcVFxcXkEiLDJdLFszLDIsInNeXFwqIiwyXSxbMSwzLCJzXFxeQSIsMSx7ImxldmVsIjoyfV1d
\begin{tikzcd}
{\X\^A} & {\X\^{\T_2A}} \\
{\X\^{\T A}} & {\X\^{\T_2A}}
\arrow["{\T\^A_2}", from=1-1, to=1-2]
\arrow[Rightarrow, no head, from=1-2, to=2-2]
\arrow["{\T\^A}"', from=1-1, to=2-1]
\arrow["{s^\*}"', from=2-1, to=2-2]
\arrow["{s\^A}"{description}, Rightarrow, from=1-2, to=2-1]
\end{tikzcd}
\end{equation*}
satisfying the following property:
\begin{equation*}
% https://q.uiver.app/#q=WzAsNSxbMCwwLCJcXFRfMlxcXkFcXG8gZl5cXCoiXSxbMSwwLCJzXlxcKlxcb1xcVFxcXkFcXG8gZl5cXCoiXSxbMSwxLCJzXlxcKlxcbyhcXFQgZileXFwqXFxvXFxUXFxee0EnfSJdLFsxLDIsIihcXFRfMmYpXlxcKlxcbyBzXlxcKlxcb1xcVFxcXntBJ30iXSxbMCwyLCIoXFxUXzJmKV5cXCpcXG9cXFRfMlxcXntBJ30iXSxbMCwxLCJzXFxeQWZeXFwqIl0sWzEsMiwic15cXCpcXGthcHBhXFxeZiJdLFsyLDMsIlxcZ2FtbWEiXSxbMCw0LCJcXGthcHBhXzJcXF5mZl5cXCoiLDJdLFs0LDMsIihcXFRfMmYpXlxcKnNcXF57QSd9IiwyXV0=
\begin{tikzcd}
{\T_2\^A\o f^\*} & {s^\*\o\T\^A\o f^\*} \\
& {s^\*\o(\T f)^\*\o\T\^{A'}} \\
{(\T_2f)^\*\o\T_2\^{A'}} & {(\T_2f)^\*\o s^\*\o\T\^{A'}}
\arrow["{s\^Af^\*}", from=1-1, to=1-2]
\arrow["{\kappa_2\^ff^\*}"', from=1-1, to=3-1]
\arrow["{s^\*\kappa\^f}", from=1-2, to=2-2]
\arrow["\Comp", from=2-2, to=3-2]
\arrow["{(\T_2f)^\*s\^{A'}}"', from=3-1, to=3-2]
\end{tikzcd}
\end{equation*}

\item[indexed vertical lift] A collection of natural transformations $l\^A\colon\T\^A\Rightarrow l^\*\o\T\^{\T A}\o\T\^A$, indexed by the objects $A$ of $\X$
\begin{equation*}
% https://q.uiver.app/#q=WzAsNSxbMCwwLCJcXFhcXF5BIl0sWzEsMCwiXFxYXFxee1xcVCBBfSJdLFswLDEsIlxcWFxcXntcXFQgQX0iXSxbMCwyLCJcXFhcXF57XFxUXjJBfSJdLFsxLDIsIlxcWFxcXntcXFQgQX0iXSxbMCwxLCJcXFRcXF5BIl0sWzEsNCwiIiwwLHsibGV2ZWwiOjIsInN0eWxlIjp7ImhlYWQiOnsibmFtZSI6Im5vbmUifX19XSxbMCwyLCJcXFRcXF5BIiwyXSxbMiwzLCJcXFRcXF57XFxUIEF9IiwyXSxbMyw0LCJsXlxcKiIsMl0sWzEsMywibFxcXkEiLDEseyJsZXZlbCI6Mn1dXQ==
\begin{tikzcd}
{\X\^A} & {\X\^{\T A}} \\
{\X\^{\T A}} \\
{\X\^{\T^2A}} & {\X\^{\T A}}
\arrow["{\T\^A}", from=1-1, to=1-2]
\arrow[Rightarrow, no head, from=1-2, to=3-2]
\arrow["{\T\^A}"', from=1-1, to=2-1]
\arrow["{\T\^{\T A}}"', from=2-1, to=3-1]
\arrow["{l^\*}"', from=3-1, to=3-2]
\arrow["{l\^A}"{description}, Rightarrow, from=1-2, to=3-1]
\end{tikzcd}
\end{equation*}
satisfying the following property:
\begin{equation*}
% https://q.uiver.app/#q=WzAsNixbMCwwLCJcXFRcXF5BXFxvIGZeXFwqIl0sWzEsMCwibF5cXCpcXG9cXFRcXF57XFxUIEF9XFxvXFxUXFxeQVxcbyBmXlxcKiJdLFsxLDEsImxeXFwqXFxvXFxUXFxee1xcVCBBfVxcbyhcXFQgZileXFwqXFxvXFxUXFxee0EnfSJdLFsxLDMsIihcXFRfMmYpXlxcKlxcbyBzXlxcKlxcb1xcVFxcXntBJ30iXSxbMCwzLCIoXFxUIGYpXlxcKlxcb1xcVFxcXntBJ30iXSxbMSwyLCJsXlxcKlxcbyhcXFReMmYpXlxcKlxcb1xcVFxcXntcXFQgQSd9XFxvXFxUXFxee0EnfSJdLFswLDEsImxcXF5BZl5cXCoiXSxbMSwyLCJsXlxcKlxcVFxcXntcXFQgQX1cXGthcHBhXFxeZiJdLFswLDQsIlxca2FwcGFcXF5mIiwyXSxbNCwzLCIoXFxUIGYpXlxcKmxcXF57QSd9IiwyXSxbNSwzLCJcXGdhbW1hIl0sWzIsNSwibF5cXCpcXGthcHBhXFxee1xcVCBmfVxcVFxcXntBJ30iXV0=
\begin{tikzcd}
{\T\^A\o f^\*} & {l^\*\o\T\^{\T A}\o\T\^A\o f^\*} \\
& {l^\*\o\T\^{\T A}\o(\T f)^\*\o\T\^{A'}} \\
& {l^\*\o(\T^2f)^\*\o\T\^{\T A'}\o\T\^{A'}} \\
{(\T f)^\*\o\T\^{A'}} & {(\T f)^\*\o l^\*\o\T\^{\T A'}\o\T\^{A'}}
\arrow["{l\^Af^\*}", from=1-1, to=1-2]
\arrow["{\kappa\^f}"', from=1-1, to=4-1]
\arrow["{l^\*\T\^{\T A}\kappa\^f}", from=1-2, to=2-2]
\arrow["{l^\*\kappa\^{\T f}\T\^{A'}}", from=2-2, to=3-2]
\arrow["\Comp", from=3-2, to=4-2]
\arrow["{(\T f)^\*l\^{A'}}"', from=4-1, to=4-2]
\end{tikzcd}
\end{equation*}

\item[indexed canonical flip] A collection of natural transformations $c\^A\colon\T\^{\T A}\o\T\^A\Rightarrow c^\*\o\T\^{\T A}\o\T\^A$, indexed by the objects $A$ of $\X$
\begin{equation*}
% https://q.uiver.app/#q=WzAsNSxbMCwwLCJcXFhcXF5BIl0sWzEsMiwiXFxYXFxee1xcVF4yQX0iXSxbMSwwLCJcXFhcXF57XFxUIEF9Il0sWzAsMiwiXFxYXFxee1xcVF4yQX0iXSxbMCwxLCJcXFhcXF57XFxUIEF9Il0sWzAsMiwiXFxUXFxeQSJdLFsyLDEsIlxcVFxcXkEiXSxbMywxLCJjXlxcKiIsMl0sWzIsMywiY1xcXkEiLDEseyJsZXZlbCI6Mn1dLFswLDQsIlxcVFxcXkEiLDJdLFs0LDMsIlxcVFxcXntcXFQgQX0iLDJdXQ==
\begin{tikzcd}
{\X\^A} & {\X\^{\T A}} \\
{\X\^{\T A}} \\
{\X\^{\T^2A}} & {\X\^{\T^2A}}
\arrow["{\T\^A}", from=1-1, to=1-2]
\arrow["{\T\^A}", from=1-2, to=3-2]
\arrow["{c^\*}"', from=3-1, to=3-2]
\arrow["{c\^A}"{description}, Rightarrow, from=1-2, to=3-1]
\arrow["{\T\^A}"', from=1-1, to=2-1]
\arrow["{\T\^{\T A}}"', from=2-1, to=3-1]
\end{tikzcd}
\end{equation*}
satisfying the following property:
\begin{equation*}
% https://q.uiver.app/#q=WzAsNyxbMCwwLCJcXFRcXF57XFxUIEF9XFxvXFxUXFxeQVxcbyBmXlxcKiJdLFsxLDAsImNeXFwqXFxvXFxUXFxee1xcVCBBfVxcb1xcVFxcXkFcXG8gZl5cXCoiXSxbMSwxLCJjXlxcKlxcb1xcVFxcXntcXFQgQX1cXG8oXFxUIGYpXlxcKlxcb1xcVFxcXntBJ30iXSxbMSwzLCIoXFxUXjJmKV5cXCpcXG8gY15cXCpcXG9cXFRcXF57XFxUIEEnfVxcb1xcVFxcXntBJ30iXSxbMCwyLCJcXFRcXF57XFxUIEF9XFxvKFxcVCBmKV5cXCpcXG9cXFRcXF57QSd9Il0sWzEsMiwiY15cXCpcXG8oXFxUXjJmKV5cXCpcXG9cXFRcXF57XFxUIEEnfVxcb1xcVFxcXntBJ30iXSxbMCwzLCIoXFxUXjJmKV5cXCpcXG9cXFRcXF57XFxUIEEnfVxcb1xcVFxcXntBJ30iXSxbMCwxLCJjXFxeQWZeXFwqIl0sWzEsMiwiY15cXCpcXFRcXF57XFxUIEF9XFxrYXBwYVxcXmYiXSxbMCw0LCJcXFRcXF57XFxUIEF9XFxrYXBwYVxcXmYiLDJdLFs1LDMsIlxcZ2FtbWEiXSxbMiw1LCJjXlxcKlxca2FwcGFcXF57XFxUIGZ9XFxUXFxee0EnfSJdLFs2LDMsIihcXFReMmYpXlxcKmNcXF57QSd9IiwyXSxbNCw2LCJcXGthcHBhXFxee1xcVCBmfVxcVFxcXntBJ30iLDJdXQ==
\begin{tikzcd}
{\T\^{\T A}\o\T\^A\o f^\*} & {c^\*\o\T\^{\T A}\o\T\^A\o f^\*} \\
& {c^\*\o\T\^{\T A}\o(\T f)^\*\o\T\^{A'}} \\
{\T\^{\T A}\o(\T f)^\*\o\T\^{A'}} & {c^\*\o(\T^2f)^\*\o\T\^{\T A'}\o\T\^{A'}} \\
{(\T^2f)^\*\o\T\^{\T A'}\o\T\^{A'}} & {(\T^2f)^\*\o c^\*\o\T\^{\T A'}\o\T\^{A'}}
\arrow["{c\^Af^\*}", from=1-1, to=1-2]
\arrow["{\T\^{\T A}\kappa\^f}"', from=1-1, to=3-1]
\arrow["{c^\*\T\^{\T A}\kappa\^f}", from=1-2, to=2-2]
\arrow["{c^\*\kappa\^{\T f}\T\^{A'}}", from=2-2, to=3-2]
\arrow["{\kappa\^{\T f}\T\^{A'}}"', from=3-1, to=4-1]
\arrow["\Comp", from=3-2, to=4-2]
\arrow["{(\T^2f)^\*c\^{A'}}"', from=4-1, to=4-2]
\end{tikzcd}
\end{equation*}
\end{description}
If the indexed tangent category has negatives, then we also have:
\begin{description}
\item[indexed negation] A collection of natural transformations $n\^A\colon\T\^A\Rightarrow n^\*\o\T\^A$, indexed by the objects $A$ of $\X$
\begin{equation*}
% https://q.uiver.app/#q=WzAsNCxbMCwwLCJcXFhcXF5BIl0sWzEsMCwiXFxYXFxee1xcVCBBfSJdLFsxLDEsIlxcWFxcXntcXFQgQX0iXSxbMCwxLCJcXFhcXF57XFxUIEF9Il0sWzAsMSwiXFxUXFxeQSJdLFswLDMsIlxcVFxcXkEiLDJdLFszLDIsIm5eXFwqIiwyXSxbMSwyLCIiLDAseyJsZXZlbCI6Miwic3R5bGUiOnsiaGVhZCI6eyJuYW1lIjoibm9uZSJ9fX1dLFsxLDMsIm5cXF5BIiwxLHsibGV2ZWwiOjJ9XV0=
\begin{tikzcd}
{\X\^A} & {\X\^{\T A}} \\
{\X\^{\T A}} & {\X\^{\T A}}
\arrow["{\T\^A}", from=1-1, to=1-2]
\arrow["{\T\^A}"', from=1-1, to=2-1]
\arrow["{n^\*}"', from=2-1, to=2-2]
\arrow[Rightarrow, no head, from=1-2, to=2-2]
\arrow["{n\^A}"{description}, Rightarrow, from=1-2, to=2-1]
\end{tikzcd}
\end{equation*}
satisfying the following property:
\begin{equation*}
% https://q.uiver.app/#q=WzAsNSxbMCwwLCJcXFRcXF5BXFxvIGZeXFwqIl0sWzEsMCwibl5cXCpcXG9cXFRcXF5BXFxvIGZeXFwqIl0sWzEsMSwibl5cXCpcXG8oXFxUIGYpXlxcKlxcb1xcVFxcXntBJ30iXSxbMCwyLCIoXFxUIGYpXlxcKlxcb1xcVFxcXntBJ30iXSxbMSwyLCIoXFxUIGYpXlxcKlxcbyBuXlxcKlxcb1xcVFxcXntBJ30iXSxbMCwzLCJcXGthcHBhXFxeZiIsMl0sWzEsMiwibl5cXCpcXGthcHBhXFxeZiJdLFswLDEsIm5cXF5BZl5cXCoiXSxbMiw0LCJcXGdhbW1hIl0sWzMsNCwiKFxcVCBmKV5cXCpuXFxee0EnfSIsMl1d
\begin{tikzcd}
{\T\^A\o f^\*} & {n^\*\o\T\^A\o f^\*} \\
& {n^\*\o(\T f)^\*\o\T\^{A'}} \\
{(\T f)^\*\o\T\^{A'}} & {(\T f)^\*\o n^\*\o\T\^{A'}}
\arrow["{n\^Af^\*}", from=1-1, to=1-2]
\arrow["{\kappa\^f}"', from=1-1, to=3-1]
\arrow["{n^\*\kappa\^f}", from=1-2, to=2-2]
\arrow["\Comp", from=2-2, to=3-2]
\arrow["{(\T f)^\*n\^{A'}}"', from=3-1, to=3-2]
\end{tikzcd}
\end{equation*}
\end{description}
Notice that all these morphisms also need to satisfy the axioms of a tangent object. We can finally state the main result of this paper.

\begin{theorem}[Grothendieck construction in the context of tangent categories]
\label{theorem:full-Grothendieck}
The Grothendieck correspondence between fibrations and indexed categories lifts to an equivalence between the $2$-category $\TngFib$ of tangent fibrations and the $2$-category $\TngIndx\=\Tng(\Indx)$ of tangent indexed categories:
\begin{align*}
&\TngFib\simeq\TngIndx
\end{align*}
Similarly, the $2$-categories $\TngFib_\co$ and $\TngIndx_\co\=\Tng_\co(\Indx)$ are equivalent:
\begin{align*}
&\TngFib_\co\simeq\TngIndx_\co
\end{align*}
Moreover, these two equivalences restrict to strong and strict morphisms, yielding two other equivalences:
\begin{align*}
&\TngFib_\cong\simeq\TngIndx_\cong\=\Tng_\cong(\Indx)\\
&\TngFib_=\simeq\TngIndx_=\=\Tng_=(\Indx)
\end{align*}
\begin{proof}
To prove this result we only need to apply the $2$-functor $\Tng$ of Proposition~\ref{proposition:Tng-functiorality} to the equivalence $\Fib\simeq\Indx$ of Proposition~\ref{proposition:classic-Grothendieck}. Notice that the $2$-functors of this equivalence preserve limits, so they are $1$-morphisms in $\twoCat_\pp$.
\end{proof}
\end{theorem}

\begin{remark}
\label{remark:Grothendieck-construction-with-fixed-base}
Notice that the equivalence of Theorem~\ref{theorem:full-Grothendieck} preserves the base tangent category. Concretely this means that if $\Pi$ is a tangent fibration, then the base of the corresponding tangent indexed category $\II(\Pi)$ is precisely the base tangent category of $\Pi$ and vice versa. This can be used to restrict the equivalence to the $2$-subcategories $\TngFib(\X,\TT)\simeq\TngIndx(\X,\TT)$ of tangent fibrations with $(\X,\TT)$ for base tangent category and of tangent indexed categories with $(\X,\TT)$ for index tangent category, respectively, for any tangent category $(\X,\TT)$. Therefore, we also obtain the equivalences:
\begin{align*}
&\TngFib(\X,\TT)\simeq\TngIndx(\X,\TT)\\
&\TngFib_\co(\X,\TT)\simeq\TngIndx_\co(\X,\TT)\\
&\TngFib_\cong(\X,\TT)\simeq\TngIndx_\cong(\X,\TT)\\
&\TngFib_=(\X,\TT)\simeq\TngIndx_=(\X,\TT)
\end{align*}
Notice that, to obtain any possible tangent fibration and tangent indexed category over any possible base tangent category it was crucial to work with the $2$-categories $\Fib$ and $\Indx$, where the base categories are not fixed. Indeed, the tangent objects over the $2$-category $\Fib(\X)$ of fibrations over a fixed base category $\X$ are tangent fibrations over the trivial tangent category over $\X$, i.e., $\Tng(\Fib(\X))\simeq\TngFib(\X,\1)$, where $\1$ denotes the trivial tangent structure. Similarly, $\Tng(\Indx(\X))\simeq\TngIndx(\X,\1)$.
\end{remark}

On one hand, in Section~\ref{subsection:reduced-Grothendieck} we proved that the correspondence between indexed tangent categories and tangent fibrations forms an adjunction which, however, is not an equivalence. On the other hand, Theorem~\ref{theorem:full-Grothendieck} shows that tangent indexed categories are equivalent to tangent fibrations. This allows one to define an adjunction between indexed tangent categories and tangent indexed categories.
\par Informally, this adjunction can be regarded as the commutator between two operations: $\Indx$ which makes objects into ``\textit{indexed objects}'' and the $2$-functor $\Tng$. With this interpretation, we can see that these two operations do not commute, in contrast with other examples, like the $2$-functor $\Mnd$ which sends a $2$-category $\CC$ to the $2$-category of monads in $\CC$ and the $2$-functor $\Tng$ (see Example~\ref{example:tangent-monads-vs-tangent-objets}).
\par Furthermore, Corollary~\ref{corollary:reduced-grothendieck-base-trivial} establishes that, when the base tangent category is trivial, the adjunction between tangent fibrations and indexed tangent categories on such a base tangent category becomes an equivalence. This, together with Theorem~\ref{theorem:full-Grothendieck}, implies that, when the base tangent category is trivial, indexed tangent categories and tangent indexed categories are equivalent.
\par A similar phenomenon was observed in~\cite{moeller:monoidal-grothendieck} by Moeller and Vasilakopoulou in the context of a Grothendieck construction for monoidal categories. In particular, they discussed two distinct notions of ``monoidal indexed categories'', one that carries a global monoidal structure, and the second one whose monoidal structure is fibrewise. They showed that, when the base monoidal category is cartesian monoidal, these two notions are equivalent. This is similar to our result: when the base tangent category is trivial, the two notions of ``tangent indexed categories'' coincide.

%__________________________________________________________________________
%__________________________________________________________________________

\section{Conclusions}
\label{section:conclusions}
In this paper, we discussed two main approaches to defining a Grothendieck construction in the context of tangent categories. In Section~\ref{section:tangent-fibrations}, we recalled Cockett and Cruttwell's result which establishes that a tangent fibration can be associated with an indexed tangent category. We partially reconstructed the tangent fibration from the corresponding indexed tangent category, by providing an adjunction of $2$-functors (Theorem~\ref{theorem:reduced-Grothendieck-construction}). We noticed that this adjunction does not provide a full equivalence, except when the base tangent category is trivial (Corollary~\ref{corollary:reduced-grothendieck-base-trivial}). In Section~\ref{subsection:pseudoalgebras}, we also compared Street's notion of internal fibrations in a $2$-category with the notion of tangent fibrations (Theorem~\ref{theorem:tangent-fibrations-are-pseudoalgebras}).
\par To sidestep this issue and find a genuine equivalence between tangent fibrations and a suitable notion of indexed categories, we introduced a new concept: the notion of tangent objects. We then employed tangent objects to lift the equivalence between fibrations and indexed categories to an equivalence between tangent fibrations and tangent indexed categories, i.e., tangent objects in the $2$-category of indexed categories (Theorem~\ref{theorem:full-Grothendieck}).
\par In this paper, we also proved the following results:
\begin{description}
\item[Example~\ref{example:tangent-cats-vs-tangent-objects}] We showed that tangent categories are precisely the tangent objects in the $2$-category $\Cat$;

\item[Example~\ref{example:tangent-monads-vs-tangent-objets}] We showed that tangent monads are precisely the tangent objects in the $2$-category $\Mnd$;

\item[Corollary~\ref{corollary:reduced-grothendieck-base-trivial}] We showed that when the base tangent category is trivial the reduced Grothendieck construction becomes an equivalence;

\item[Proposition~\ref{proposition:tangent-fibrations-are-tangent-objects}] We showed that tangent fibrations are precisely the tangent objects in the $2$-category $\Fib$.
\end{description}
We also introduced these new concepts:
\begin{description}
\item[Definition~\ref{definition:indexed-tangent-category}] We introduced the notion of an indexed tangent category;

\item[Definition~\ref{definition:tangent-object}] We introduced the concept of a tangent object in a strict $2$-category;

\item[Example~\ref{example:tangent-monoidal-categories}] We introduced the concept of a tangent monoidal category;

\item[Example~\ref{example:tangent-enriched-categories}] We introduced the concept of a tangent enriched category;

\item[Definition~\ref{definition:tangent-indexed-category}] We defined the concept of a tangent indexed category.
\end{description}

%__________________________________________________________________________
\subsection{Future work}
\label{subsection:future-work}
In this paper, we introduced a new concept: the notion of tangent objects. We have already shown that tangent categories, tangent monads, and tangent fibrations are all examples of tangent objects in suitable $2$-categories and we employed tangent objects to introduce the new notions of tangent monoidal categories, tangent enriched categories, and tangent indexed categories.
\par Developing a formal theory of tangent objects will have important benefits in understanding all these notions. In particular, we are interested in extending notions of tangent category theory, such as vector fields, differential object, differential bundles, or connections, to the formal context of tangent objects.
\par We also expect the Grothendieck construction for tangent categories to play an important role in tangent category theory. In particular, in~\cite{cruttwell:reverse-tangent-cats}, Cruttwell and Lemay introduced the notion of a reverse tangent category and they employed the notion of a dual fibration. We are interested in understanding what kind of structure is carried by the dual of a tangent fibration.

%__________________________________________________________________________
%__________________________________________________________________________

\printbibliography

\end{document}